\documentclass[preprint,3p,12pt]{elsarticle}

\usepackage{booktabs}
\usepackage{color}
\usepackage{amsmath}
\usepackage{amsfonts}
\usepackage{amssymb}
\usepackage{latexsym}
\usepackage{array}
\usepackage{amsthm}
\usepackage{graphics}
\usepackage{comment}
\usepackage{epsf,graphicx,subfigure}
\usepackage{epsfig}
\usepackage{epstopdf}
\usepackage{stmaryrd}
\usepackage{tikz}
\usepackage{hyperref}

\usetikzlibrary{mindmap}

\usetikzlibrary{arrows,positioning,shapes.geometric}

\usepackage{longtable}
\usepackage{rotating}
\usepackage{lscape}
\usepackage[normalem]{ulem}

\usepackage{verbatim}
\usetikzlibrary{positioning,shadows,arrows}

\newtheorem{Theorem}{Theorem}
\newtheorem{Proposition}{Proposition}
\newtheorem{Lemma}{Lemma}
\newtheorem{Corollary}{Corollary}
\newtheorem{Assumption}{Assumption}
\newtheorem{Remark}{Remark}
\newtheorem{Definition}{Definition}
\newtheorem{Example}{Example}

\begin{document}

\begin{frontmatter}

\title{Gaussian Approximation for the Moving Averaged Modulus Wavelet Transform and its Variants}
\author[GRLiu,GRLiu2]{Gi-Ren Liu\corref{mycorrespondingauthor}}
\cortext[mycorrespondingauthor]{Gi-Ren Liu}
\ead{girenliu@mail.ncku.edu.tw}

\author[YCSheu]{Yuan-Chung Sheu}
\author[HTWu]{Hau-Tieng Wu}

\address[GRLiu]{Department of Mathematics, National Cheng-Kung University, Tainan, Taiwan}
\address[GRLiu2]{National Center for Theoretical Sciences, National Taiwan University, Taipei, Taiwan}
\address[YCSheu]{Department of Applied Mathematics, National Yang Ming Chiao Tung University, Hsinchu, Taiwan}
\address[HTWu]{Courant Institute of Mathematical Sciences, New York University, New York, United States of America
}

\begin{abstract}
The moving average of the complex modulus of the analytic wavelet transform provides a robust time-scale representation for signals to small time shifts and deformation.
In this work, we derive the Wiener chaos expansion of this representation for stationary Gaussian processes by the Malliavin calculus and combinatorial techniques.
The expansion allows us to obtain a lower bound for the Wasserstein distance between
the time-scale representations of two long-range dependent Gaussian processes in terms of Hurst indices.
Moreover, we apply the expansion to establish an upper bound for the smooth Wasserstein distance and the Kolmogorov distance
between the distributions of a random vector derived from the time-scale representation
and its normal counterpart.
It is worth mentioning that the expansion consists of infinite Wiener chaos and
the projection coefficients converge to zero slowly as the order of the Wiener chaos increases, and
we provide a rational-decay upper bound for these distribution distances, the rate of which depends on the nonlinear transformation of the amplitude of the complex wavelet coefficients.

\end{abstract}

\end{frontmatter}

\makeatletter
\def\ps@pprintTitle{%
  \let\@oddhead\@empty
  \let\@evenhead\@empty
  \let\@oddfoot\@empty
  \let\@evenfoot\@oddfoot
}
\makeatother

%% \linenumbers

\section{Introduction}
The wavelet transform (WT) is one of the most useful tools in signal processing \cite{daubechies1992ten,mallat1999wavelet}.
It performs a scale decomposition of an input signal $X$ by convolving
it with a family of functions $\{\psi_{j}\}_{j\in \mathbb{Z}}$, which are generated by dilating a mother wavelet $\psi$
in the way $\psi_{j}(t)= 2^{-j}\psi(2^{-j}t)$, where $t\in \mathbb{R}$.
The WT of $X$, denoted by $\{W[j]X(t)\}_{t\in \mathbb{R},j\in \mathbb{Z}}$, where $W[j]X(t):= X\star \psi_j(t)$, provides both the magnitude and phase information of $X$ in the time-scale domain when the mother wavelet $\psi$ is complex-valued.
This information has been applied to analyze the heart rate variability \cite{pichot1999wavelet},
detect the seizure through the Electroencephalography (EEG) signals \cite{faust2015wavelet},
and prove the existence of intermittency in the local field potentials recorded from patients with Parkinson's disease \cite{sen2007evidence}.
%By the inverse of WT, these information can also be used to simulate stationary/nonstationary Gaussian/non-Gaussian processes, including
%wind pressures and seismic ground motions,
%matching the prescribed marginal probability distribution function of $X$ and time-scale power spectral density function %\cite{bruna2019multiscale,hong2022simulating,zhang2021maximum}.
%The simulated records can be further used to assess the reliability of structures and infrastructure \cite{balafas2018wavelet}.
The WT also serves as a feature extractor in the field of machine learning
\cite{burhan2016feature,liu2021large,liu2020diffuse,subasi2021eeg,taran2020automatic},  particularly when the invariance to small temporal shifts is desired.
Because a small temporal shift of $X$ produces a phase shift on $\{W[j]X(t)\}_{t\in \mathbb{R},j\in \mathbb{Z}}$,
usually the modulus of the wavelet coefficients  (i.e., $\{U^{A}[j]X\}_{j\in\mathbb{Z}}:=\{A(|W[j]X(t)|)\}_{t\in \mathbb{R},j\in \mathbb{Z}}$ with
$A(r)=r$) and its variants (e.g., $A(r)=r^{2}$ \cite{balestriero2017linear} and $A(r)=\ln(r)$ \cite{beltran2010estimation,lardies2004modal})
are used \cite[Theorem 1]{guth2022phase}.
To make the features stable to actions of small shifts and deformations to $X$,
practitioners consider $S_{J}^{A}[j]$, which comes from
convolving the output of $U^{A}[j]$ with a low-pass filter $\phi_{J}$  \cite{mallat2012group}:
\begin{equation}\label{def:SAJ_in_introduction}
S^{A}_{J}[j]X(t) := \left(U^{A}[j]X\right)\star \phi_{J}(t),
\end{equation}
where $J$ is an integer used to control the width of the low-pass filter $\phi_{J}$ through
$\phi_{J}(t) = 2^{-J}\phi\left(2^{-J}t\right)$
and $\phi$ is a real-valued function, usually chosen to be the father wavelet associated with
$\psi$.
In the work \cite{mallat2012group}, $\{S^{A}_{J}[j]\}_{j\in \mathbb{Z}}$
is called the first-order scattering transform, and its higher-order counterparts
$\left(U^{A}[j_{n}]\cdots U^{A}[j_{2}]U^{A}[j_{1}]X\right)\star \phi_{J}$, where $n\in \mathbb{N}$ and $j_{1},\ldots,j_{n}\in \mathbb{Z}$,
are introduced to improve the mathematical understanding of convolutional neural networks.

In the presence of noise and interference, the input signal $X$ is usually viewed as a random process.
Various properties of random processes, such as stationarity and self-similarity, are characterized by their wavelet transforms, as shown in \cite{averkamp1998some,cambanis1995continuous,li2002wavelet,masry1993wavelet}.
If $X$ is a stationary Gaussian process, its spectral density can be estimated from the sample variance of empirical wavelet coefficients
based on observations of a single path of $X$ at random times \cite{bardet2010non}.
If the covariance function of $X$ exhibits hyperbolically decaying oscillations,
its spectral density contains singularities at non-zero frequencies.
In \cite{ayache2022asymptotic}, simultaneous estimators for the singularity location and long-memory parameters were proposed using quadratic functionals of filtered transforms of $X$.
In the context of non-stationary random processes, such as those arising from time-warped stationary Gaussian processes (also denoted by $X$),
the instantaneous wavelet spectrum (i.e., $\mathbb{E}\left[U^{A}[j]X(t)\right]$ with $A(r) = r^{2}$) can be applied to infer the spectral densities of the underlying stationary Gaussian processes \cite{meynard2018spectral}.
For multifractional processes, where the self-similarity property of their probability distributions is allowed to change locally over time,
the wavelet decomposition of their deterministic kernel functions provides a way to sharply capture their path behavior \cite{ayache2018multifractional}.
In \cite{roueff2009central}, the authors considered a discrete-time analog of the transformation $W[j]X$
\begin{align}\label{def:ZiJk}
\mathbf{Z}_{j,t}= \left\{\underset{k\in \mathbb{Z}}{\sum}v_{i,j}(\gamma_{j} t-k)X(k)\right\}_{i=1,2,\ldots,d},
\end{align}
where $d$ denotes the number of different wavelets,  $i$ indexes the components of the $d$-dimensional vector $\mathbf{Z}_{j,t}$,
%a fixed positive integer corresponding to the number of WT coefficients we are interested in,
$\{\gamma_{j}\}_{j\in \mathbb{N}}$ is a divergent sequence of positive integers corresponding to the scale, e.g., $\gamma_{j}=2^{j}$,
$t\in\mathbb{Z}$ corresponding to the time,
$v_{i,j}:\mathbb{Z}\rightarrow \mathbb{R}$  corresponding to the discretized wavelet satisfies $\underset{k\in \mathbb{Z}}{\sum}v_{i,j}^{2}(k)<\infty$,
and $\{X(k)\}_{k\in \mathbb{Z}}$ is a sequence of independent and identically distributed (i.i.d.) real-valued
random variables.
When $j\rightarrow\infty$, they proved that
the sequence of $d$-dimensional vectors
\begin{align}\label{empirical_mean}
n_{j}^{-1/2}\overset{n_{j}-1}{\underset{t=0}{\sum}}\left(\mathbf{Z}_{j,t}^{2}-\mathbb{E}\left[\mathbf{Z}_{j,t}^{2}\right]\right),
\end{align}
converges to a $d$-dimensional normal random vector in distribution sense under some conditions on $v_{i,j}$ and width of the averaging $n_{j}$, where $\mathbf{Z}_{j,t}^{2}$ represents the entrywise square of $\mathbf{Z}_{j,t}$.
The i.i.d. random sequence $\{X(k)\}_{k\in \mathbb{Z}}$ in (\ref{def:ZiJk}) was extended to long-range dependent Gaussian sequences in \cite{moulines2007spectral} as well as
subordinated Gaussian sequences with long-term or short-term memory in \cite{clausel2012large,serroukh2000statistical}.
Due to the difference in the strength of dependence of the sequence $X$,
the limiting theorem for the wavelet coefficients in \cite{clausel2012large} was proved by the so-called large scale reduction principle,
while the central limit theorem for (\ref{empirical_mean}) in \cite[Theorem 1]{serroukh2000statistical}
was proved by showing that the sequence $\{\mathbf{Z}_{j,t}^{2}-\mathbb{E}\left[\mathbf{Z}_{j,t}^{2}\right]\}_{t}$
is strong mixing.
Especially when the sequence $X$ is generated from a Hermite polynomial of a stationary Gaussian sequence with long memory,
the authors of \cite{clausel2014wavelet} showed that after suitable renormalization,
(\ref{empirical_mean}) converges either to a Gaussian process or to a Rosenblatt process when $j\rightarrow\infty$
depending on the order of Hermite polynomials.
For the discrete wavelet packet decomposition, which performs the discrete WT iteratively without taking account of
any nonlinear activation functions of stationary random processes, \cite{atto2009central}
proves that the outputs of decomposition converge in distribution to white Gaussian processes when the resolution level of decomposition increases
by showing that the third and higher-order cumulants
of the outputs converge to zero.
In \cite{bruna2015intermittent},
the authors considered a transformation similar to (\ref{def:SAJ_in_introduction}) as follows
\begin{align}\label{def:U:secondorder}
2^{j/2}|X\star\psi|\star \psi_{j}(t),\ t\in \mathbb{R}.
\end{align}
When $X$ is a fractional Brownian motion, the authors of \cite{bruna2015intermittent} applied a central limit theorem for sums of locally dependent random variables to show that
the marginal distribution of the random process defined in  (\ref{def:U:secondorder})
converges to a complex normal distribution when $j\rightarrow\infty$ \cite[Lemma 3.3]{bruna2015intermittent}.
In our previous work \cite{liu2023central}, we showed that non-Gaussian random fields, indexed by both time and scale,
may arise in the limit of (\ref{def:U:secondorder})
when $X$ is a stationary Gaussian process with long-range dependence and
$\psi$ is a real-valued wavelet.

In this work, we analyze the smooth Wasserstein distance and the Kolmogorov distance between the finite-dimensional distributions
of the centralized $S_{J}^{A}[j]X$
and its Gaussian counterpart, which is a normal random vector with the same covariance structure as $S_{J}^{A}[j]X$,
in the case when $X$ is a stationary Gaussian process and
$\psi$ is an analytic wavelet.
The first contribution of this work is to provide spectral representations for $U^{A}[j]X$ and $S_{J}^{A}[j]X$
for homogeneous or logarithmic functions $A$.
The representation is a linear combination of finite or infinite orthogonal Wiener chaos, depending on function $A$.
According to our literature survey,
no previous work has established the spectral representation of the modulus, or more general transform, of WT of Gaussian processes.
It is indispensable for calculating the convolution $U^{A}[j]X\star\phi_{J}$.
The second contribution of this work is to provide a quantitative central limit theorem
for $S_{J}^{A}[j]X$.
We note that the quantitative central limit theorem for random processes expressed as a linear combination of
finite Wiener chaos has been well analyzed \cite{nourdin2009stein,nourdin2012normal,nualart2006malliavin}.
However, it remains unclear how to address the case involving infinite Wiener chaos.

To achieve these goals, we apply the Malliavin calculus and combinatorial techniques to analyze
the complex modulus of the wavelet coefficients $W[j]X$.
This part has never been considered in existing works, including our previous work \cite{liu2022asymptotic,liu2024scattering}
in which the wavelet $\psi$ is real-valued and $A(r)=r^{2}$.
It is worth mentioning that the orthogonal expansion of $S_{J}^{A}[j]X$ is an infinite sum of
Wiener chaos, particularly for the practical cases $A(r)= r$ and $A(r)=\ln(r)$.
%As pointed out in \cite[Section 6.3]{nourdin2012normal}, for the case of infinite sum of Wiener chaos,
%it is difficult to straightly apply \cite[Theorem 6.1.2]{nourdin2012normal} and the results in \cite[Section 7]{nourdin2010invariance}
%to obtain an estimate of the distribution distance due to the complexity of the expression of the upper bound in \cite[Proposition %3.7]{nourdin2009stein} (see also \cite[(6.3.2)]{nourdin2012normal}).
In order to apply the general upper bound for the Gaussian approximation error in
\cite[Proposition 3.7]{nourdin2009stein} and \cite[Theorem 6.1.2 and (6.3.2)]{nourdin2012normal}
to the case of an infinite sum of Wiener chaos, we simplify the complexity of the expression of the upper bound and provide a more concise expression.
From the more concise expression for the upper bound (Proposition \ref{prop:multidim_stein_RHS}),
we obtain upper bounds for the smooth Wasserstein distance and the Kolmogorov distance between the laws of the centralized $S_{J}^{A}[j]X$
and its Gaussian counterpart. The upper bounds depend on the window length of the moving average and the function $A$
(Theorem \ref{theorem:convergence_rate:rational} and Corollary \ref{corollary:kol}).

The rest of the paper is organized as follows.
In Section \ref{sec:preliminary}, we {summarize necessary material for WT when the wavelet is analytic and} present some preliminaries about the Wiener-It$\hat{\textup{o}}$ integrals
and the Malliavin calculus. In Section \ref{sec:mainresult}, we state our main
results, including Theorems \ref{thm:U1_chaosexpansion} and \ref{theorem:convergence_rate:rational}.
The proofs of our main results and some technical lemmas are given
in the appendix.
%Section \ref{sec:proof}.
%Table \ref{List_symbols} below contains a list of frequently used symbols and abbreviations.

\section{Preliminaries}\label{sec:preliminary}
\subsection{Wavelet transform with an analytic wavelet}
Let $\psi$ be a complex-valued function in  $L^{1}(\mathbb{R})\cap L^{2}(\mathbb{R})$ defined through two real-valued mother wavelets $\psi_{R}$ and $\psi_{I}$ as follows
\begin{equation*}
\psi(t)=\psi_{R}(t)+i\psi_{I}(t),\ t\in \mathbb{R}.
\end{equation*}
The function $\psi$ is called an analytic wavelet if $\psi_{R}$ and $\psi_{I}$
meet the Hilbert pair requirement
\begin{equation}\label{def:Hilbert_transform}
\psi_{I}(t) = \frac{1}{\pi}\ \textup{p.v.}\int_{-\infty}^{\infty}\frac{\psi_{R}(\tau)}{t-\tau}d\tau,
\end{equation}
where p.v. means the principle value.
Given a signal $X: \mathbb{R}\rightarrow\mathbb{R}$, the family of analytic wavelets
\begin{equation*}
\psi_{j}(t) = \frac{1}{2^{j}}\psi\left(\frac{t}{2^{j}}\right), \ j\in \mathbb{Z},
\end{equation*}
defines the WT of $X$ \cite{lilly2010analytic} through the convolution of $\psi_{j}$
with the input $X$:
\begin{equation*}
W[j]X(t) = X\star\psi_{j}(t) = \int_{\mathbb{R}}X(s)\psi_{j}(t-s)ds,\ j\in \mathbb{Z},
\end{equation*}
when it is well defined.
The magnitude information returned by the analytic wavelet transform
\begin{align}\label{def:U1}
U^{A}\left[j\right]X :=& A(\left|W[j]X(t)\right|),\ j\in \mathbb{Z},
\end{align}
describes the envelope of wavelet coefficients of $X$, where $A(r) = r$ for $r\in [0,\infty)$.
When $A(r)= r^{2}$ (resp. $\ln(r)$), the time-scale representation $\{U^{A}\left[j\right]X(t)\}_{t\in \mathbb{R}, j\in \mathbb{Z}}$
is the so-called scalogram (resp. logarithmic scalogram) of $X$ in the field of signal processing \cite[Figure 8]{anden2014deep}.
We note that the scalogram typically refers to the squared magnitude of the WT with continuous scales.
Although this paper deals only with discrete scales, we still refer to $\{\left|W\left[j\right]X(t)\right|^{2}\}_{t\in \mathbb{R}, j\in \mathbb{Z}}$ as the scalogram of $X$.
In the following, we consider two types of functions for $A$: $A(r)=r^{\nu}$ for some $\nu\in(0,\infty)$ or $A(r)= \ln(r)$.

In practice, in order to make the features stable to actions of small diffeomorphisms, such as deformations, to $X$,
the transformation $U^{A}[j]$ is followed by a convolution with a low-pass filter $\phi_{J}$  \cite{anden2014deep,mallat2012group}:
\begin{equation*}%\label{def:S1}
S^{A}_{J}[j]X(t) := \left(U^{A}[j]X\right)\star \phi_{J}(t) = \int_{\mathbb{R}}U^{A}[j]X(s)\phi_{J}(t-s)ds,\ J\in \mathbb{Z},
\end{equation*}
where
\begin{equation*}
\phi_{J}(t) = \frac{1}{2^{J}}\phi\left(\frac{t}{2^{J}}\right)
\end{equation*}
%and $J$ is an integer used to control the width of the low-pass filter $\phi_{J}$.
and $\phi$ is a real-valued function
with its spectrum concentrated around zero (i.e.,  a low-pass filter) \cite{mallat1999wavelet}.
We denote the Fourier transform of $\psi_{R}$ by $\widehat{\psi_{R}}$ , i.e.,
\begin{align*}
\widehat{\psi_{R}}(\lambda) = \int_{\mathbb{R}}e^{-it\lambda}\psi_{R}(t)dt.
\end{align*}
Similarly, the Fourier transform of $\psi_{I}$ is denoted as $\widehat{\psi_{I}}$.
Because $\psi_{R}$ and $\psi_{I}$ are a Hilbert pair (\ref{def:Hilbert_transform}),
\begin{align}\label{def:Hilbert_pair}
\widehat{\psi_{I}}(\lambda) =  -i\ \textup{sgn}(\lambda) \widehat{\psi_{R}}(\lambda),
\end{align}
where
\begin{align*}
\textup{sgn}(\lambda)=\left\{\begin{array}{ll}-1 &\  \textup{if}\ \lambda<0,\\
0 &\ \textup{if}\  \lambda = 0,\\
+1 &\ \textup{if}\ \lambda>0.\end{array}
\right.
\end{align*}
For each $j\in \mathbb{Z}$, the Fourier transform of the scaled wavelets
\begin{equation*}
\psi_{R,j}(t) := \frac{1}{2^{j}}\psi_{R}\left(\frac{t}{2^{j}}\right)\ \textup{and}\
\psi_{I,j}(t) := \frac{1}{2^{j}}\psi_{I}\left(\frac{t}{2^{j}}\right)
\end{equation*}
are $\widehat{\psi_{R}}(2^{j}\cdot)$ and $\widehat{\psi_{I}}(2^{j}\cdot)$, respectively.
In terms of notation, $\widehat{\psi_{R,j}}(\cdot)=\widehat{\psi_{R}}(2^{j}\cdot)$ and $\widehat{\psi_{I,j}}(\cdot)=\widehat{\psi_{I}}(2^{j}\cdot)$.
All assumptions needed in this paper about the analytic wavelet $\psi$ and the low-pass filter $\phi_{J}$
are summarized as follows.
%Note that the larger $J$, the wider the support of $\phi_{J}$.
%\begin{align}
%|\widehat{\phi}(\lambda)|^{2}=\int_{1}^{\infty}|\widehat{\psi}(s\lambda)|^{2}\frac{ds}{s}.
%\end{align}

\begin{Assumption}\label{Assumption:wavelet}
For the real part $\psi_{R}$ of the analytic wavelet  $\psi$,
we assume that $\widehat{\psi_{R}}\in L^{1}(\mathbb{R})\cap L^{2}(\mathbb{R})$
and
there exists a bounded and continuous function $C_{\psi_{R}}:\mathbb{R} \rightarrow \mathbb{C}$
with $C_{\psi_{R}}(0)\neq0$ such that
\begin{equation*}%\label{eq:Psi:large_a}
\widehat{\psi_{R}}(\lambda) = C_{\psi_{R}}(\lambda)|\lambda|^{\alpha},\ \lambda\in \mathbb{R},
\end{equation*}
for a certain $\alpha>0$.
For the averaging function $\phi$, we assume that it is real-valued and $\widehat{\phi}\in L^{1}(\mathbb{R})\cap L^{2}(\mathbb{R})$.
%Without loss of generality, we further assume that $\int_{\mathbb{R}}\phi(t)dt=1$.
\end{Assumption}

\begin{Example}
The generalized Morse wavelets \cite{daubechies1988time} satisfy Assumption \ref{Assumption:wavelet}.
These wavelets are characterized by the following form for their Fourier transform:
\begin{align}\label{morse_wavelet}
\widehat{\psi}(\lambda) = \lambda^{\alpha} e^{-\lambda^{\gamma}} 1_{[0,\infty)}(\lambda),\ \alpha>0,\ \gamma>0.
\end{align}
Because the Fourier transform of $\psi_{R}$ is related to $\widehat{\psi}$ as follows
\begin{align*}%\label{psi_real_hat}
\widehat{\psi_{R}}(\lambda) = \left\{\begin{array}{lr}
2^{-1}\widehat{\psi}(\lambda) & \textup{if}\ \lambda>0,\\
0 & \textup{if}\ \lambda=0,\\
2^{-1}\overline{\widehat{\psi}(-\lambda)} & \textup{if}\ \lambda<0,
\end{array}\right.
\end{align*}
we have
$$C_{\psi_{R}}(\lambda) = 2^{-1}e^{-|\lambda|^{\gamma}},\ \lambda\in \mathbb{R}.$$
When $\gamma=1$, the wavelets of the form (\ref{morse_wavelet}) are called Cauchy wavelets of order $\alpha$.
When $(\alpha,\gamma)=(2,2)$, $\psi$ is the complex Mexican hat wavelet.
For the averaging function $\phi$, several suitable choices exist, including
\begin{align*}
\phi(t) = e^{-t^{2}},\ \phi(t) = e^{-|t|},\
\phi(t) = (1+t^{2})^{-1},
\end{align*}
where $t\in \mathbb{R}$.
%Each of these functions satisfies the condition $\widehat{\phi}\in L^{1}(\mathbb{R})\cap L^{2}(\mathbb{R})$.
\end{Example}

\subsection{Stationary Gaussian processes}
Given the relative breadth of analytical tools for Gaussian processes \cite{major1981lecture,nourdin2012normal,nualart2006malliavin}
and their ubiquity in applications \cite{makowiec2006long}, we consider stationary Gaussian processes as a model for the
input $X$ of the analytic wavelet transform.
Let $W$ be a complex-valued Gaussian random measure on $\mathbb{R}$
satisfying
\begin{align*}%\label{ortho}
W(\Delta_{1})=\overline{W(-\Delta_{1})},\ \
\mathbb{E}[W(\Delta_{1})]=0,\
\end{align*}
and
$$\mathbb{E}\left[W(\Delta_{1})
\overline{W(\Delta_{2})}\right]=\textup{Leb}(\Delta_{1}\cap\Delta_{2})$$
for any
 $\Delta_{1},\Delta_{2}\in
\mathcal{B}(\mathbb{R})$, where Leb is the Lebesgue measure on $\mathbb{R}$ {and $\mathcal{B}(\mathbb{R})$ is the
Borel {$\sigma$-}algebra on $\mathbb{R}$}.
Let $(\Omega, \mathcal{F}_{W}, \mathbb{P})$ be a
probability space, where the $\sigma$-algebra $\mathcal{F}_{W}$ is generated by $W$.

If $X$ is a mean-square continuous and
stationary real Gaussian random process with constant mean $\mu_{X}$ and
covariance function $R_{X}$,
by the Bochner-Khinchin theorem \cite[Chapter 4]{krylov2002introduction},
there exists a unique nonnegative measure $F_{X}:\mathcal{B}(\mathbb{R})\rightarrow [0,\infty)$ such that $F_{X}(\Delta) = F_{X}(-\Delta)$ for any $\Delta\in\mathcal{B}(\mathbb{R})$
and
\begin{equation*}%\label{eq:BK}
R_{X}(t) = \int_{\mathbb{R}}e^{i\lambda t}F_{X}(d\lambda),\ t\in \mathbb{R}.
\end{equation*}
The measure $F_{X}$ is called the spectral measure of the covariance function $R_{X}$.
\begin{Assumption}\label{Assumption:spectral}
The spectral measure $F_{X}$ is absolutely continuous with respect to the Lebesgue measure with a density function
$f_{X}\in L^{1}(\mathbb{R})$. It has one of the forms:
(a)
$f_{X}(\lambda) = C_{X}(\lambda)$ or
(b)
\begin{equation*}%\label{condition_f}
f_{X}(\lambda) = \frac{C_{X}(\lambda)}{|\lambda|^{1-\beta}},
\end{equation*}
where $\beta\in(0,1)$ is the long-memory
parameter and $C_{X}:\mathbb{R}\rightarrow[0,\infty)$ is bounded and continuous.
According to the convention $0^{0}=1$,
the case (a) is henceforth interpreted as the limiting case $\beta=1$. For the case (b), we further assume that $C_{X}(0)>0$.
\end{Assumption}
%For convenience of presentation, we set $0^{0}=1$, by which
%Assumption \ref{Assumption:spectral}(a) can be viewed as a limiting case \beta=1 of Assumption \ref{Assumption:spectral}(b).
\begin{Example}
The stationary Gaussian process $X$ is called an Ornstein-Uhlenbeck process \cite{krylov2002introduction} if
\begin{align*}
\textup{Cov}\left(X(t),X(t')\right) =  e^{-c|t-t'|},\ t,t'\in \mathbb{R},
\end{align*}
where $c$ is a positive constant. In this case,
\begin{align*}
f_{X}(\lambda) =  \frac{c}{\pi} \frac{1}{\lambda^{2}+c^{2}},\ \lambda\in \mathbb{R}.
\end{align*}
This spectral density satisfies Assumption \ref{Assumption:spectral}(a).
Methods for estimating the spectral density $f_{X}$
from data can be found in \cite{bardet2010non} and the references therein.
\end{Example}

\begin{Example}
For the generalized Linnik
covariance function defined as
\begin{align*}%\label{Linnik}
\textup{Cov}\left(X(t),X(t')\right)=(1+|t-t'|^{\beta_{1}})^{-\beta_{2}},\ \beta_{1}\in(0,2],\ \beta_{2}>0,
\end{align*}
\cite{lim2010analytic} shows that
the corresponding spectral density $f_{X}$ is infinitely differentiable
on $\mathbb{R}\setminus\{0\}$.
Furthermore, \cite[Corollary 3.10]{lim2010analytic} shows that the behavior of $f_{X}$ near the origin depends on the value of $\beta_{1}\beta_{2}$
as follows.
\begin{itemize}
\item For $\beta_{1}\beta_{2}>1$,
\begin{align*}
f_{X}(\lambda) \sim \left[\beta_{1}\pi\Gamma(\beta_{2})\right]^{-1} \Gamma\left(\frac{1}{\beta_{1}}\right)\Gamma\left(\frac{\beta_{1}\beta_{2}-1}{\beta_{1}}\right)
\ \textup{as}\ |\lambda|\rightarrow 0.
\end{align*}
\item For $\beta_{1}\beta_{2}<1$,
\begin{align*}
f_{X}(\lambda) \sim \left[2^{\beta_{1}\beta_{2}}\pi^{\frac{1}{2}}\Gamma\left(\frac{\beta_{1}\beta_{2}}{2}\right)\right]^{-1}\Gamma\left(\frac{1-\beta_{1}\beta_{2}}{2}\right) |\lambda|^{\beta_{1}\beta_{2}-1}\
\textup{as}\ |\lambda|\rightarrow 0.
\end{align*}
\end{itemize}
Thus, the spectral density of the generalized Linnik
covariance function satisfies Assumption \ref{Assumption:spectral}.
\end{Example}

Under Assumption \ref{Assumption:spectral}, $X$ can be expressed as a Wiener integral as follows
\begin{align*}%\label{spectral_X}
X(t) = \mu_{X}+ \int_{\mathbb{R}}e^{it\lambda}\sqrt{f_{X}(\lambda)}W(d\lambda),\ t\in \mathbb{R}.
\end{align*}
Because $\int_{\mathbb{R}}\psi_{j}(s) ds = 0$, by the stochastic Fubini theorem \cite[Theorem 2.1]{pipiras2010regularization},
\begin{align}\label{spectral_Xstarpsia}
X\star \psi_{j}(t) = \int_{\mathbb{R}}e^{it\lambda} \widehat{\psi_{R}}(2^{j}\lambda)\sqrt{f_{X}(\lambda)} W(d\lambda)
+i\int_{\mathbb{R}}e^{it\lambda} \widehat{\psi_{I}}(2^{j}\lambda)\sqrt{f_{X}(\lambda)} W(d\lambda).
\end{align}
Both the real and imaginary parts of  $X\star \psi_{j}(t)$ are normal random variables with mean zero and variance
\begin{align}\label{def:sigma_j}
\sigma_{j}^{2}
=\int_{\mathbb{R}}|\widehat{\psi_{R}}(2^{j}\lambda)|^{2}f_{X}(\lambda)d\lambda=\int_{\mathbb{R}}|\widehat{\psi_{I}}(2^{j}\lambda)|^{2}f_{X}(\lambda)d\lambda.
\end{align}

For $k\in \mathbb{N}\cup \{0\}$, denote $L_{k}\left(u\right)$ to be the Laguerre polynomial of degree $k$ with the formula
\begin{align*}
L_{k}\left(u\right) = \frac{e^{u}}{k!}\frac{d^{k}}{du^{k}}\left(e^{-u}u^{k}\right),\ k\in \mathbb{N}\cup\{0\}.
\end{align*}

\begin{Lemma}\label{lemma:Cmn}
For $x_{1},x_{2}\in \mathbb{R}$ and function $A: (0,\infty)\rightarrow \mathbb{R}$ satisfying
$$
\int_{0}^{\infty}\left(A(\sqrt{2r})\right)^{2}e^{-r}dr<\infty,
$$
we have the expansion
\begin{align}\label{hermite_expansion}
A\left(\left|x_{1}+ix_{2}\right|\right)
=\underset{\begin{subarray}{c}m,n\in \mathbb{N}\cup \{0\}\end{subarray}}{\sum}\ C_{m,n}\frac{H_{m}(x_{1})}{\sqrt{m!}}
\frac{H_{n}(x_{2})}{\sqrt{n!}},
\end{align}
where
\[
H_{m}(x)=(-1)^{m}e^{\frac{x^{2}}{2}}\frac{d^{m}}{dx^{m}}e^{-\frac{x^{2}}{2}},\ m=0,1,2,\ldots,
\]
are the (probabilistic) Hermite polynomials,
%A complete orthogonal basis of the Gaussian-Hilbert space
%is given by the Hermite polynomials $\{H_{\ell}(y)\}_{\ell=0,1,2,\ldots}$, which are
%For example, $H_{0}(y) = 1$, $H_{1}(y) = y$, and $H_{2}(y) = y^{2}-1$.
\begin{align}\label{special_Cmn}
C_{m,n} =
\left\{\begin{array}{ll}h_{m}h_{n}c_{A,(m+n)/2}&\ \textup{for}\ m,n\in 2\mathbb{N}\cup\{0\},\\
0 & \ \textup{otherwise,}
\end{array}\right.
\end{align}
\begin{align}\label{def:hmhn}
h_{m}=
(-1)^{\frac{m}{2}}\frac{\sqrt{m!}}{2^{\frac{m}{2}}\left(\frac{m}{2}\right)!},
\end{align}
and
\begin{align}\label{def:can}
c_{A,(m+n)/2} =
\int_{0}^{\infty} A(\sqrt{2u})L_{\frac{m+n}{2}}\left(u\right)e^{-u}du.
\end{align}

\end{Lemma}
The proof of Lemma \ref{lemma:Cmn} is provided in \ref{appendix_sec:Hermite_expansion_modulus}.
The constants $\{c_{A,(m+n)/2}\}_{m,n\in 2\mathbb{N}\cup\{0\}}$ in Lemma \ref{lemma:Cmn} for practical cases are shown as follows.
\begin{itemize}
\item $A(r) = r^{\nu}$,\ $\nu\in(0,\infty)$:
\begin{align}\label{def:can:r_nu}
c_{A,\frac{\ell}{2}} =
2^{\frac{\nu}{2}}\Gamma(\frac{\nu}{2}+1)
\binom{\frac{\ell}{2}-\frac{\nu}{2}-1}{\frac{\ell}{2}},\ \ell\in 2\mathbb{N}.
\end{align}
For the binomial coefficient above,
\begin{align*}
\binom{\frac{\ell}{2}-\frac{\nu}{2}-1}{\frac{\ell}{2}}
=\left(\frac{\ell}{2}-\frac{\nu}{2}-1\right) \left(\frac{\ell}{2}-\frac{\nu}{2}-2\right)\cdots\left(1-\frac{\nu}{2}\right) \left(-\frac{\nu}{2}\right)\left(\frac{\ell}{2}!\right)^{-1}.
\end{align*}
Especially, when $\nu=2$,
\begin{align*}%\label{def:cA:r2}
c_{A,0} = 2,\ c_{A,1} = -2,
\end{align*}
and $c_{A,\frac{\ell}{2}}=0$ for $\ell\in\{4,6,8,\ldots\}$.

\item $A(r) = \ln(r)$: $c_{A,0} = \frac{1}{2}\ln2-\frac{1}{2}\gamma$, where $\gamma$ is the Euler-Mascheroni constant, and
\begin{align*}%\label{def:can:log}
c_{A,\frac{\ell}{2}} = -\frac{1}{\ell}
\end{align*}
for $\ell\in 2\mathbb{N}$.
%$\ln(r) = -\gamma-\overset{\infty}{\underset{n=1}{\sum}}\frac{1}{n}L_{n}(r),$
%where $\gamma$ is the Euler-Mascheroni constant.

\end{itemize}
By (\ref{spectral_Xstarpsia}) and Lemma \ref{lemma:Cmn}, for the case $A(r)=r^{\nu}$ with $\nu>0$,
\begin{align}\notag
U^{A}[j]X(t) = &\sigma_{j}^{\nu}\underset{\begin{subarray}{c}m,n\in \mathbb{N}\cup \{0\}\end{subarray}}{\sum}\ \frac{C_{m,n}}{\sqrt{m!n!}}H_{m}\left(\frac{1}{\sigma_{j}}\int_{\mathbb{R}}e^{it\lambda} \widehat{\psi_{R}}(2^{j}\lambda)\sqrt{f_{X}(\lambda)} W(d\lambda)\right)
\\\label{Hermite_expansion_U}&\times H_{n}\left(\frac{1}{\sigma_{j}}\int_{\mathbb{R}}e^{it\lambda} \widehat{\psi_{I}}(2^{j}\lambda)\sqrt{f_{X}(\lambda)} W(d\lambda)\right).
\end{align}
For the case $A(r)=\ln(r)$,
\begin{align}\notag
U^{A}[j]X(t) = &\ln(\sigma_{j})+\underset{\begin{subarray}{c}m,n\in \mathbb{N}\cup \{0\}\end{subarray}}{\sum}\ \frac{C_{m,n}}{\sqrt{m!n!}}H_{m}\left(\frac{1}{\sigma_{j}}\int_{\mathbb{R}}e^{it\lambda} \widehat{\psi_{R}}(2^{j}\lambda)\sqrt{f_{X}(\lambda)} W(d\lambda)\right)
\\\label{Hermite_expansion_U:ln}&\times H_{n}\left(\frac{1}{\sigma_{j}}\int_{\mathbb{R}}e^{it\lambda} \widehat{\psi_{I}}(2^{j}\lambda)\sqrt{f_{X}(\lambda)} W(d\lambda)\right).
\end{align}

Let $\overline{H}=\{f\in L^{2}(\mathbb{R})\mid f(-\lambda)=\overline{f(\lambda)}\ \textup{for all}\ \lambda\in \mathbb{R}\}$ be a complex Hilbert space
with {the} inner product
$
\langle f,g\rangle = \int_{\mathbb{R}}f(\lambda)\overline{g(\lambda)}d\lambda
$
and $\|f\|_{\overline{H}}^{2}=\langle f,f\rangle$ for any $f,g\in \overline{H}$.
Given an integer $m\geq1$, we denote the $m$-th tensor product of the Hilbert space $\overline{H}$ by $\overline{H}^{\otimes m}$.
The $m$-th symmetric tensor product of $\overline{H}$
is denoted by $\overline{H}^{\odot m}$, which contains those functions $f\in \overline{H}^{\otimes m}$
satisfying
$f(\lambda_{p(1)},\ldots,\lambda_{p(m)}) = f(\lambda_{1},\ldots,\lambda_{m})$
for any permutation $(p(1),p(2),\ldots,p(m))$ of the set $\{1,2,\ldots,m\}$.
For any $f\in \overline{H}^{\otimes m}$,
the {\em $m$-fold Wiener-It$\hat{o}$ integrals} of $f$ with respect to the random measure $W$ is defined by
\begin{align*}%\label{def:multiple_wiener}
I_{m}(f) = \int_{\mathbb{R}^{m}}^{'}f(\lambda_{1},\ldots,\lambda_{m})W(d\lambda_{1})\cdots W(d\lambda_{m}),
\end{align*}
where
$\int^{'}$ means that the integral excludes the diagonal hyperplanes
$\lambda_{k}=\mp \lambda_{k^{'}}$ for $k, k^{'}\in\{1,\ldots,m\}$ and  $k\neq k^{'}$ \cite{major1981lecture}.
Lemma \ref{lemma:itoformula} below, which comes from \cite[Theorem 4.3 and Proposition 5.1]{major1981lecture} and
\cite[Theorems 2.7.7 and 2.7.10]{nourdin2012normal},
provides a significant link between nonlinear functions of normal random variables and Wiener-It$\hat{\textup{o}}$ integrals.

\begin{Lemma}[It$\hat{\textup{o}}$'s formula and product formula \cite{major1981lecture,nourdin2012normal}]\label{lemma:itoformula}
Let $f\in \overline{H}$ be such that $\|f\|_{\overline{H}}=1$. Then, for any integer $m\geq 1$, we have
\begin{align*}
H_{m}\left(\int_{\mathbb{R}}f(\lambda)W(d\lambda)\right)
=I_{m}\left(f^{\otimes m}\right).
\end{align*}
For any $m,n\geq 1$, if $f\in \overline{H}^{\odot m}$ and $g\in \overline{H}^{\odot n}$, then
\begin{align*}
I_{m}(f)I_{n}(g) = \overset{m\wedge n}{\underset{r=0}{\sum}}r!\binom{m}{r}\binom{n}{r}I_{m+n-2r}\left(f\otimes_{r} g\right),
\end{align*}
where $f\otimes_{r} g$ is
the $r$th contraction of $f$ and $g$ defined as
\begin{align}\label{def:contraction}
f\otimes_{r} g(\lambda_{1:m+n-2r})
= \int_{\mathbb{R}^{r}} f(\tau_{1:r},\lambda_{1:m-r})
g(-\tau_{1:r},\lambda_{m-r+1:m+n-2r})d\tau_{1}\cdots d\tau_{r}
\end{align}
for $r=1,2,\ldots,m\wedge n$.
{When $r=0$, set} $f\otimes_{0} g = f\otimes g$.
\end{Lemma}

Here, for any integers $p_1<p_2$, we denote $(\lambda_{p_1},\,\lambda_{p_1+1},\ldots,\lambda_{p_2})$ by $\lambda_{p_1:p_2}$ to simplify the lengthy expressions.

\begin{Lemma}\label{lemma:B_identity}
For  $\ell\in 2\mathbb{N}$, let $P[\ell]$ represent the set of permutations of $\{1,2,\ldots,\ell\}$.
For $\{\lambda_{1},\lambda_{2},\ldots,\lambda_{\ell}\}\subset \mathbb{R}\setminus \{0\}$,
define
\begin{align}\label{def:B_inLemma}
B(\ell,\lambda_{1:\ell}) = \frac{1}{\ell!} \underset{p\in P[\ell]}{\sum}\ \underset{\begin{subarray}{c}m,n\in 2\mathbb{N}\cup \{0\}\\ m+n=\ell\end{subarray}}{\sum}\
\left[\left(\frac{m}{2}\right)!\left(\frac{n}{2}\right)!\right]^{-1}
(-1)^{\frac{n}{2}}\overset{\ell}{\underset{k=\ell-n+1}{\prod}} \textup{sgn}(\lambda_{p(k)}),
\end{align}
where sgn is the sign function, indicating the sign of the non-zero input.
The equality
\begin{align}\label{exact:B}
B(\ell,\lambda_{1:\ell})
=\left\{\begin{array}{ll}
 2^\ell\   \frac{\ell}{2}!\ (\ell!)^{-1}& \ \textup{if}\ N(\lambda_{1:\ell})= \ell/2,\\
0 &\ \textup{if}\ N(\lambda_{1:\ell})\neq \ell/2,
\end{array}
\right.
\end{align}
holds, where $N(\lambda_{1:\ell})$ is the number of negative elements in $\{\lambda_{1},\lambda_{2},\ldots,\lambda_{\ell}\}$.
\end{Lemma}
By Lemma \ref{lemma:itoformula}, $U^{A}[j]X$ in
(\ref{Hermite_expansion_U}) and (\ref{Hermite_expansion_U:ln}) can be further expressed as triple summations of Wiener-It$\hat{\textup{o}}$ integrals.
However, the obtained decomposition of $U^{A}[j]X$ is not orthogonal due to the double sum over
$m$ and $n$ in (\ref{Hermite_expansion_U}) and (\ref{Hermite_expansion_U:ln}).
Because of the Hilbert pair relation (\ref{def:Hilbert_pair}) between the real and imaginary parts of the analytic wavelet,
the weighted sum of products of sign functions in (\ref{def:B_inLemma}) naturally arises.
For obtaining an orthogonal Wiener chaos decomposition of $U^{A}[j]X$ (Theorem \ref{thm:U1_chaosexpansion} below),
Lemma \ref{lemma:B_identity} plays a key role in merging non-orthogonal terms in the triple summations obtained from (\ref{Hermite_expansion_U}), (\ref{Hermite_expansion_U:ln}),
and Lemma \ref{lemma:itoformula}.
The proof of Lemma \ref{lemma:B_identity} is provided in \ref{appendix_sec:B}.

Finally, we wrap up this section by providing some properties of the Malliavin calculus
\cite[Chapter 2]{nourdin2012normal},  which are essential for the analysis of Wiener chaos in this work.
\begin{Lemma}[\cite{nourdin2012normal}]\label{lemma:deltaD} %\cite[Proposition 2.8.8]
Let $D$ denote the Malliavin derivative with respect to the Gaussian white noise random measure $W$,
and let  $L^{-1}$ denote the pseudo-inverse of the infinitesimal generator of the Ornstein-Uhlenbeck semigroup.
For every integer $p\in \mathbb{N}$ and $f\in \overline{H}^{\odot p}$,
\begin{align*}%\label{lemma:DIp}
DI_{p}(f)
= pI_{p-1}(f)
\end{align*}
and
\begin{equation}\label{def:inverseL}
L^{-1}I_{p}(f) = -\frac{1}{p}\ I_{p}\left(f\right),
\end{equation}
where $I_{p-1}(f)\in \overline{H}$ and
\begin{align*}\notag
I_{p-1}(f)(\cdot) = \int^{'}_{\mathbb{R}^{p-1}}f(\lambda_{1},\ldots,\lambda_{p-1},\cdot) W(d\lambda_{1})\cdots W(d\lambda_{p-1}).
\end{align*}
For any  $q\in[1,\infty]$,
\begin{align}\label{norm_estimate1}
\|I_{p}(f)\|_{\mathbb{D}^{1,q}}\leq c_{p,q} \|f\|_{\overline{H}^{\otimes p}},
\end{align}
where $c_{p,q}>0$ is an universal constant,
\begin{align}\label{norm_estimate2}
\|I_{p}(f)\|_{\mathbb{D}^{1,q}}:=\left(\mathbb{E}[|I_{p}(f)|^{q}]
+\mathbb{E}[\|DI_{p}(f)\|_{\overline{H}}^{q}]
\right)^{1/q},
\end{align}
and
\begin{align*}
\|f\|_{\overline{H}^{\otimes p}} = \left[\int_{\mathbb{R}^{p}}|f(\lambda_{1},\ldots,\lambda_{p})|^{2}d\lambda_{1}\cdots\lambda_{p}\right]^{1/2}.
\end{align*}
\end{Lemma}

%Due to
%the Hermite expansion of $A(|\cdot|)$ in Lemma \ref{lemma:Cmn}, the product formula in Lemma \ref{lemma:itoformula}, and
%the Hilbert pair relation (\ref{def:Hilbert_pair}) between the real and imaginary parts of the analytic wavelet,
%the weighted sum of products of sign functions in (\ref{def:B_inLemma}) will naturally pop out.
%We need the equality (\ref{exact:B}) to simplify such a complicated form. See the proof of Theorem \ref{thm:U1_chaosexpansion} for details.

\section{Main results}\label{sec:mainresult}
\begin{Theorem}\label{thm:U1_chaosexpansion}
Under Assumptions \ref{Assumption:wavelet} and \ref{Assumption:spectral},
the process $U^{A}[j]X$ defined in (\ref{def:U1})
can be expressed as a series of Wiener-It$\hat{o}$ integrals as follows.
\\
(a) For $A(r)=r^{\nu}$, where $\nu\in(0,\infty)$,
\begin{align}\label{thm:representationU}
U^{A}[j]X(t)=
\mathbb{E}\left[U^{A}[j]X(t)\right]+\sigma_{j}^{\nu}
\underset{\ell\in 2\mathbb{N}}{\sum}\
\int_{\mathbb{R}^{\ell}}^{'}Q^{(\ell)}_{t,j}(\lambda_{1},\ldots,\lambda_{\ell})W(d\lambda_{1})\cdots W(d\lambda_{\ell}),
\end{align}
where
$\mathbb{E}\left[U^{A}[j]X(t)\right] = \sigma_{j}^{\nu}c_{A,0}$,
$\sigma_{j}$ is defined in (\ref{def:sigma_j}),
\begin{align*}%\label{def:Q_in_theorem}
Q^{(\ell)}_{t,j}(\lambda_{1},\ldots,\lambda_{\ell})
=&\,c_{\ell}
\left[\sigma_{j}^{-\ell}\overset{\ell}{\underset{k=1}{\prod}}e^{it\lambda_{k}}\widehat{\psi_{R}}(2^{j}\lambda_{k})\sqrt{f_{X}(\lambda_{k})}\right]
1_{\{N(\lambda_{1:\ell})= \ell/2\}},
\end{align*}
\begin{align}\label{def:c_ell}
c_{\ell} = (-2)^{\frac{\ell}{2}}\left(\frac{\ell}{2}!\right)(\ell!)^{-1}\ c_{A,\frac{\ell}{2}},
\end{align}
$N(\lambda_{1:\ell})$ is the number of negative elements in $\{\lambda_{k}\}_{k=1}^{\ell}$,
and $c_{A,\frac{\ell}{2}}$ is defined in (\ref{def:can}).
\\
(b) The representation (\ref{thm:representationU}) can also be applied to the case $A(r)=\ln(r)$
with slight notation modification:  $\nu=0$ and $\mathbb{E}\left[U^{A}[j]X(t)\right] = c_{A,0}+\ln(\sigma_{j})$.
\end{Theorem}
The proof of Theorem \ref{thm:U1_chaosexpansion} is provided in \ref{sec:proof:thm:U1_chaosexpansion}.
Corollary \ref{prop:application} below is an easy implication of Theorem \ref{thm:U1_chaosexpansion}, whose proof is provided in \ref{sec:proof:prop:application}.
\begin{Definition}
The Wasserstein metric is defined by
\begin{align*}
d_{\textup{W}}(Z_{1},Z_{2}) = {\sup}\left\{\left|\mathbb{E}\left[h(Z_{1})\right]- \mathbb{E}\left[h(Z_{2})\right]\right|
\mid h: \mathbb{R}\rightarrow \mathbb{R}\ \textup{is Lipschitz and}\ \|h\|_{\textup{Lip}}\leq 1\ \right\}
\end{align*}
for any random variables $Z_{1}$ and $Z_{2}$.
\end{Definition}

\begin{Corollary}\label{prop:application}
Given two stationary Gaussian processes $X_{1}$ and $X_{2}$ with spectral densities
$f_{X_{1}}$ and $f_{X_{2}}$, denote
\begin{align*}%\label{def:sigma_application}
\sigma_{p,j}^{2}= \mathbb{E}\left[|W[j]X_{p}|^{2}\right],\ p=1,2.
\end{align*}
Suppose that the wavelet $\psi$ satisfies Assumption \ref{Assumption:wavelet}.
For $j\in \mathbb{Z}$, we have
\begin{align*}
d_{\textup{W}}(U^{A}[j]X_{1},U^{A}[j]X_{2})\geq
\left\{\begin{array}{ll}
2^{\frac{\nu}{2}}\Gamma(\frac{\nu}{2}+1) \left|\sigma_{1,j}^{\nu}-\sigma_{2,j}^{\nu}\right| &\ \textup{if}\ A(r) = r^{\nu}\ \textup{with}\ \nu>0;
\\\left|\ln \sigma_{1,j}-\ln \sigma_{2,j}\right| &\ \textup{if}\ A(r) = \ln(r).
\end{array}\right.
\end{align*}
Especially, if $A(r)=\ln(r)$ and
\begin{equation}\label{condition_f:v2}
f_{X_{p}}(\lambda) = \frac{C_{X_{p}}(\lambda)}{|\lambda|^{1-\beta_{p}}},\ \lambda\in \mathbb{R}\setminus\{0\},\ p=1,2,
\end{equation}
for some  long-memory
parameters $\beta_{1},\beta_{2}\in(0,1)$ and $C_{X_{1}},C_{X_{2}}\in C_{b}(\mathbb{R},[0,\infty))$, then
\begin{align*}%\label{prop:application:Hurst}
2(\ln2)^{-1}\underset{j\rightarrow\infty}{\underline{\lim}} j^{-1}d_{\textup{W}}(U^{A}[j]X_{1},U^{A}[j]X_{2})\geq  \left|\beta_{1}-\beta_{2}\right|.
\end{align*}
\end{Corollary}
In view of that $\beta_{1}$ and $\beta_{2}$ in (\ref{condition_f:v2}) give the strength of long-range dependence
of $X_{1}$ and $X_{2}$,
Corollary \ref{prop:application} shows that if two stationary Gaussian processes have remarkable difference
on the strength of long-range dependence, this discrepancy will also be reflected in the distribution distance between their logarithmic scalograms, especially at large-scales.
%The relationship, similar to (\ref{prop:application:Hurst}), can't be obtained from the homogeneous case $A(r)=r^{\nu}$, where $\nu>0$.

By Theorem \ref{thm:U1_chaosexpansion} and  the stochastic Fubini theorem \cite{pipiras2010regularization},
for the case $A(r)=r^{\nu}$, where $\nu\in(0,\infty)$, we have
\begin{align}\notag
S^{A}_{J}[j]X(t)=&
\mathbb{E}\left[S^{A}_{J}[j]X(t)\right]+\underset{\ell\in 2\mathbb{N}}{\sum}\
\int_{\mathbb{R}^{\ell}}^{'}s^{(\ell)}_{t,j}(\lambda_{1},\ldots,\lambda_{\ell})W(d\lambda_{1})\cdots W(d\lambda_{\ell})
\\\label{thm:S1_chaosexpansion}=&
\mathbb{E}\left[S^{A}_{J}[j]X(t)\right]+\underset{\ell\in 2\mathbb{N}}{\sum}\
I_{\ell}\left(s^{(\ell)}_{t,j}\right),
\end{align}
where
\begin{align}\label{integrand_s}
s^{(\ell)}_{t,j}(\lambda_{1},\ldots,\lambda_{\ell})
=&\ \sigma_{j}^{\nu}\ Q^{(\ell)}_{t,j}(\lambda_{1},\ldots,\lambda_{\ell})\ \widehat{\phi_{J}}(\lambda_{1}+\cdots+\lambda_{\ell})
\end{align}
and
\begin{align*}
\mathbb{E}\left[S^{A}_{J}[j]X(t)\right] =\int_{\mathbb{R}} \mathbb{E}\left[U^{A}[j]X(s)\right]\phi_{J}(t-s)ds = \sigma_{j}^{\nu}c_{A,0}\ \widehat{\phi}(0)
\end{align*}
for all $j,J\in \mathbb{Z}$ and $t\in \mathbb{R}$.
The representation (\ref{thm:S1_chaosexpansion}) can also be applied to the case $A(r)=\ln(r)$ with slight notation modification: $\nu=0$
and $\mathbb{E}\left[S^{A}_{J}[j]X(t)\right]=\left(c_{A,0}+\ln(\sigma_{j})\right)\widehat{\phi}(0).$

For any $d\in \mathbb{N}$, $j_{1},j_{2},\ldots,j_{d}\in \mathbb{Z}$, and $t_{1},t_{2},\ldots,t_{d}\in \mathbb{R}$,
let $\mathbf{F}=(F_{1},\ldots,F_{d})$ with
\begin{equation*}%\label{def:Fi}
F_{m} = 2^{\frac{J}{2}}\left\{S^{A}_{J}[j_{m}]X(2^{J}t_{m})-\mathbb{E}\left[S^{A}_{J}[j_{m}]X(2^{J}t_{m})\right]\right\},\ m=1,2,\ldots,d.
\end{equation*}
From (\ref{thm:S1_chaosexpansion}) and (\ref{integrand_s}), for each $K\in 2\mathbb{N}$, we can express $\mathbf{F}$ as
$$\mathbf{F} = \mathbf{F}_{\leq K} + \mathbf{F}_{> K},$$
where  $\mathbf{F}_{\leq K} = (F_{1,\leq K},\ldots,F_{d,\leq K})$
and
$\mathbf{F}_{> K} = (F_{1,> K},\ldots,F_{d,> K})$, and
the components are defined as follows
\begin{equation}\label{thm:S1_chaosexpansion_part1}
F_{m,\leq K}= 2^{\frac{J}{2}}\underset{\ell\in\{2,4,\ldots,K\}}{\sum}\ I_{\ell}\left(s^{(\ell)}_{2^{J}t_{m},j_{m}}\right)
\end{equation}
and
\begin{equation}\label{thm:S1_chaosexpansion_part2}
F_{m,>K} = 2^{\frac{J}{2}}\underset{\ell\in\{K+2,K+4,\ldots\}}{\sum}\ I_{\ell}\left(s^{(\ell)}_{2^{J}t_{m},j_{m}}\right)
\end{equation}
for $m=1,2,\ldots,d.$

Let $\mathbf{N}$ (resp. $\mathbf{N}_{\leq K}$) be a $d$-dimensional normal random vector
with the same covariance matrix as that of $\mathbf{F}$ (resp. $\mathbf{F}_{\leq K}$).
Below, we first explore the smooth Wasserstein distance \cite{gaunt2023bounding} between $\mathbf{F}$ and $\mathbf{N}$, and %. Because some of the readers might be more familiar with the Kolmogorov distance, we
conclude this section with a corollary concerning the Kolmogorov distance \cite[Definition C.2.1]{nourdin2012normal} between $\mathbf{F}$ and $\mathbf{N}$.
\begin{Definition}\label{def:sWasserstein}
The smooth Wasserstein distance between the distribution of
$\mathbb{R}^{d}$-valued random variables $\mathbf{F}$ and $\mathbf{N}$  is denoted and defined by
\begin{align}\notag
d_{\mathcal{H}_{2}}(\mathbf{F},\mathbf{N})
= \underset{h\in \mathcal{H}_{2}}{\textup{sup}}\big|\mathbb{E}\left[h\left(\mathbf{F}\right)\right]-\mathbb{E}\left[h\left(\mathbf{N}\right)\right]\big|.
\end{align}
Here, the class $\mathcal{H}_{2}$ of test functions is defined as
\begin{align*}
\mathcal{H}_{2}= \left\{h\in C^{1}(\mathbb{R}^{d})\mid  \|h\|_{\infty}\leq 1,  \|h^{'}\|_{\infty}\leq 1, h^{'}\ \textup{is Lipschitz, and}\ \|h^{'}\|_{\textup{Lip}}\leq1.\right\},
\end{align*}
where $C^{1}(\mathbb{R}^{d})$ is the space of continuously differentiable functions on $\mathbb{R}^{d}$,
$$
\|h^{'}\|_{\infty}=\underset{1\leq m\leq d}{\max}\ \underset{\mathbf{x}\in \mathbb{R}^{d}}{\sup}\ \Big|\frac{\partial h}{\partial x_{m}}(\mathbf{x})\Big|,
$$
and
\begin{align}\label{def:Lip:dh}
\|h^{'}\|_{\textup{Lip}}=\underset{1\leq m\leq d}{\max}\ \underset{\begin{subarray}{c}\mathbf{x},\mathbf{y}\in \mathbb{R}^{d}\\ \mathbf{x}\neq\mathbf{y}\end{subarray}}{\sup}\ \frac{\Big|\frac{\partial h}{\partial x_{m}}(\mathbf{x})-\frac{\partial h}{\partial x_{m}}(\mathbf{y})\Big|}{\|\mathbf{x}-\mathbf{y}\|}.
\end{align}
\end{Definition}

%For any Lipschitz differentiable function $h:\mathbb{R}^{d}\rightarrow \mathbb{R}$ with Lipschitz constant $\|h^{'}\|_{\infty}$,
%To explore the smooth Wasserstein distance between $\mathbf{F}$ and $\mathbf{N}$, we bound $d_{\mathcal{H}_{2}}(\mathbf{F},\mathbf{N})$
%the smooth Wasserstein distance between $\mathbf{F}$ and $\mathbf{N}$
By the triangle inequality
\begin{align}\notag
&\left|\mathbb{E}\left[h\left(\mathbf{F}\right)\right]-\mathbb{E}\left[h\left(\mathbf{N}\right)\right]\right|
\\\notag\leq& \left|\mathbb{E}\left[h\left(\mathbf{F}\right)\right]-\mathbb{E}\left[h\left(\mathbf{F}_{\leq K}\right)\right]\right|
\\\notag&+\left|\mathbb{E}\left[h\left(\mathbf{F}_{\leq K}\right)\right]-\mathbb{E}\left[h\left(\mathbf{N}_{\leq K}\right)\right]\right|
\\\label{triangle_inequality}&+\left|\mathbb{E}\left[h\left(\mathbf{N}_{\leq K}\right)\right]-\mathbb{E}\left[h\left(\mathbf{N}\right)\right]\right|.
\end{align}
To estimate the smooth Wasserstein distance between $\mathbf{F}$ and $\mathbf{N}$,
we estimate the first and third terms on the right-hand side of (\ref{triangle_inequality}) in Lemma \ref{lemma:gp_bound} below.
The estimate for the second term is provided in Proposition \ref{prop:multidim_stein_RHS}.

\begin{Lemma}\label{lemma:gp_bound}
For the case $A(r)=r^{\nu}$ with $\nu\in(0,\infty)\setminus 2\mathbb{N}$,
under Assumptions \ref{Assumption:wavelet} and \ref{Assumption:spectral} with $2\alpha+\beta\geq1$, for any $h\in \mathcal{H}_{2}$ and $\varepsilon>0$,
there exists a constant $C>0$ such that
\begin{equation}\label{norm_F>K}
\left|\mathbb{E}\left[h\left(\mathbf{F}\right)\right]-\mathbb{E}\left[h\left(\mathbf{F}_{\leq K}\right)\right]\right| \leq
dC K^{-\frac{\nu}{2}-\frac{1}{4}+\varepsilon}
 \end{equation}
and
\begin{equation}\label{norm_N-N}
\left|\mathbb{E}\left[h\left(\mathbf{N}_{\leq K}\right)\right]-\mathbb{E}\left[h\left(\mathbf{N}\right)\right] \right|
\leq dC K^{-\frac{\nu}{2}-\frac{1}{4}+\varepsilon}
\end{equation}
for all $J\in \mathbb{Z}$ and $K\in 2\mathbb{N}$.
The constant $C$
only depends
on $\varepsilon$, $j_{1},\ldots,j_{d}$, $f_{X}$, $\psi$, $\phi$, and $A$.
The inequalities (\ref{norm_F>K}) and (\ref{norm_N-N})
can also be applied to the case $A(r)=\ln(r)$ with slight notation modification: $\nu=0$ and $\varepsilon=0$.
\end{Lemma}

Lemma \ref{lemma:gp_bound} states that when $A(r) = r^\nu$ with $\nu\in(0,\infty)\setminus 2\mathbb{N}$ or
$A(r)=\ln(r)$,
both of which result in the expansion of $\mathbf{F}$ consisting of infinite Wiener chaos,
the first and third terms on the right-hand side of (\ref{triangle_inequality})
converge to zero when $K\rightarrow\infty$, uniformly with respect to $J$.
Moreover, the smaller $\nu$ is, the slower the convergence is.
The explicit form of the constant $C$ in Lemma \ref{lemma:gp_bound} can be found in (\ref{norm:F>K:estimate:end}).
The proof of Lemma \ref{lemma:gp_bound} is provided in \ref{sec:proof:lemma:gp_bound}.
%More details about the constant $C$ in (\ref{norm_N-N}) can be found in (\ref{norm:F>K:estimate:end}).

In the following, we apply Stein's method to analyze the second term in (\ref{triangle_inequality}).
%For readers' convenience, we recall some definition, {\color{red}like the pseudo-inverse of the infinitesimal generator of the Ornstein-Uhlenbeck %semigroup $L^{-1}$
%and the Malliavin derivative $D$,} and properties about the Malliavin calculus in \ref{appendix_sec:Stein}.
Lemma \ref{cited_thm:multiGA} below is a slight modification of \cite[Theorem 6.1.2]{nourdin2012normal},
in which the test function $h$ is twice differentiable.
How to relax the twice differentiability assumption on the test function $h$ is described in \ref{sec:proof:lemma:cited_thm:multiGA}.
%{{\color{red}Let $\mathcal{F}_{W}$ be the sigma algebra generated by the standard Brownian motion on $\mathbb{R}$. [HT: it has not been %defined yet unless I missed it...]}

\begin{Lemma}[\cite{nourdin2012normal}]\label{cited_thm:multiGA}
%Let $L^{-1}$ denote the pseudo-inverse of the infinitesimal generator of the Ornstein-Uhlenbeck semigroup,
%and let $D$ be the Malliavin derivative defined in \ref{appendix_sec:Stein}.
For any integer $d\geq2$, let $\mathbf{S}=(S_{1},\ldots,S_{d})$ be a $\mathcal{F}_{W}$-measurable random vector such that $\mathbb{E}[|S_{m}|^{4}]+\mathbb{E}[\|DS_{m}\|_{\overline{H}}^{4}]<\infty$
and $\mathbb{E}[S_{m}]=0$ for $m=1,\ldots,d$. Let $[\mathbf{C}(m,n)]_{1\leq m,n\leq d}$ be a non-negative definite matrix in $\mathbb{R}^{d\times d}$,
and let $\mathbf{N}_{\mathbf{C}}$ be a $d$-dimensional normal random vector with mean zero and covariance matrix $\mathbf{C}$.
Then, for any function $h\in\mathcal{H}_{2}$,
%$h:\mathbb{R}^{d}\rightarrow \mathbb{R}$,
{we have}
\begin{align}\label{multidim_stein}
|\mathbb{E}[h(\mathbf{S})]-\mathbb{E}[h(\mathbf{N}_{\mathbf{C}})]|\leq
&\ \frac{1}{2}
\|h^{'}\|_{\textup{Lip}}\rho,
\end{align}
where $\|h^{'}\|_{\textup{Lip}}$ is defined in (\ref{def:Lip:dh}),
%\begin{align*}%\label{def:Lip:dh}
%\|h^{'}\|_{\textup{Lip}}=\underset{1\leq m\leq d}{\max}\ \underset{\begin{subarray}{c}\mathbf{x},\mathbf{y}\in %\mathbb{R}^{d}\\ \mathbf{x}\neq\mathbf{y}\end{subarray}}{\sup}\ \frac{\Big|\frac{\partial h}{\partial %x_{m}}(\mathbf{x})-\frac{\partial h}{\partial x_{m}}(\mathbf{y})\Big|}{\|\mathbf{x}-\mathbf{y}\|}
%\end{align*}
and
\begin{align*}
\rho = \sqrt{\overset{d}{\underset{m,n=1}{\sum}}
\mathbb{E}\left[\left(\mathbf{C}(m,n)-\Big\langle DS_{n},-DL^{-1}S_{m}\Big\rangle\right)^{2}\right]}.
\end{align*}

\end{Lemma}

If $S_{m} = F_{m,\leq K}$, by Meyer's inequality (see (\ref{norm_estimate1}) and (\ref{norm_estimate2})), the condition $\mathbb{E}[|S_{m}|^{4}]+\mathbb{E}[\|DS_{m}\|_{\overline{H}}^{4}]<\infty$ holds.
On the other hand, because $F_{m,\leq K}$ is a linear combination of Wiener chaos of order greater than or equal to 2,
$\mathbb{E}[F_{m,\leq K}]=0$.
Denote the covariance matrix of $\mathbf{F}_{\leq K}$ by $\mathbf{C}_{J,K}\in \mathbb{R}^{d\times d}$; that is,
\begin{equation*}%\label{def:Fi_cov}
\mathbf{C}_{J,K}(m,n) := \mathbb{E}[F_{m,\leq K}F_{n,\leq K}],\ 1\leq m,n\leq d,
\end{equation*}
which satisfies the requirement of non-negative definite for the matrix $\mathbf{C}$ in Lemma \ref{cited_thm:multiGA}.
Because all conditions in Lemma \ref{cited_thm:multiGA} are satisfied,
the inequality (\ref{multidim_stein}) holds with $\mathbf{S} = \mathbf{F}_{\leq K}$ and $\mathbf{C}=\mathbf{C}_{J,K}$.
By making use of the explicit structure of the Wiener chaos decomposition of $\mathbf{F}_{\leq K}$ derived from Theorem \ref{thm:U1_chaosexpansion} (see also (\ref{thm:S1_chaosexpansion}) and (\ref{thm:S1_chaosexpansion_part1})),
we get the following results.

\begin{Proposition}\label{prop:multidim_stein_RHS}
Under Assumptions \ref{Assumption:wavelet} and \ref{Assumption:spectral} with $2\alpha+\beta\geq1$, for any $t_{1},...,t_{d}\in \mathbb{R}$, $j_{1},...,j_{d}\in \mathbb{Z}$,
and $h\in \mathcal{H}_{2}$,
there exists a constant $C$ such that
\begin{align}\label{Prop:statement:stein:nu}
|\mathbb{E}[h(\mathbf{F}_{\leq K})]-\mathbb{E}[h(\mathbf{N}_{\mathbf{C}_{J,K}})]|\leq  C 2^{-\frac{J}{2}}\left(\underset{\ell\in\{2,4,...,K\}}{\sum}|c_\ell| \sqrt{\ell!}3^{\frac{\ell}{2}}\right)^{2}
\end{align}
for all $J\in \mathbb{Z}$ and $K\in 2\mathbb{N}$, where $c_{\ell}$ is defined in (\ref{def:c_ell}).
The constant $C$ in (\ref{Prop:statement:stein:nu})
only depends on
$j_{1},\ldots,j_{d}$, $f_{X}$, $\psi$, $\phi$, and $A$.
\end{Proposition}
The proof of Proposition \ref{prop:multidim_stein_RHS} is provided in \ref{sec:proof:prop:multidim_stein_RHS}.
The explicit form of the constant $C$ in (\ref{Prop:statement:stein:nu}) can be found in (\ref{proof:part(b):C}).

\begin{Theorem}\label{theorem:convergence_rate:rational}
If Assumptions \ref{Assumption:wavelet} and \ref{Assumption:spectral} hold and $2\alpha+\beta\geq1$,
then for any $d\in\mathbb{N}$, $t_{1},\ldots,t_{d}\in \mathbb{R}$, and $j_{1},\ldots,j_{d}\in \mathbb{Z}$,
the sequence of random vectors
\begin{equation}\notag
\mathbf{F}=2^{J/2}\left[S^{A}_{J}[j_{1}]X(2^{J}t_{1})-\mathbb{E}\left[S^{A}_{J}[j_{1}]X(2^{J}t_{1})\right],\ldots,
S^{A}_{J}[j_{d}]X(2^{J}t_{d})-\mathbb{E}\left[S^{A}_{J}[j_{d}]X(2^{J}t_{d})\right]\right],
\end{equation}
satisfies
\begin{align}\label{sWdistance:F_N_O}
d_{\mathcal{H}_{2}}\left(\mathbf{F},\mathbf{N}\right) \leq
\left\{\begin{array}{ll}
\mathcal{O}(2^{-\frac{J}{2}}) &\  \textup{for}\ A(r)=r^{\nu}\ \textup{with}\ \nu\in2\mathbb{N},\\
\mathcal{O}(J^{-\frac{\nu}{2}-\frac{1}{4}+\varepsilon})&\  \textup{for}\ A(r)=r^{\nu}\ \textup{with}\ \nu\in(0,\infty)\setminus2\mathbb{N},\\
\mathcal{O}(J^{-\frac{1}{4}})&\  \textup{for}\ A(r)=\ln(r)
\end{array}\right.
\end{align}
for any $\varepsilon>0$ when $J\rightarrow\infty$, where
$\mathbf{N}$ is a zero-mean $d$-dimensional normal random vector with the same covariance matrix as $\mathbf{F}$.
Furthermore,
\begin{equation}\label{limit_cov}
\underset{J\rightarrow\infty}{\lim}\mathbb{E}\left[\mathbf{F}^{\textup{T}}\mathbf{F}\right] =
\left[\kappa_{m,n} \int_{\mathbb{R}}
e^{i\lambda (t_{m}-t_{n})}
|\widehat{\phi}(\lambda)|^{2}d\lambda\right]_{1\leq m,n\leq d},
\end{equation}
where
\begin{align*}%\label{thm1:kappa_j1fixed}
\kappa_{m,n}= \frac{1}{2\pi}\int_{\mathbb{R}}\textup{Cov}\big(U^{A}[j_{m}]X(\tau),U^{A}[j_{n}]X(0)\big)d\tau.
\end{align*}

\end{Theorem}
The proof of Theorem \ref{theorem:convergence_rate:rational}
is provided in \ref{sec:proof:theorem:convergence_rate:rational}.

\begin{Remark}
Before explaining why we obtain different convergence rates in Theorem 2 for  different types of function $A$, let us first recall that $\mathbf{F}$ is a $d$-dimensional random vector constructed from $S^{A}_{J}[j_{1}]X$,...,$S^{A}_{J}[j_{d}]X$
through appropriate centering and rescaling. For any $K\in 2\mathbb{N}$, $\mathbf{F}_{\leq K}$ is the truncation of the Wiener chaos expansion of $\mathbf{F}$ up to the $K$th Wiener chaos.
In Proposition 1, we estimate the smooth Wasserstein distance between $\mathbf{F}_{\leq K}$ and its limiting distribution
without specifying the form of $A$.
From the inequality in Proposition 1, it is evident that the distance between $\mathbf{F}$ and its limiting distribution decays at a rate of  $\mathcal{O}(2^{-J/2})$ as $J\rightarrow\infty$,
provided that the series on the right-hand side of (\ref{Prop:statement:stein:nu}),  given by
\begin{equation}\label{series}
\underset{\ell\in\{2,4,...,K\}}{\sum}|c_\ell| \sqrt{\ell!}3^{\frac{\ell}{2}},
\end{equation}
 converges when the truncation index  $K$ tends to infinity. Hence, the key point is the decay rate of the coefficients $c_{\ell}$,
 which is related to the Laguerre polynomial expansion of the function $A$.
\begin{itemize}
\item For $A(r)=r^{\nu}$ with $\nu\in 2\mathbb{N}$, the convergence condition on the series (\ref{series}) is satisfied because the coefficients $c_{\ell}$ in the series vanish for $\ell>\nu/2$. This leads to an exponential decay in the Wasserstein distance between $\mathbf{F}$ and its limiting distribution.

\item For $A(r)=\ln(r)$ or $A(r)=r^{\nu}$ with $\nu\in (0,\infty)\setminus 2\mathbb{N}$,
the decay of $c_{\ell}$ is not fast enough. It results that the series (\ref{series}) does not converge.
Hence, the truncation index $K$ in Proposition 1 needs to be finite.
Consequently, the smooth Wasserstein distance between $\mathbf{F}$ and its limiting distribution is mainly contributed from the truncation error,
which does not exhibit exponential decay as $J\rightarrow\infty$, as estimated in Lemma \ref{lemma:gp_bound}.
\end{itemize}
\end{Remark}

Recall that the Kolmogorov distance between the distribution of
$\mathbb{R}^{d}$-valued random variables $\mathbf{F}$ and $\mathbf{N}$  is denoted and defined by
\begin{align}\notag
d_{\textup{Kol}}(\mathbf{F},\mathbf{N})
= \underset{z_{1},\ldots,z_{d}\in \mathbb{R}}{\textup{sup}}\big|&\mathbb{P}\left(\mathbf{F}\in(-\infty,z_{1}]\times\ldots\times(-\infty,z_{d}]\right)
\\\notag&-\mathbb{P}\left(\mathbf{N}\in(-\infty,z_{1}]\times\ldots\times(-\infty,z_{d}]\right)\big|.
\end{align}
By Proposition 2.6 in \cite{gaunt2023bounding},
Corollary \ref{corollary:kol} below
shows that the Kolmogorov distance bounds between $\mathbf{F}$ and its normal counterpart $\mathbf{N}$
can be obtained from the smooth Wasserstein distance bounds in Theorem \ref{theorem:convergence_rate:rational}.
To keep the paper self-contained, the proof of  Corollary \ref{corollary:kol} is provided in \ref{sec:proof:kol}.

\begin{Corollary}\label{corollary:kol}
Let the assumptions and notation of Theorem \ref{theorem:convergence_rate:rational} prevail.
The $J$-dependent random vectors $\mathbf{F}$
satisfies
\begin{align*}%\label{triangle_inequality_v3s}
d_{\textup{Kol}}\left(\mathbf{F},\mathbf{N}\right) \leq
\left\{\begin{array}{ll}
\mathcal{O}(2^{-\frac{J}{6}}) &\  \textup{for}\ A(r)=r^{\nu}\ \textup{with}\ \nu\in2\mathbb{N},\\
\mathcal{O}(J^{-\frac{\nu}{6}-\frac{1}{12}+\varepsilon})&\  \textup{for}\ A(r)=r^{\nu}\ \textup{with}\ \nu\in(0,\infty)\setminus2\mathbb{N},\\
\mathcal{O}(J^{-\frac{1}{12}})&\  \textup{for}\ A(r)=\ln(r)
\end{array}\right.
\end{align*}
for any $\varepsilon>0$ when $J\rightarrow\infty$, where
$\mathbf{N}$ is a zero-mean $d$-dimensional normal random vector with the same covariance matrix as the corresponding $\mathbf{F}$.
\end{Corollary}

%We remark that for the case $A(r)= r^{\nu}$ with $\nu\in 2\mathbb{N}$,
%because the Wiener chaos expansion of $\mathbf{F}$ is only comprised of finite terms,
%the smooth Wasserstein distance
%between $\mathbf{F}$ and its Gaussian counterpart in Theorem \ref{theorem:convergence_rate:rational}
%can be replaced by the convex distance \cite{nourdin2022multivariate}.  (This part has not been finished.)

\begin{Remark}
Under the condition $2\alpha+\beta\geq1$, the process $X$, which exhibits long-range dependence when $\beta\in (0,1)$ and weak dependence when $\beta=1$,
is transformed to a weakly dependent process $U^{A}[j]X$
through the transformation $U^{A}[j]$.
This is the primary reason we can obtain a central limit theorem from the moving average of
$U^{A}[j]X(t)-\mathbb{E}\left[U^{A}[j]X(t)\right]$
as the window length of the averaging tends to infinity.
If $\alpha$ and $\beta$ are very small, the leading term in the Wiener chaos decomposition of $U^{A}[j]X$ will be a long-range dependent non-Gaussian process. Inspired by the results in \cite{clausel2014wavelet,leonenko2012limit,liu2024convergence,major1981lecture},
one can expect non-central limit theorems for the moving average of $U^{A}[j]X$, where the limit will be determined by the leading term in the Wiener chaos decomposition of $U^{A}[j]X$.

More specifically, (\ref{complex_U1X}) in the proof of Theorem \ref{thm:U1_chaosexpansion} shows that
the leading term in the Wiener chaos decomposition of $U^{A}[j]X(t)-\mathbb{E}\left[U^{A}[j]X(t)\right]$
is the sum of two second Wiener chaos $V_{1}$ and $V_{2}$, where
\begin{align*}
V_{1}(t) = H_{2}\left(\frac{1}{\sigma_{j}}\textup{Re}\left(W[j]X(t)\right)\right)
\end{align*}
and
\begin{align*}
V_{2}(t) = H_{2}\left(\frac{1}{\sigma_{j}}\textup{Im}\left(W[j]X(t)\right)\right),
\end{align*}
multiplied by a constant that depends on $\sigma_{j}$ and the function $A$.
From the spectral representation
\begin{align*}
V_{1}(t) = \sigma_{j}^{-2}
\int_{\mathbb{R}^{2}}^{'}e^{it(\lambda_{1}+\lambda_{2})}\left[\overset{2}{\underset{k=1}{\prod}}\widehat{\psi_{R}}(2^{j}\lambda_{k})\
\sqrt{f_{X}(\lambda_{k})}\right]W(d\lambda_{1})W(d\lambda_{2}),
\end{align*}
the covariance between $V_{1}(t)$ and $V_{1}(t')$, where $t,t'\in \mathbb{R}$, can be expressed as
\begin{align}\notag
&\textup{Cov}\left(V_{1}(t),V_{1}(t')\right)
\\\notag=& \sigma_{j}^{-4}
\int_{\mathbb{R}^{2}}e^{i(t-t')(\lambda_{1}+\lambda_{2})}|\widehat{\psi_{R}}(2^{j}\lambda_{1})\widehat{\psi_{R}}(2^{j}\lambda_{2})|^{2}
f_{X}(\lambda_{1})f_{X}(\lambda_{2})d\lambda_{1}d\lambda_{2}
\\\notag=&\sigma_{j}^{-4}
\int_{\mathbb{R}}e^{i(t-t')\lambda}f_{V_{1}}(\lambda)d\lambda,
\end{align}
where $f_{V_{1}}$ is the spectral density of $V_{1}$ and
\begin{align*}
f_{V_{1}}(\lambda)  = \int_{\mathbb{R}}\left[|\widehat{\psi_{R}}(2^{j}(\lambda-\xi))|^{2}f_{X}(\lambda-\xi)\right]
\left[|\widehat{\psi_{R}}(2^{j}\xi)|^{2}
f_{X}(\xi)\right]d\xi.
\end{align*}
The condition $2\alpha+\beta\geq1$ ensures that  the spectral densities of $V_{1}$ and $V_{2}$ do not exhibit singularities at the origin.
In other words, under this condition,
the leading term in the Wiener chaos decomposition of $U^{A}[j]X(t)-\mathbb{E}\left[U^{A}[j]X(t)\right]$ is a weakly dependent process.
This is the key reason why a central limit theorem can be derived from the averaging of $U^{A}[j]X$.
In applications of the analytic wavelet transform, as discussed in \cite{holighaus2019characterization,lilly2010analytic},
$\alpha$ is typically set to be greater than or equal to one.
%However, if both $\alpha$ and $\beta$ are very small,
%$V_{1}$ and $V_{2}$ will be long-range dependent non-Gaussian processes.
%Inspired by the literature \cite{clausel2014wavelet,leonenko2012limit}, one can expect non-central limit theorems in such cases, where the limiting distribution is %determined by the leading term in the Wiener chaos decomposition of $U^{A}[j]X$.
\end{Remark}

\section{Discussion and conclusions}
We derived the Wiener chaos decomposition of the modulus of the analytic wavelet transform and its variants of stationary Gaussian processes
and proved a quantitative central limit theorem for its moving average.
Because the complex modulus performs a square root on the scalogram,
we observed that Wiener chaos decompositions of the modulus wavelet transform and the scalogram
have a significant difference. The former consists of infinite Wiener chaos, while the latter consists of only finite Wiener chaos.
Such differences affect the convergence speed of the Gaussian approximation error of their respective moving averages.

The modulus wavelet transform, which may be further transformed by a nonlinear function $A$, is a core component of the scattering transform \cite{mallat2012group},
in which the composition of the modulus wavelet transform coupled with the moving average
\begin{align}\label{def:secondST}
S^{A}_{J}[j,j+\delta]X(t) = \int_{\mathbb{R}}U^{A}[j+\delta]U^{A}[j]X(s)\phi_{J}(t-s)ds,\ j\in \mathbb{Z}, \delta\in \mathbb{N}_{+},
\end{align}
was proposed to extract more detailed features from $X$.
In order to make sure that $S^{A}_{J}[j,j+\delta]$ is a non-expansive map \cite[Proposition 2.5]{mallat2012group},
the complex modulus, which corresponds to the case $A(r)=r$, is used as the {\em activation function}.
 %It has been applied to various applications, including audio classification  \cite{anden2011multiscale}, fetal heart rate analysis %\cite{chudavcek2013scattering}, music genre classification \cite{chen2013music}, heart sound classification \cite{li2019heart}, and sleep stage %classification \cite{liu2020hospitals,liu2020diffuse}.
As we have observed in this work, $U^{A}[j]X$ in (\ref{def:secondST})
is a non-Gaussian process, which consists of infinite Wiener chaos for the non-expansive case $A(r)=r$.
Therefore, to further analyze the second-layer modulus wavelet transform of $U^{A}[j]X$, developing new techniques to handle the nonlinear interaction across different layers and analytic wavelet is necessary, and we will report the results in our future work.
%Hence, the Gaussian approximation for the outputs of the scattering transform $S_{J}^{A}[j,j+\delta]$ remains an open question for future work.
To sum up, our current work is not only interesting from the wavelet transform perspective but also paves a way toward a theoretical understanding of the scattering transform.

\section*{Acknowledgement}
%The authors would like to thank the editors for handling
%their paper and the anonymous reviewers for their valuable
%comments. Their comments have significantly improved the
%quality of their paper.
The authors acknowledge the anonymous reviewers'
constructive comments and critiques that help improve the manuscript.
This work benefited from support of the National Center for Theoretical Science (NCTS, Taiwan).
G. R. Liu's work was supported by the National Science and Technology Council
under grant number 110-2628-M-006-003-MY3 and 113-2628-M-006 -002 -MY3.
Y. C. Sheu's work was supported by the National Science and Technology Council
under grant number NSTC 112-2115-M-A49-007 and NSTC 113-2115-M-A49-009.

\phantomsection
\addcontentsline{toc}{section}{References}
\bibliographystyle{is-abbrv}
\small
\bibliography{reference20241004.bib}

\begin{thebibliography}{10}

\bibitem{anden2014deep}
J.~And{\'e}n and S.~Mallat.
\newblock Deep scattering spectrum.
\newblock {\em IEEE Trans. Signal Process.}, 62\penalty0 (16):\penalty0
  4114--4128, 2014.

\bibitem{atto2009central}
A.~M. Atto and D.~Pastor.
\newblock Central limit theorems for wavelet packet decompositions of
  stationary random processes.
\newblock {\em IEEE Trans. Signal Process.}, 58\penalty0 (2):\penalty0
  896--901, 2009.

\bibitem{averkamp1998some}
R.~Averkamp and C.~Houdre.
\newblock Some distributional properties of the continuous wavelet transform of
  random processes.
\newblock {\em IEEE Trans. Inform. Theory}, 44\penalty0 (3):\penalty0
  1111--1124, 1998.

\bibitem{ayache2018multifractional}
A.~Ayache.
\newblock {\em Multifractional stochastic fields: wavelet strategies in
  multifractional frameworks}.
\newblock World Scientific, 2018.

\bibitem{ayache2022asymptotic}
A.~Ayache, M.~Fradon, R.~Nanayakkara, and A.~Olenko.
\newblock Asymptotic normality of simultaneous estimators of cyclic long-memory
  processes.
\newblock {\em Electron. J. Stat.}, 16\penalty0 (1):\penalty0 84--115, 2022.

\bibitem{balestriero2017linear}
R.~Balestriero and H.~Glotin.
\newblock Linear time complexity deep {Fourier} scattering network and
  extension to nonlinear invariants.
\newblock {\em arXiv preprint arXiv:1707.05841}, 2017.

\bibitem{bardet2010non}
J.-M. Bardet and P.~R. Bertrand.
\newblock A non-parametric estimator of the spectral density of a
  continuous-time {Gaussian} process observed at random times.
\newblock {\em Scand. J. Stat.}, 37\penalty0 (3):\penalty0 458--476, 2010.

\bibitem{beltran2010estimation}
J.~R. Beltr{\'a}n and J.~P. De~Leon.
\newblock Estimation of the instantaneous amplitude and the instantaneous
  frequency of audio signals using complex wavelets.
\newblock {\em Signal Process.}, 90\penalty0 (12):\penalty0 3093--3109, 2010.

\bibitem{bruna2015intermittent}
J.~Bruna, S.~Mallat, E.~Bacry, and J.-F. Muzy.
\newblock Intermittent process analysis with scattering moments.
\newblock {\em Ann. Statist.}, 43\penalty0 (1):\penalty0 323--351, 2015.

\bibitem{burhan2016feature}
N.~Burhan, R.~Ghazali, et~al.
\newblock Feature extraction of surface electromyography (s{EMG}) and signal
  processing technique in wavelet transform: A review.
\newblock In {\em 2016 IEEE International Conference on Automatic Control and
  Intelligent Systems (I2CACIS)}, pages 141--146. IEEE, 2016.

\bibitem{cambanis1995continuous}
S.~Cambanis and C.~Houdr{\'e}.
\newblock On the continuous wavelet transform of second-order random processes.
\newblock {\em IEEE Trans. Inform. Theory}, 41\penalty0 (3):\penalty0 628--642,
  1995.

\bibitem{chernozhukov2017central}
V.~Chernozhukov, D.~Chetverikov, and K.~Kato.
\newblock Central limit theorems and bootstrap in high dimensions.
\newblock {\em Ann. Probab.}, 45\penalty0 (4):\penalty0 2309--2352, 2017.

\bibitem{clausel2012large}
M.~Clausel, F.~Roueff, M.~S. Taqqu, and C.~Tudor.
\newblock Large scale behavior of wavelet coefficients of non-linear
  subordinated processes with long memory.
\newblock {\em Appl. Comput. Harmon. Anal.}, 32\penalty0 (2):\penalty0
  223--241, 2012.

\bibitem{clausel2014wavelet}
M.~Clausel, F.~Roueff, M.~S. Taqqu, and C.~Tudor.
\newblock Wavelet estimation of the long memory parameter for {Hermite}
  polynomial of {Gaussian} processes.
\newblock {\em ESAIM Probab. Stat.}, 18:\penalty0 42--76, 2014.

\bibitem{daubechies1988time}
I.~Daubechies.
\newblock Time-frequency localization operators: a geometric phase space
  approach.
\newblock {\em IEEE Trans. Inform. Theory}, 34\penalty0 (4):\penalty0 605--612,
  1988.

\bibitem{daubechies1992ten}
I.~Daubechies.
\newblock {\em Ten lectures on wavelets}.
\newblock SIAM, 1992.

\bibitem{faust2015wavelet}
O.~Faust, U.~R. Acharya, H.~Adeli, and A.~Adeli.
\newblock Wavelet-based {EEG} processing for computer-aided seizure detection
  and epilepsy diagnosis.
\newblock {\em Seizure}, 26:\penalty0 56--64, 2015.

\bibitem{gaunt2023bounding}
R.~E. Gaunt and S.~Li.
\newblock Bounding {Kolmogorov} distances through {Wasserstein} and related
  integral probability metrics.
\newblock {\em J. Math. Anal. Appl.}, 522\penalty0 (1):\penalty0 126985, 2023.

\bibitem{guth2022phase}
F.~Guth, J.~Zarka, and S.~Mallat.
\newblock Phase collapse in neural networks.
\newblock In {\em International Conference on Learning Representations}, 2022.

\bibitem{holighaus2019characterization}
N.~Holighaus, G.~Koliander, Z.~Prusa, and L.~D. Abreu.
\newblock Characterization of analytic wavelet transforms and a new phaseless
  reconstruction algorithm.
\newblock {\em IEEE Trans. Signal Process.}, 67\penalty0 (15):\penalty0
  3894--3908, 2019.

\bibitem{krylov2002introduction}
N.~V. Krylov.
\newblock {\em Introduction to the theory of random processes}, volume~43.
\newblock American Mathematical Soc., 2002.

\bibitem{lardies2004modal}
J.~Lardies, M.-N. Ta, and M.~Berthillier.
\newblock Modal parameter estimation based on the wavelet transform of output
  data.
\newblock {\em Arch. Appl. Mech.}, 73\penalty0 (9):\penalty0 718--733, 2004.

\bibitem{leonenko2012limit}
N.~Leonenko.
\newblock {\em Limit theorems for random fields with singular spectrum}, volume
  465.
\newblock Springer Science \& Business Media, 2012.

\bibitem{li2002wavelet}
T.-H. Li and H.-S. Oh.
\newblock Wavelet spectrum and its characterization property for random
  processes.
\newblock {\em IEEE Trans. Inform. Theory}, 48\penalty0 (11):\penalty0
  2922--2937, 2002.

\bibitem{lilly2010analytic}
J.~M. Lilly and S.~C. Olhede.
\newblock On the analytic wavelet transform.
\newblock {\em IEEE Trans. Inform. Theory}, 56\penalty0 (8):\penalty0
  4135--4156, 2010.

\bibitem{lim2010analytic}
S.~Lim and L.~Teo.
\newblock Analytic and asymptotic properties of multivariate generalized
  {L}innik's probability densities.
\newblock {\em J. Fourier Anal. Appl.}, 16:\penalty0 715--747, 2010.

\bibitem{liu2024convergence}
G.-R. Liu.
\newblock Convergence rate analysis in limit theorems for nonlinear functionals
  of the second {Wiener} chaos.
\newblock {\em Stochastic Process. Appl.}, 178:\penalty0 104477, 2024.

\bibitem{liu2021large}
G.-R. Liu, T.-Y. Lin, H.-T. Wu, Y.-C. Sheu, C.-L. Liu, W.-T. Liu, M.-C. Yang,
  Y.-L. Ni, K.-T. Chou, C.-H. Chen, et~al.
\newblock Large-scale assessment of consistency in sleep stage scoring rules
  among multiple sleep centers using an interpretable machine learning
  algorithm.
\newblock {\em J. Clin. Sleep Med.}, 17\penalty0 (2):\penalty0 159--166, 2021.

\bibitem{liu2020diffuse}
G.-R. Liu, Y.-L. Lo, J.~Malik, Y.-C. Sheu, and H.-T. Wu.
\newblock Diffuse to fuse {EEG} spectra--intrinsic geometry of sleep dynamics
  for classification.
\newblock {\em Biomed. Signal Process. Control}, 55:\penalty0 101576, 2020.

\bibitem{liu2022asymptotic}
G.-R. Liu, Y.-C. Sheu, and H.-T. Wu.
\newblock Asymptotic analysis of higher-order scattering transform of
  {Gaussian} processes.
\newblock {\em Electron. J. Probab.}, 27:\penalty0 1--27, 2022.

\bibitem{liu2023central}
G.-R. Liu, Y.-C. Sheu, and H.-T. Wu.
\newblock Central and noncentral limit theorems arising from the scattering
  transform and its neural activation generalization.
\newblock {\em SIAM J. Math. Anal.}, 55\penalty0 (2):\penalty0 1170--1213,
  2023.

\bibitem{liu2024scattering}
G.-R. Liu, Y.-C. Sheu, and H.-T. Wu.
\newblock When scattering transform meets non-{G}aussian random processes, a
  double scaling limit result.
\newblock {\em Bernoulli}, 30\penalty0 (3):\penalty0 2346--2371, 2024.

\bibitem{major1981lecture}
P.~Major.
\newblock {\em Muliple {Wiener-It$\hat{o}$} integrals}, volume 849.
\newblock Springer, 2nd edition, 2014.

\bibitem{makowiec2006long}
D.~Makowiec, R.~Ga{\l}a, A.~Dudkowska, A.~Rynkiewicz, M.~Zwierz, et~al.
\newblock Long-range dependencies in heart rate signals-revisited.
\newblock {\em Phys. A}, 369\penalty0 (2):\penalty0 632--644, 2006.

\bibitem{mallat1999wavelet}
S.~Mallat.
\newblock {\em A wavelet tour of signal processing}.
\newblock Elsevier, 1999.

\bibitem{mallat2012group}
S.~Mallat.
\newblock Group invariant scattering.
\newblock {\em Comm. Pure Appl. Math.}, 65\penalty0 (10):\penalty0 1331--1398,
  2012.

\bibitem{masry1993wavelet}
E.~Masry.
\newblock The wavelet transform of stochastic processes with stationary
  increments and its application to fractional {Brownian} motion.
\newblock {\em IEEE Trans. Inform. Theory}, 39\penalty0 (1):\penalty0 260--264,
  1993.

\bibitem{massey2005polar}
R.~Massey and A.~Refregier.
\newblock Polar shapelets.
\newblock {\em Monthly Not. Roy. Astr. Soc.}, 363\penalty0 (1):\penalty0
  197--210, 2005.

\bibitem{meynard2018spectral}
A.~Meynard and B.~Torr{\'e}sani.
\newblock Spectral analysis for nonstationary audio.
\newblock {\em IEEE/ACM Trans. Audio Speech Lang. Process.}, 26\penalty0
  (12):\penalty0 2371--2380, 2018.

\bibitem{moulines2007spectral}
E.~Moulines, F.~Roueff, and M.~S. Taqqu.
\newblock On the spectral density of the wavelet coefficients of long-memory
  time series with application to the log-regression estimation of the memory
  parameter.
\newblock {\em J. Time Ser. Anal.}, 28\penalty0 (2):\penalty0 155--187, 2007.

\bibitem{nourdin2009stein}
I.~Nourdin and G.~Peccati.
\newblock Stein's method on {Wiener} chaos.
\newblock {\em Probab. Theory Relat. Fields}, 145\penalty0 (1-2):\penalty0
  75--118, 2009.

\bibitem{nourdin2012normal}
I.~Nourdin and G.~Peccati.
\newblock {\em Normal approximations with {Malliavin} calculus: from {Stein's}
  method to universality}.
\newblock Number 192. Cambridge {University} Press, 2012.

\bibitem{nualart2006malliavin}
D.~Nualart.
\newblock {\em The {Malliavin} calculus and related topics}, volume 1995.
\newblock Springer, 2006.

\bibitem{pichot1999wavelet}
V.~Pichot, J.-M. Gaspoz, S.~Molliex, A.~Antoniadis, T.~Busso, F.~Roche,
  F.~Costes, L.~Quintin, J.-R. Lacour, and J.-C. Barth{\'e}l{\'e}my.
\newblock Wavelet transform to quantify heart rate variability and to assess
  its instantaneous changes.
\newblock {\em J. Appl. Physiol.}, 86\penalty0 (3):\penalty0 1081--1091, 1999.

\bibitem{pipiras2010regularization}
V.~Pipiras and M.~S. Taqqu.
\newblock Regularization and integral representations of {Hermite} processes.
\newblock {\em Statist. Probab. Lett.}, 80\penalty0 (23-24):\penalty0
  2014--2023, 2010.

\bibitem{roueff2009central}
F.~Roueff and M.~S. Taqqu.
\newblock Central limit theorems for arrays of decimated linear processes.
\newblock {\em Stochastic Process. Appl.}, 119\penalty0 (9):\penalty0
  3006--3041, 2009.

\bibitem{sen2007evidence}
A.~K. Sen and J.~O. Dostrovsky.
\newblock Evidence of intermittency in the local field potentials recorded from
  patients with {Parkinson}'s disease: a wavelet-based approach.
\newblock {\em Comput. Math. Methods Med.}, 8\penalty0 (3):\penalty0 165--171,
  2007.

\bibitem{serroukh2000statistical}
A.~Serroukh, A.~T. Walden, and D.~B. Percival.
\newblock Statistical properties and uses of the wavelet variance estimator for
  the scale analysis of time series.
\newblock {\em J. Amer. Statist. Assoc.}, 95\penalty0 (449):\penalty0 184--196,
  2000.

\bibitem{spivey2019art}
M.~Z. Spivey.
\newblock {\em The art of proving binomial identities}.
\newblock CRC Press, 2019.

\bibitem{subasi2021eeg}
A.~Subasi, T.~Tuncer, S.~Dogan, D.~Tanko, and U.~Sakoglu.
\newblock {EEG}-based emotion recognition using tunable {Q} wavelet transform
  and rotation forest ensemble classifier.
\newblock {\em Biomed. Signal Process. Control}, 68:\penalty0 102648, 2021.

\bibitem{sved1984counting}
M.~Sved.
\newblock Counting and recounting: the aftermath.
\newblock {\em Math. Intelligencer}, 6\penalty0 (4):\penalty0 44--46, 1984.

\bibitem{taran2020automatic}
S.~Taran, P.~C. Sharma, and V.~Bajaj.
\newblock Automatic sleep stages classification using optimize flexible
  analytic wavelet transform.
\newblock {\em Knowl. Based Syst.}, 192:\penalty0 105367, 2020.

\end{thebibliography}

%\bibliographystyle{is-abbrv}
%\bibliography{reference.bib}

\vskip 20 pt
%%%%%%%%%%%%%%%%%%%%%%%%%%%%%%%%%%%%%%%5
{\center{\Large  \textbf{Appendices\\}}}

\hspace{-6cm}
\numberwithin{equation}{section}
\numberwithin{figure}{section}
\begin{appendix}

%\section{Proofs}\label{sec:proof}

\section{Proof of Lemma \ref{lemma:Cmn}}\label{appendix_sec:Hermite_expansion_modulus}

Because
\begin{align*}
\int_{\mathbb{R}}\int_{\mathbb{R}}A\left(\sqrt{x_{1}^{2}+x_{2}^{2}}\right)\frac{1}{2\pi}e^{-\frac{x_{1}^{2}+x_{2}^{2}}{2}}dx_{1}dx_{2}
=\int_{0}^{\infty}A(r)e^{-\frac{r^{2}}{2}}rdr=\int_{0}^{\infty}A(\sqrt{2r})e^{-r}dr<\infty
\end{align*}
and
the set of normalized  probabilist's Hermite polynomials $\{\frac{1}{\sqrt{m!}}H_{m}\}_{m\in \mathbb{N}\cup \{0\}}$ forms
an orthonormal basis for the Gaussian Hilbert space
$L^{2}\left(\mathbb{R}, \frac{1}{\sqrt{2\pi}}e^{-\frac{x^{2}}{2}}dx\right)$,
the expansion (\ref{hermite_expansion}) holds with
\begin{align}\label{def:Cmn}
C_{m,n} = \int_{\mathbb{R}}\int_{\mathbb{R}}A\left(\sqrt{y_{1}^{2}+y_{2}^{2}}\right) \frac{H_{m}(y_{1})}{\sqrt{m!}}\frac{H_{n}(y_{2})}{\sqrt{n!}}
\frac{1}{2\pi}e^{-\frac{y_{1}^{2}+y_{2}^{2}}{2}}dy_{1}
dy_{2}
\end{align}
for $m,n\in \mathbb{N}\cup\{0\}$.
Because $A(\sqrt{y_{1}^{2}+y_{2}^{2}})$ is an even function of $y_{1}$ and $y_{2}$
and $H_{m}(-y)  = (-1)^{m}H_{m}(y)$ for any $y\in \mathbb{R}$ and $m\in \mathbb{N}\cup\{0\}$,
$C_{m,n} = 0$ if $m$ or $n$ is odd.
Hence, we only need to compute $C_{m,n}$ for the cases $m,n\in 2\mathbb{N}\cup\{0\}$.
By converting to polar coordinates, (\ref{def:Cmn}) can be rewritten as
\begin{align}\label{polar_integral}
C_{m,n} =
\int_{0}^{\infty} A(r)\left[ \int_{0}^{2\pi}\frac{H_{m}(r\cos\theta)}{\sqrt{m!}}
\frac{H_{n}(r\sin\theta)}{\sqrt{n!}}d\theta\right]\frac{1}{2\pi}e^{-\frac{r^{2}}{2}}
rdr.
\end{align}
From \cite{massey2005polar}, in which the physicists' Hermite polynomials were used,
\begin{align}\label{formula_polar}
\int_{0}^{2\pi}\frac{H_{m}(r\cos\theta)}{\sqrt{m!}}
\frac{H_{n}(r\sin\theta)}{\sqrt{n!}}d\theta
=2\pi h_{m}h_{n}
L_{\frac{m+n}{2}}\left(\frac{r^{2}}{2}\right),
\end{align}
where $h_{m}$ is defined in (\ref{def:hmhn}) and
$L_{\frac{m+n}{2}}\left(z\right)$ is the Laguerre polynomial of degree $\frac{m+n}{2}$.
We obtain (\ref{special_Cmn}) by substituting (\ref{formula_polar}) into (\ref{polar_integral}) as follows
\begin{align}\notag
C_{m,n} =&\ 2\pi h_{m}h_{n}
 \int_{0}^{\infty} A(r)L_{\frac{m+n}{2}}\left(\frac{r^{2}}{2}\right)\frac{1}{2\pi}e^{-\frac{r^{2}}{2}}rdr
 \\\notag=&\ h_{m}h_{n}\int_{0}^{\infty} A(\sqrt{2u})L_{\frac{m+n}{2}}\left(u\right)e^{-u}du.
\end{align}

\section{Proof of Lemma \ref{lemma:B_identity}}\label{appendix_sec:B}
For any $n\in\{0,2,\ldots,\ell\}$, where $\ell\in 2\mathbb{N}$,
 and $p\in P[\ell]$, if $\lambda_{k}\neq0$ for $k=1,2,\ldots,\ell$, then
\begin{align}\label{indicator}
\overset{\ell}{\underset{k=\ell-n+1}{\prod}} \textup{sgn}(\lambda_{p(k)})
=1_{\{N(\lambda_{1:\ell},n,p)\in 2\mathbb{Z}\}}-1_{\{N(\lambda_{1:\ell},n,p)\notin 2\mathbb{Z}\}},
\end{align}
where $N(\lambda_{1:\ell},n,p)$ is the cardinality of the set
$\{k\mid \ell-n+1\leq k \leq \ell, \lambda_{p(k)}<0 \}$.
Let $N=N(\lambda_{1:\ell})\leq \ell$ be the number of negative elements in $\{\lambda_{k}\}_{k=1}^{\ell}$.
By (\ref{indicator}) and the probability mass function of hypergeometric random variables,
\begin{align}\notag
&\frac{1}{\ell!}\underset{p\in P[\ell]}{\sum}\ \overset{\ell}{\underset{k=\ell-n+1}{\prod}} \textup{sgn}(\lambda_{p(k)})
\\\notag=\ & \frac{1}{\ell!}\underset{p\in P[\ell]}{\sum}1_{\{N(\lambda_{1:\ell},n,p)\in 2\mathbb{Z}\}}
-\frac{1}{\ell!}\underset{p\in P[\ell]}{\sum}1_{\{N(\lambda_{1:\ell},n,p)\notin 2\mathbb{Z}\}}
=
\binom{\ell}{n}^{-1}a_{n},
\end{align}
where
\begin{equation}\label{def:an}
a_{n} = \overset{n}{\underset{q=0}{\sum}}(-1)^{q}\binom{N}{q}\binom{\ell-N}{n-q}.
\end{equation}
Hence, (\ref{def:B_inLemma}) can be rewritten as
\begin{align}\label{def:B_v2}
B(\ell,\lambda_{1:\ell}) =&\underset{n\in\{0,2,\ldots,\ell\}}{\sum}\
\left[\left(\frac{\ell}{2}-\frac{n}{2}\right)!\left(\frac{n}{2}\right)!\right]^{-1}
(-1)^{\frac{n}{2}}\binom{\ell}{n}^{-1}a_{n}
=\underset{n\in\{0,2,\ldots,\ell\}}{\sum}\
w_{n} a_{n},
\end{align}
where
\begin{align}\label{def:weight}
w_{n}= \frac{1}{\ell!}\frac{n! (\ell-n)!}{(n/2)!(\ell/2-n/2)!}
(-1)^{\frac{n}{2}}.
\end{align}

{\bf Case 1 ($N=\ell/2$):}
In this case, for $n\in\{0,2,\ldots,\ell\}$,
\begin{align}\label{hyp_case1}
a_{n}=\overset{n}{\underset{q=0}{\sum}}(-1)^{q}\binom{\ell/2}{q}\binom{\ell/2}{n-q}
= (-1)^{n/2}\binom{\ell/2}{n/2},
\end{align}
where the last equality follows from \cite[Identity 81 in Page 61]{spivey2019art}.
By substituting (\ref{hyp_case1}) into (\ref{def:B_v2}), $B(\ell,\lambda_{1:\ell})$ can be rewritten as follows
\begin{align}\notag
B(\ell,\lambda_{1:\ell})
=&\left(\frac{\ell}{2}!\right)\frac{1}{\ell!}
\underset{n\in\{0,2,\ldots,\ell\}}{\sum}
\binom{n}{n/2}\binom{\ell-n}{\ell/2-n/2}
= \left(\frac{\ell}{2}!\right)\frac{1}{\ell!}2^{\ell},
\end{align}
where the last equality follows from \cite{sved1984counting}.

{\bf Case 2 ($N>\ell/2$):}
First of all, since the constant $a_{n}$ in (\ref{def:an}) is recognized to be the
coefficient of  $x^{n}$ in the polynomial
$
P_{\ell,N}(x) = (1-x)^{N}(1+x)^{\ell-N}
$
(in terms of notation, $a_{n}= [x^{n}] P_{\ell,N}(x)$),
from (\ref{def:B_v2}), $B(\ell,\lambda_{1:\ell})$ can be viewed as a weighted sum of
the coefficients of  $\{x^{0},x^{2},\ldots,x^{\ell}\}$ in $P_{\ell,N}(x)$,
i.e.,
\begin{align}\label{B:weighted_sum_coeff1}
B(\ell,\lambda_{1:\ell}) =
\underset{n\in\{0,2,\ldots,\ell\}}{\sum}w_{n} a_{n}
=\underset{n\in\{0,2,\ldots,\ell\}}{\sum}w_{n}[x^{n}] P_{\ell,N}(x).
\end{align}
Because $P_{\ell,N}(x)$ can be rewritten as
\begin{align}\notag
P_{\ell,N}(x) =& (1-x)^{N}\left[\overset{\ell-N}{\underset{k=0}{\sum}}b_{k}(1-x)^{k}\right]
=\overset{\ell-N}{\underset{h=0}{\sum}}b_{\ell-N-h}(1-x)^{\ell-h}
\end{align}
for some constants $b_{0},b_{1},\ldots,b_{\ell-N}$, (\ref{B:weighted_sum_coeff1}) can be rewritten as
\begin{align}\notag
B(\ell,\lambda_{1:\ell}) =&
\underset{n\in\{0,2,\ldots,\ell\}}{\sum}w_{n}[x^{n}] \left[\overset{\ell-N}{\underset{h=0}{\sum}}b_{\ell-N-h}(1-x)^{\ell-h}\right]
\\\label{B:weighted_sum_coeff2}=&\overset{\ell-N}{\underset{h=0}{\sum}}b_{\ell-N-h}A_{h},
\end{align}
where
\begin{equation}\label{def:Ah1}
A_{h}:=\underset{n\in\{0,2,\ldots,\ell\}}{\sum}w_{n}[x^{n}] (1-x)^{\ell-h}.
\end{equation}
By (\ref{def:weight}), for each $h\in\{0,1,\ldots,\ell-N\}$,
\begin{align}\notag
A_{h}
=&\frac{1}{\ell!}\underset{\begin{subarray}{c}n\in\{0,2,\ldots,\ell\}\\ n\leq \ell-h\end{subarray}}{\sum} \frac{n! (\ell-n)!}{(n/2)!(\ell/2-n/2)!}
(-1)^{\frac{n}{2}}\binom{\ell-h}{n}
\\\label{def:Ah2}=&
\frac{1}{\ell!}(\ell-h)!
\underset{\begin{subarray}{c}n\in\{0,2,\ldots,\ell\}\\ n\leq \ell-h\end{subarray}}{\sum} \frac{(\ell-n)!}{(n/2)!(\ell/2-n/2)!}
(-1)^{\frac{n}{2}}\frac{1}{(\ell-n-h)!}.
\end{align}
Obviously, $A_{0} = 0$.
For $h\in\{1,\ldots,\ell-N\}\cap 2\mathbb{N}$, from (\ref{def:Ah2}),
\begin{align}\notag
A_{h}=& \frac{1}{\ell!}(\ell-h)!\ 2^{\frac{h}{2}}
\underset{\begin{subarray}{c}n\in\{0,2,\ldots,\ell\}\\ n\leq \ell-h\end{subarray}}{\sum} \frac{(\ell-n-1)(\ell-n-3)\cdots (\ell-n-h+1)}{(n/2)!(\ell/2-n/2-h/2)!}
(-1)^{\frac{n}{2}}
\\\notag=& \frac{1}{\ell!}\frac{(\ell-h)!\ 2^{\frac{h}{2}}}{(\ell/2-h/2)!}
\overset{(\ell-h)/2}{\underset{k=0}{\sum}}
\underset{h/2\ \textup{terms}}{\underbrace{(\ell-2k-1)(\ell-2k-3)\cdots (\ell-2k-h+1)}}\binom{(\ell-h)/2}{k}
(-1)^{k}.
\end{align}
Under the condition $h<\frac{\ell}{2}$, which is satisfied because  $N>\frac{\ell}{2}$, we have
$$
\overset{(\ell-h)/2}{\underset{k=0}{\sum}}
k^{q}\binom{(\ell-h)/2}{k}
(-1)^{k} = 0,\ \ q=0,1,\ldots,\frac{h}{2}.
$$
Hence, we obtain $A_{h}=0$ for $h\in\{1,\ldots,\ell-N\}\cap 2\mathbb{N}$.

For $h\in\{1,\ldots,\ell-N\}\setminus 2\mathbb{N}$, from (\ref{def:Ah2}),
\begin{align}\notag
A_{h}=&
\frac{1}{\ell!}(\ell-h)!\ 2^{\frac{h+1}{2}}
\underset{\begin{subarray}{c}n\in\{0,2,\ldots,\ell\}\\ n\leq \ell-h\end{subarray}}{\sum} \frac{(\ell-n-1)(\ell-n-3)\cdots (\ell-n-h+2)}{(n/2)!(\ell/2-n/2-h/2-1/2)!}
(-1)^{\frac{n}{2}}
\\\notag=&\frac{1}{\ell!}\frac{(\ell-h)!\ 2^{\frac{h+1}{2}}}{(\ell/2-h/2-1/2)!}
\hspace{-0.2cm}\overset{(\ell-h-1)/2}{\underset{k=0}{\sum}}\hspace{-0.3cm}
\underset{(h-1)/2\ \textup{terms}}{\underbrace{(\ell-2k-1)(\ell-2k-3)\cdots (\ell-2k-h+2)}}
\\\notag&\times\binom{(\ell-h-1)/2}{k}
(-1)^{k}.
\end{align}
Under the condition $h<\frac{\ell}{2}$, which is satisfied again because  $N>\frac{\ell}{2}$, we have
$$
\overset{(\ell-h-1)/2}{\underset{k=0}{\sum}}
k^{q}\binom{(\ell-h-1)/2}{k}
(-1)^{k} = 0,\ \ q=0,1,\ldots,\frac{h-1}{2}.
$$
Hence,
we also have $A_{h}=0$ for $h\in\{1,\ldots,\ell-N\}\setminus 2\mathbb{N}$.

In summary, we have $A_{h}=0$ for $h\in\{0,1,\ldots,\ell-N\}$. Therefore, by
(\ref{B:weighted_sum_coeff2}),
$B(\ell,\lambda_{1:\ell})=0$  if $N>\ell/2$.

{\bf Case 3 ($N<\ell/2$):}
For this case, we consider the following decomposition of $P_{\ell,N}(x)$:
\begin{align}\notag
P_{\ell,N}(x) =& (1+x)^{\ell-N}\left[\overset{N}{\underset{k=0}{\sum}}\widetilde{b}_{k}(1+x)^{k}\right]
=\overset{N}{\underset{h=0}{\sum}}\widetilde{b}_{N-h}(1+x)^{\ell-h},
\end{align}
where $\widetilde{b}_{0},\ldots,\widetilde{b}_{N}\in \mathbb{R}.$
Similar to (\ref{B:weighted_sum_coeff2}), we have
\begin{align*}%\label{B:weighted_sum_coeff3}
B(\ell,\lambda_{1:\ell}) =\overset{N}{\underset{h=0}{\sum}}\widetilde{b}_{N-h}\widetilde{A}_{h},
\end{align*}
where
$$\widetilde{A}_{h}:=\underset{n\in\{0,2,\ldots,\ell\}}{\sum}w_{n}[x^{n}] (1+x)^{\ell-h}.$$
Because $$[x^{n}] (1+x)^{\ell-h}=\binom{\ell-h}{n} = [x^{n}] (1-x)^{\ell-h}$$ for any nonnegative even integer $n$,
$\widetilde{A}_{h}=A_{h}=0$ for any nonnegative integer $h<\ell/2$, where $A_{h}$ is defined in (\ref{def:Ah1}).
Therefore, $B(\ell,\lambda_{1:\ell})=0$ for the case $N<\ell/2$.

\section{Proof of Theorem \ref{thm:U1_chaosexpansion}}\label{sec:proof:thm:U1_chaosexpansion}
We only prove the result for the case $A(r)=r^{\nu}$, where $\nu>0$. The proof for the case $A(r)=\ln(r)$ is similar, so we omit it.
First of all, we denote
\begin{align*}%\label{def:qRqI}
q_{R}(\lambda) = \widehat{\psi_{R}}(2^{j}\lambda)\sqrt{f_{X}(\lambda)}, \ q_{I}(\lambda) = \widehat{\psi_{I}}(2^{j}\lambda)\sqrt{f_{X}(\lambda)},
\end{align*}
$q_{R,t}(\lambda)= e^{it\lambda}q_{R}(\lambda)$,
and $q_{I,t}(\lambda)= e^{it\lambda}q_{I}(\lambda) $.
By It$\hat{\textup{o}}$'s formula in Lemma \ref{lemma:itoformula},
for any $m,n\in \mathbb{N}$,
\begin{align}\label{proof:HmHn}
H_{m}\left(\frac{1}{\sigma_{j}}\int_{\mathbb{R}}q_{R,t}(\lambda) W(d\lambda)\right)
H_{n}\left(\frac{1}{\sigma_{j}}\int_{\mathbb{R}}q_{I,t}(\lambda) W(d\lambda)\right)
=
\sigma_{j}^{-(m+n)}I_{m}\left(q_{R,t}^{\otimes m}\right)I_{n}\left(q_{I,t}^{\otimes n}\right),
\end{align}
where $q_{R,t}^{\otimes m}$ is the $m$-fold tensor product of $q_{R,t}$, i.e.,
\begin{equation*}
q_{R,t}^{\otimes m}(\lambda_{1},\ldots,\lambda_{m}) = \overset{m}{\underset{k=1}{\prod}} q_{R,t}(\lambda_{k}),\ \lambda_{1},\ldots,\lambda_{m}\in \mathbb{R}.
\end{equation*}
By default, $I_{0}=1$.
By substituting (\ref{proof:HmHn}) into (\ref{Hermite_expansion_U}), we obtain
\begin{align}\label{complex_U1X}
U^{A}[j]X(t) = \sigma_{j}^{\nu}\underset{\begin{subarray}{c}m,n\in \mathbb{N}\cup \{0\}\end{subarray}}{\sum}\ \frac{C_{m,n}}{\sqrt{m!n!}}\sigma_{j}^{-(m+n)}I_{m}\left(q_{R,t}^{\otimes m}\right)I_{n}\left(q_{I,t}^{\otimes n}\right).
\end{align}
By the product formula in Lemma \ref{lemma:itoformula} (see also \cite[Proposition A.1]{clausel2014wavelet}),
\begin{align}\label{product_formula_IRII}
I_{m}\left(q_{R,t}^{\otimes m}\right)I_{n}\left(q_{I,t}^{\otimes n}\right)
= \overset{m\wedge n}{\underset{r=0}{\sum}}r! \binom{m}{r}\binom{n}{r}
I_{m+n-2r}(q_{R,t}^{\otimes m}\otimes_{r}q_{I,t}^{\otimes n}),
\end{align}
where $\otimes_{\ell}$ is the contraction operator defined in Lemma \ref{lemma:itoformula}, or more precisely,
\begin{align}\notag
q_{R,t}^{\otimes m}\otimes_{r}q_{I,t}^{\otimes n}(\lambda_{1},\ldots,\lambda_{m+n-2r})
=\int_{\mathbb{R}^{r}}q_{R,t}^{\otimes m}(\lambda_{1},\ldots,\lambda_{m-r},&u_{1},u_{2},\ldots,u_{r})
\\\notag q_{I,t}^{\otimes n}(\lambda_{m-r+1},\ldots,\lambda_{m+n-2r},-u_{1},&-u_{2},\ldots,-u_{r})du_{1}\cdots du_{r}
\end{align}
for every $r\in\{1,2,\ldots,m\wedge n\}$.
For $r=0$,
\begin{align}\notag
q_{R,t}^{\otimes m}\otimes_{0}q_{I,t}^{\otimes n}(\lambda_{1},\ldots,\lambda_{m+n})
=&q_{R,t}^{\otimes m}(\lambda_{1},\ldots,\lambda_{m})
q_{I,t}^{\otimes n}(\lambda_{m+1},\ldots,\lambda_{m+n}).
\end{align}
Because (\ref{def:Hilbert_pair}) implies  $q_{I}(\lambda)=-i\ \textup{sgn}(2^{j}\lambda) q_{R}(\lambda)$,
\begin{align*}%\label{inner_qR_qI}
\int_{\mathbb{R}}q_{R}(\lambda)q_{I}(-\lambda)d\lambda
=i\int_{\mathbb{R}}|\widehat{\psi_{R}}(2^{j}\lambda)|^{2}f_{X}(\lambda)\ \textup{sgn}(\lambda)d\lambda = 0.
\end{align*}
Hence, for any $m,n\in \mathbb{N}$, $r\in\{1,2,\ldots,m\wedge n\}$, and $\lambda_{1},\ldots,\lambda_{m+n-2r}\in \mathbb{R}$,
\begin{align}\notag
q_{R,t}^{\otimes m}\otimes_{r}q_{I,t}^{\otimes n}(\lambda_{1},\ldots,\lambda_{m+n-2r})=0.
\end{align}
It implies that (\ref{product_formula_IRII}) can be simplified as
\begin{equation}\label{proof:ImIn}
I_{m}\left(q_{R,t}^{\otimes m}\right)I_{n}\left(q_{I,t}^{\otimes n}\right)
= I_{m+n}(q_{R,t}^{\otimes m} \otimes q_{I,t}^{\otimes n}).
\end{equation}
By (\ref{special_Cmn}) and (\ref{proof:ImIn}), the series (\ref{complex_U1X}) can be rewritten as
\begin{align}\notag
U^{A}[j]X(t)=&
\sigma_{j}^{\nu}\underset{m,n\in 2\mathbb{N}\cup \{0\}}{\sum}\
c_{A,\frac{m+n}{2}} \frac{h_{m}h_{n}}{\sqrt{m!n!}}\sigma_{j}^{-(m+n)}I_{m+n}(q_{R,t}^{\otimes m}\otimes q_{I,t}^{\otimes n})
\\\label{series:U1}=&
\sigma_{j}^{\nu}c_{A,0}+\sigma_{j}^{\nu}
\underset{\ell\in 2\mathbb{N}}{\sum}
I_{\ell}\left(\mathbb{Q}^{(\ell)}_{t,j}\right),
\end{align}
where $\mathbb{Q}^{(\ell)}_{t,j}: \mathbb{R}^{\ell}\rightarrow \mathbb{C}$ is defined as
\begin{align*}
\mathbb{Q}^{(\ell)}_{t,j} = c_{A,\frac{\ell}{2}} \sigma_{j}^{-\ell}2^{-\frac{\ell}{2}}(-1)^{\frac{\ell}{2}}\underset{\begin{subarray}{c}m,n\in 2\mathbb{N}\cup \{0\}\\ m+n=\ell\end{subarray}}{\sum}\
\left[\left(\frac{m}{2}\right)!\left(\frac{n}{2}\right)!\right]^{-1}
q_{R,t}^{\otimes m}\otimes q_{I,t}^{\otimes n}.
\end{align*}
Because
$q_{I,t}(\cdot)
=-i\ \textup{sgn}(2^{j}\cdot)q_{R,t}(\cdot)=-i\ \textup{sgn}(\cdot)q_{R,t}(\cdot)$, we have
\begin{align}\notag
\mathbb{Q}^{(\ell)}_{t,j}(\lambda_{1},\ldots,\lambda_{\ell}) =&c_{A,\frac{\ell}{2}} \sigma_{j}^{-\ell} 2^{-\frac{\ell}{2}}(-1)^{\frac{\ell}{2}}q_{R,t}^{\otimes \ell}(\lambda_{1},\ldots,\lambda_{\ell})
\\\label{proof:Q}&\times\underset{\begin{subarray}{c}m,n\in 2\mathbb{N}\cup \{0\}\\ m+n=\ell\end{subarray}}{\sum}\
\left[\left(\frac{m}{2}\right)!\left(\frac{n}{2}\right)!\right]^{-1}
(-1)^{\frac{n}{2}}\overset{\ell}{\underset{k=\ell-n+1}{\prod}} \textup{sgn}(\lambda_{k})
\end{align}
for all $(\lambda_{1},\ldots,\lambda_{\ell})\in \mathbb{R}^{\ell}$. By default,
$\overset{\ell}{\underset{k=\ell+1}{\prod}} \textup{sgn}(\lambda_{k})=1$.

Denote the symmetrization of $\mathbb{Q}^{(\ell)}_{t,j}$ with respect to $(\lambda_{1},\ldots,\lambda_{\ell})$ by $Q^{(\ell)}_{t,j}$, which is defined as
\begin{align*}%\label{def:Qpsisym}
Q^{(\ell)}_{t,j}(\lambda_{1},\ldots,\lambda_{\ell})
=\frac{1}{\ell!} \underset{p\in P[\ell]}{\sum}\mathbb{Q}^{(\ell)}_{t,j}(\lambda_{p(1)},\lambda_{p(2)},\ldots,\lambda_{p(\ell)}),
\end{align*}
where $P[\ell]$ represents the set of permutations of $\{1,2,\ldots,\ell\}$.
By the property
\begin{equation*}%\label{relation:Q_Qsym}
I_{\ell}\left(\mathbb{Q}^{(\ell)}_{t,j}\right)=
 I_{\ell}\left(Q^{(\ell)}_{t,j}\right),
\end{equation*}
the series representation (\ref{series:U1}) for $U^{A}[j]X(t)$ can be rewritten as
\begin{align}\label{series:U1_v2}
U^{A}[j]X(t)=
\sigma_{j}^{\nu}c_{A,0}+\sigma_{j}^{\nu}
\underset{\ell\in 2\mathbb{N}}{\sum}
I_{\ell}\left(Q^{(\ell)}_{t,j}\right).
\end{align}
From (\ref{proof:Q}),
\begin{align}\label{Qpsisym_comp1}
Q^{(\ell)}_{t,j}(\lambda_{1},\lambda_{2},\ldots,\lambda_{\ell})
=c_{A,\frac{\ell}{2}} \sigma_{j}^{-\ell}2^{-\frac{\ell}{2}}(-1)^{\frac{\ell}{2}}
q_{R,t}^{\otimes \ell}(\lambda_{1:\ell})
B(\ell,\lambda_{1:\ell}),
\end{align}
where $\lambda_{1:\ell}=(\lambda_{1},\ldots,\lambda_{\ell})$  and
\begin{align*}%\label{def:B}
B(\ell,\lambda_{1:\ell}) = \frac{1}{\ell!} \underset{p\in P[\ell]}{\sum}\ \underset{\begin{subarray}{c}m,n\in 2\mathbb{N}\cup \{0\}\\ m+n=\ell\end{subarray}}{\sum}\
\left[\left(\frac{m}{2}\right)!\left(\frac{n}{2}\right)!\right]^{-1}
(-1)^{\frac{n}{2}}\overset{\ell}{\underset{k=\ell-n+1}{\prod}} \textup{sgn}(\lambda_{p(k)}).
\end{align*}
For any positive even integer $\ell$ and $\{\lambda_{1},\lambda_{2},\ldots,\lambda_{\ell}\}\subset \mathbb{R}\setminus \{0\}$,
Lemma \ref{lemma:B_identity} shows that
\begin{align}\label{exact:B_inproof}
B(\ell,\lambda_{1:\ell})
=\left\{\begin{array}{ll}
 2^\ell\   \frac{\ell}{2}!\ (\ell!)^{-1}& \ \textup{if}\ N(\lambda_{1:\ell})= \ell/2,\\
0 &\ \textup{if}\ N(\lambda_{1:\ell})\neq \ell/2,
\end{array}
\right.
\end{align}
where $N(\lambda_{1:\ell})$ is the number of negative elements in $\{\lambda_{k}\}_{k=1}^{\ell}$.
By substituting (\ref{exact:B_inproof}) into (\ref{Qpsisym_comp1}),
\begin{align}\label{pre_mergeQ}
Q^{(\ell)}_{t,j}(\lambda_{1},\lambda_{2},\ldots,\lambda_{\ell})
=(-2)^{\frac{\ell}{2}}\ (\frac{\ell}{2}!)(\ell!)^{-1}\ c_{A,\frac{\ell}{2}}
\sigma_{j}^{-\ell}q_{R,t}^{\otimes \ell}(\lambda_{1:\ell})
1_{\{N(\lambda_{1:\ell})= \ell/2\}}.
\end{align}
The proof of Theorem \ref{thm:U1_chaosexpansion}
is concluded by
substituting (\ref{pre_mergeQ}) into (\ref{series:U1_v2})
and noticing that the Lebesque measure of
$\{(\lambda_{1},\lambda_{2},\ldots,\lambda_{\ell})\in \mathbb{R}^{\ell}\mid \lambda_{k}=0\ \textup{for some}\ k\}$ is zero.

\section{Proof of Corollary \ref{prop:application}}\label{sec:proof:prop:application}
In the following, we only give details for the case $A(r)=\ln(r)$. The proof for $A(r)=r^\nu$ is similar and we omit it.
Because $h(x) = x$ is one of the test functions for the Wasserstein metric,
\begin{align}\label{proof:application:beta1beta2}
d_{\textup{W}}\left(U^{A}[j]X_{1},U^{A}[j]X_{2}\right)
\geq& \left|\mathbb{E}\left[U^{A}[j]X_{1}\right]-\mathbb{E}\left[U^{A}[j]X_{2}\right]\right|
=
\left|\ln \sigma_{1,j}-\ln \sigma_{2,j}\right|,
\end{align}
where
the last equality follows from Theorem \ref{thm:U1_chaosexpansion}.
If the spectral densities of $X_{1}$ and $X_{2}$ have the form (\ref{condition_f:v2}) and
the wavelet function $\psi_{R}$ satisfies Assumption \ref{Assumption:wavelet},
then by the dominated convergence theorem,
\begin{align*}
\underset{j\rightarrow\infty}{\lim}2^{-\beta_{p} j/2}\sigma_{p,j} = \left[C_{X_{p}}(0)\int_{\mathbb{R}}|\widehat{\psi_{R}}(\lambda)|^{2}|\lambda|^{\beta_{p}-1}d\lambda\right]^{\frac{1}{2}},
\end{align*}
that is,
\begin{align}\label{j:limit}
\underset{j\rightarrow\infty}{\lim}-\frac{1}{2}(\ln2) \beta_{p} j +\ln\sigma_{p,j} = \frac{1}{2}\ln\left[C_{X_{p}}(0)\int_{\mathbb{R}}|\widehat{\psi_{R}}(\lambda)|^{2}|\lambda|^{\beta_{p}-1}d\lambda\right]
\end{align}
for $p=1,2$.
By (\ref{proof:application:beta1beta2}), we obtain
\begin{align}\notag
\underset{j\rightarrow\infty}{\underline{\lim}} j^{-1}d_{\textup{W}}(U^{A}[j]X_{1},U^{A}[j]X_{2})
\geq \underset{j\rightarrow\infty}{\underline{\lim}}
\left|\frac{\ln \sigma_{1,j}}{j}-\frac{\ln \sigma_{2,j}}{j}\right|
=\frac{1}{2}(\ln2)\left|\beta_{1}-\beta_{2}
\right|,
\end{align}
where the equality follows from (\ref{j:limit}).

\section{Proof of Lemma \ref{lemma:gp_bound}}\label{sec:proof:lemma:gp_bound}

{\bf(a)}
First of all, because $\|h^{'}\|_{\infty}\leq1$,
\begin{align*}
\left|\mathbb{E}\left[h\left(\mathbf{F}\right)\right]-\mathbb{E}\left[h\left(\mathbf{F}_{\leq K}\right)\right]\right|
\leq
\mathbb{E}\left[\overset{d}{\underset{m=1}{\sum}}|F_{m,> K}|\right].
\end{align*}
By the orthogonal property
\begin{align*}%\label{ortho_Ip_Iq}
\mathbb{E}\left[I_{\ell}(s^{(\ell)}_{2^{J}t_{m},j_{m}})I_{\ell^{'}}(s^{(\ell^{'})}_{2^{J}t_{m},j_{m}})\right] = \left\{\begin{array}{ll}\ell! \|s^{(\ell)}_{2^{J}t_{m},j_{m}}\|^{2}_{2}
&  \textup{if}\  \ell=\ell^{'},\\
0 &  \textup{if}\  \ell\neq\ell^{'},
\end{array}
\right.
\end{align*}
and (\ref{thm:S1_chaosexpansion_part2}), we have
\begin{align}\notag
\left(\mathbb{E}|F_{m,>K}|\right)^{2}\leq&\mathbb{E}\left[|F_{m,>K}|^{2}\right]
\\\label{norm:F>K:estimate:initial}=&2^{J}\underset{\ell\in\{K+2,K+4,\ldots\}}{\sum}
\ell!\ \|s^{(\ell)}_{2^{J}t_{m},j_{m}}\|^{2}_{2}.
\end{align}
For any $t\in \mathbb{R}$ and $j\in \mathbb{Z}$, according to the definition of $s^{(\ell)}_{t,j}$ in  (\ref{integrand_s}),
\begin{align}\notag
\|s^{(\ell)}_{t,j}\|_{2}^{2}
\leq&\sigma_{j}^{2\nu-2\ell}c_{\ell}^{2}\int_{\mathbb{R}^{\ell}} \left[\overset{\ell}{\underset{k=1}{\prod}}f_{X\star\psi_{R,j}}(\lambda_{k})\right]|\widehat{\phi_{J}}(\lambda_{1}+\cdots+\lambda_{\ell})|^{2}d\lambda_{1}\cdots\lambda_{\ell}
\\\label{bound_gp_v1}=&
2^{-J}\sigma_{j}^{2\nu-2\ell}c_{\ell}^{2}
\int_{\mathbb{R}} f^{\star\ell}_{X\star\psi_{R,j}}(2^{-J}\eta)|\widehat{\phi}(\eta)|^{2}d\eta,
\end{align}
where
\begin{align*}
f_{X\star\psi_{R,j}}(\lambda_{k})
=\left|\widehat{\psi_{R}}(2^{j}\lambda_{k})\right|^{2}
f_{X}(\lambda_{k})
\end{align*}
and $f^{\star\ell}_{X\star\psi_{R,j}}$ is the $\ell$-fold convolution of $f_{X\star\psi_{R,j}}$
with itself.
By the nonnegativity of $f_{X\star\psi_{R,j}}$ and $\sigma_{j}^{2} = \int_{\mathbb{R}}f_{X\star\psi_{R,j}}(\lambda)d\lambda$,
\begin{align}\notag
f^{\star\ell}_{X\star\psi_{R,j}}(\eta) = \int_{\mathbb{R}}f^{\star(\ell-1)}_{X\star\psi_{R,j}}(\eta-\zeta)f_{X\star\psi_{R,j}}(\zeta)d\zeta
\leq \|f^{\star(\ell-1)}_{X\star\psi_{R,j}}\|_{\infty} \sigma_{j}^{2}.
\end{align}
It implies that for any $\ell\in \mathbb{N}$,
\begin{align}\label{bound_fconvp}
\|f^{\star\ell}_{X\star\psi_{R,j}}\|_{\infty}
\leq \|f_{X\star\psi_{R,j}}\|_{\infty} \sigma_{j}^{2(\ell-1)}.
\end{align}
Note that $\|f_{X\star\psi_{R,j}}\|_{\infty}<\infty$ under the assumption $2\alpha+\beta\geq1$.
By using (\ref{bound_fconvp}) to bound the integrand in (\ref{bound_gp_v1}),
\begin{align}\label{estimate_s_ell_proof}
\|s^{(\ell)}_{t,j}\|_{2}^{2}
\leq&
2^{-J}\sigma_{j}^{2\nu-2}c_{\ell}^{2}
\|f_{X\star\psi_{R,j}}\|_{\infty} \|\widehat{\phi}\|^{2}_{2}.
\end{align}
By applying (\ref{estimate_s_ell_proof}) to (\ref{norm:F>K:estimate:initial}),
\begin{align}\label{norm:F>K:estimate:final}
\left(\mathbb{E}|F_{m,>K}|\right)^{2}
\leq
\sigma_{j_{m}}^{2\nu-2}\|f_{X\star\psi_{R,j_{m}}}\|_{\infty}\|\widehat{\phi}\|^{2}_{2}
\underset{\ell\in\{K+2,K+4,\ldots\}}{\sum}\ell!
c_{\ell}^{2}.
\end{align}
From the definition of $c_{\ell}$ in (\ref{def:c_ell}),
\begin{align}\label{proof:sum_cell_factorial}
\underset{\ell\in\{K+2,K+4,\ldots\}}{\sum}\ell!
c_{\ell}^{2} =
\underset{\ell\in\{K+2,K+4,\ldots\}}{\sum}2^{\ell}\left(\frac{\ell}{2}!\right)^{2} (\ell!)^{-1}c_{A,\frac{\ell}{2}}^{2}.
\end{align}
For $A(r)=r^{\nu}$ with $\nu>0$,
$$
c_{A,\frac{\ell}{2}} = 2^{\frac{\nu}{2}}\Gamma(\frac{\nu}{2}+1)\binom{\frac{\ell}{2}-\frac{\nu}{2}-1}{\frac{\ell}{2}}.
$$
By L'Hospital's rule,
\begin{align}\notag
\frac{\ln |c_{A,\frac{\ell}{2}+1}|-\ln |c_{A,\frac{\ell}{2}}|}{\ln(\frac{\ell}{2}+1)-\ln \frac{\ell}{2}}
=&
\frac{\ln |\binom{\frac{\ell}{2}-\frac{\nu}{2}}{\frac{\ell}{2}+1}|-\ln |\binom{\frac{\ell}{2}-\frac{\nu}{2}-1}{\frac{\ell}{2}}|}{\ln(\frac{\ell}{2}+1)-\ln \frac{\ell}{2}}
\\\label{limit_binom}=&\frac{\ln(\frac{\ell}{2}-\frac{\nu}{2})-\ln(\frac{\ell}{2}+1)}{\ln(\frac{\ell}{2}+1)-\ln \frac{\ell}{2}}
\rightarrow
-\frac{\nu}{2}-1
\end{align}
as $\ell\rightarrow\infty$.
From (\ref{limit_binom}), we know that for any $\varepsilon>0$, there exists a constant $C_{1}(\nu,\varepsilon)>0$
such that
\begin{equation}\label{limit_binom2}
|c_{A,\frac{\ell}{2}}| \leq C_{1}(\nu,\varepsilon) \ell^{-\frac{\nu}{2}-1+\varepsilon}
\end{equation}
for any integer $\ell\in 2\mathbb{N}$. For $A(r)=\ln(r)$, $c_{A,\frac{\ell}{2}}=-\ell^{-1}$ for $\ell\in 2\mathbb{N}$,
which can be viewed as a special case of (\ref{limit_binom2}) with $\nu=0$ and $\varepsilon=0$.
On the other hand,
by Stirling's formula, or more precisely
\begin{align}\label{Stirling_formula}
\sqrt{2\pi n}\left(\frac{n}{e}\right)^{n} e^{\frac{1}{12n+1}} <n!<\sqrt{2\pi n}\left(\frac{n}{e}\right)^{n} e^{\frac{1}{12n}}
\end{align}
for all $n\in \mathbb{N}$,
\begin{equation}\label{bound_stirling}
2^{\ell}\left(\frac{\ell}{2}!\right)^{2} (\ell!)^{-1} \leq (2\pi \ell)^{\frac{1}{2}}
\end{equation}
for any integer $\ell\in 2\mathbb{N}$.
By applying the inequalities  (\ref{limit_binom2}) and (\ref{bound_stirling}) to (\ref{proof:sum_cell_factorial}),
there exists a constant $C_{2}(\nu,\varepsilon)>0$ such that
\begin{align}\notag
\underset{\ell\in\{K+2,K+4,\ldots\}}{\sum}\ell!
c_{\ell}^{2} \leq&
(C_{1}(\nu,\varepsilon))^{2}\underset{\ell\in\{K+2,K+4,\ldots\}}{\sum}(2\pi \ell)^{\frac{1}{2}} \ell^{-\nu-2+2\varepsilon}
\\\label{proof:sum_cell_factorial_v2}\leq&
C_{2}(\nu,\varepsilon) K^{-\nu-\frac{1}{2}+2\varepsilon}
\end{align}
for all $K\in 2\mathbb{N}$.
By combining (\ref{norm:F>K:estimate:final}) and (\ref{proof:sum_cell_factorial_v2}), we obtain
\begin{align*}
\left(\mathbb{E}|F_{m,>K}|\right)^{2}
\leq
\sigma_{j_{m}}^{2\nu-2}\|f_{X\star\psi_{R,j_{m}}}\|_{\infty}\|\widehat{\phi}\|^{2}_{2}
C_{2}(\nu,\varepsilon) K^{-\nu-\frac{1}{2}+2\varepsilon},
\end{align*}
which leads to the inequality (\ref{norm_F>K}).
The constant $C$ in (\ref{norm_F>K}) is given by
\begin{align}\label{norm:F>K:estimate:end}
C = \left[\underset{1\leq m\leq d}{\max} \sigma_{j_{m}}^{2\nu-2}\|f_{X\star\psi_{R,j_{m}}}\|_{\infty}\|\widehat{\phi}\|^{2}_{2}
C_{2}(\nu,\varepsilon)\right]^{1/2}.
\end{align}

%\subsection{Proof of Lemma \ref{lemma:norm_N-N}}\label{sec:proof:lemma:N-N}
{\bf(b)}  To explore $\mathbf{N}_{\leq K}$, recall that $\mathbf{N}_{\leq K}$ has the
same covariance structure as $\mathbf{F}_{\leq K}$.
By the Wiener chaos decomposition of $\mathbf{F}_{\leq K}$ in (\ref{thm:S1_chaosexpansion_part1}),
%$S_{J}[j]X$ in (\ref{thm:S1_chaosexpansion}),
for any $m,n\in\{1,\ldots,d\}$,
\begin{align*}
\mathbb{E}\left[F_{m,\leq K},F_{n,\leq K}\right]
= 2^{J}
\underset{\ell\in\{2,4,\ldots,K\}}{\sum}\ \ell! \langle s^{(\ell)}_{2^{J}t_{m},j_{m}},s^{(\ell)}_{2^{J}t_{n},j_{n}}\rangle,
\end{align*}
where
\begin{align*}
\langle s^{(\ell)}_{2^{J}t_{m},j_{m}},s^{(\ell)}_{2^{J}t_{n},j_{n}}\rangle
=\int_{\mathbb{R}^{\ell}}s^{(\ell)}_{2^{J}t_{m},j_{m}}(\lambda_{1},\ldots,\lambda_{\ell})
\overline{s^{(\ell)}_{2^{J}t_{n},j_{n}}(\lambda_{1},\ldots,\lambda_{\ell})}d\lambda_{1}\cdots d\lambda_{\ell}.
\end{align*}
Let
$$\left\{\mathbf{N}^{(\ell)}\right\}_{\ell\in 2\mathbb{N}}=\left\{\left(N^{(\ell)}_{1},\ldots,N^{(\ell)}_{d}\right)\right\}_{\ell\in 2\mathbb{N}}$$
be a sequence of independent $d$-dimensional normal random vectors
with mean zero and
$$\mathbb{E}\left[N^{(\ell)}_{m}N^{(\ell)}_{n}\right] =2^{J}\ell! \langle s^{(\ell)}_{2^{J}t_{m},j_{m}},s^{(\ell)}_{2^{J}t_{n},j_{n}}\rangle$$
for $m,n\in\{1,\ldots,d\}$.
Then,
$$\mathbf{N}\overset{d}{=}\underset{\ell\in 2\mathbb{N}}{\sum}\mathbf{N}^{(\ell)}$$
and
$$\mathbf{N}_{\leq K}\overset{d}{=}\underset{\ell\in \{2,4,\ldots,K\}}{\sum}\mathbf{N}^{(\ell)},$$
where $\overset{d}{=}$ means the equality is in the distribution sense.
For any Lipschitz differentiable function $h:\mathbb{R}^{d}\rightarrow \mathbb{R}$ with Lipschitz constant $\|h^{'}\|_{\infty}\leq1$,
\begin{align}\notag
&\left|\mathbb{E}\left[h\left(\mathbf{N}_{\leq K}\right)\right]-\mathbb{E}\left[h\left(\mathbf{N}\right)\right] \right|
\leq
\overset{d}{\underset{m=1}{\sum}}\mathbb{E}\left|\underset{\ell\in\{K+2,K+4,\ldots\}}{\sum}N^{(\ell)}_{m}\right|
\\\label{proof:N-N:final}=& \sqrt{\frac{2}{\pi}}
\overset{d}{\underset{m=1}{\sum}}\left[\underset{\ell\in\{K+2,K+4,\ldots\}}{\sum}\mathbb{E} |N^{(\ell)}_{m}|^{2}\right]^{\frac{1}{2}}
=\sqrt{\frac{2}{\pi}}
\overset{d}{\underset{m=1}{\sum}}\left[\underset{\ell\in\{K+2,K+4,\ldots\}}{\sum}
2^{J}\ell! \|s^{(\ell)}_{2^{J}t_{m},j_{m}}\|^{2}_{2}
\right]^{\frac{1}{2}}.
\end{align}
%where the first inequality follows from
%the independence of the sequence $\{N^{(\ell)}_{m}\}_{\ell\in 2\mathbb{N}}$ for each fixed $m\in\{1,\ldots,d\}$.
The summation inside the square brackets in the last term is the same as the upper bound of $\left(\mathbb{E}|F_{m,>K}|\right)^{2}$
in (\ref{norm:F>K:estimate:initial}).
Hence, we can directly apply (\ref{norm_F>K}) to the right hand side of (\ref{proof:N-N:final}) to obtain
\begin{align}\notag
&\left|\mathbb{E}\left[h\left(\mathbf{N}_{\leq K}\right)\right]-\mathbb{E}\left[h\left(\mathbf{N}\right)\right] \right|
\leq d C K^{-\frac{\nu}{2}-\frac{1}{4}+\varepsilon}.
\end{align}

\section{Proof of Lemma \ref{cited_thm:multiGA}}\label{sec:proof:lemma:cited_thm:multiGA}
For any twice differentiable function $h$,  \cite[Theorem 6.1.2]{nourdin2012normal} shows that
\begin{align}\label{multidim_stein_v2}
|\mathbb{E}[h(\mathbf{S})]-\mathbb{E}[h(\mathbf{N}_{\mathbf{C}})]|\leq
\frac{1}{2}
\left[\underset{1\leq m,n\leq d}{\max}\ \underset{\mathbf{x}\in \mathbb{R}^{d}}{\sup}\ \Big|\frac{\partial^{2}h}{\partial x_{m} \partial x_{n}}(\mathbf{x})\Big|\right]\rho.
\end{align}
Because
\begin{align}\label{relation:C2norm_Lip}
\underset{1\leq m,n\leq d}{\max}\ \underset{\mathbf{x}\in \mathbb{R}^{d}}{\sup}\ \Big|\frac{\partial^{2}h}{\partial x_{m} \partial x_{n}}(\mathbf{x})\Big|\leq\|h^{'}\|_{\textup{Lip}}
\end{align}
for any twice differentiable function $h$ in $\mathcal{H}_{2}$,
it suffices to show that (\ref{multidim_stein}) holds for non-twice differentiable functions in $\mathcal{H}_{2}$.
For such class of functions $h$,
we define
\begin{align*}
h_{\varepsilon}(\mathbf{x}) = \int_{\mathbb{R}^{d}}\frac{1}{(2\pi \varepsilon)^{\frac{d}{2}}}e^{-\frac{\|\mathbf{z}\|^{2}}{2\varepsilon}}
h(\mathbf{x}-\mathbf{z})d\mathbf{z},\ \mathbf{x}\in \mathbb{R}^{d}.
\end{align*}
%which is a smooth approximation of $h$.
For any $\varepsilon>0$, $h_{\varepsilon}$ is twice differentiable, so (\ref{multidim_stein_v2}) and (\ref{relation:C2norm_Lip})
imply that
\begin{align}\label{Wasserstein_smooth_approx}
|\mathbb{E}[h_{\varepsilon}(\mathbf{S})]-\mathbb{E}[h_{\varepsilon}(\mathbf{N}_{\mathbf{C}})]|\leq
&\frac{1}{2}
\|h_{\varepsilon}^{'}\|_{\textup{Lip}}\rho.
\end{align}
Because $h$ is continuous and bounded, by the dominated convergence theorem,
\begin{align}\label{Wasserstein_smooth_approx_LHS}
\underset{\varepsilon\rightarrow 0}{\lim} |\mathbb{E}[h_{\varepsilon}(\mathbf{S})]-\mathbb{E}[h_{\varepsilon}(\mathbf{N}_{\mathbf{C}})]|
= |\mathbb{E}[h(\mathbf{S})]-\mathbb{E}[h(\mathbf{N}_{\mathbf{C}})]|.
\end{align}
On the other hand,
because $h$ is continuously differentiable and $\|\frac{\partial h}{\partial x_{m}}\|_{\infty}\leq1$ for $m\in \{1,2,\ldots,d\}$,
\begin{align*}
\frac{\partial h_{\varepsilon}}{\partial x_{m}}(\mathbf{x}) = \int_{\mathbb{R}^{d}}\frac{1}{(2\pi \varepsilon)^{\frac{d}{2}}}e^{-\frac{\|\mathbf{z}\|^{2}}{2\varepsilon}}
\frac{\partial h}{\partial x_{m}}(\mathbf{x}-\mathbf{z})d\mathbf{z}.
\end{align*}
This implies that for any $\mathbf{x},\mathbf{y}\in \mathbb{R}^{d}$ with $\mathbf{x}\neq\mathbf{y}$,
\begin{align}\notag
\left|\frac{\partial h_{\varepsilon}}{\partial x_{m}}(\mathbf{x}) -\frac{\partial h_{\varepsilon}}{\partial x_{m}}(\mathbf{y})\right|
&\leq
\int_{\mathbb{R}^{d}}\frac{1}{(2\pi \varepsilon)^{\frac{d}{2}}}e^{-\frac{\|\mathbf{z}\|^{2}}{2\varepsilon}}
\left|\frac{\partial h}{\partial x_{m}}(\mathbf{x}-\mathbf{z})-\frac{\partial h}{\partial x_{m}}(\mathbf{y}-\mathbf{z})\right|d\mathbf{z},
\\\label{Lip:derivative_hsmooth}&\leq
\|\frac{\partial h}{\partial x_{m}}\|_{\textup{Lip}}\|\mathbf{x}-\mathbf{y}\|.
\end{align}
The inequality (\ref{Lip:derivative_hsmooth}) implies that the term $\|h_{\varepsilon}^{'}\|_{\textup{Lip}}$
on the right hand side of (\ref{Wasserstein_smooth_approx}) can be bounded as follows
\begin{align}\label{Wasserstein_smooth_approx_RHS}
\underset{\varepsilon\rightarrow 0}{\lim\sup}\|h_{\varepsilon}^{'}\|_{\textup{Lip}}=
\underset{1\leq m\leq d}{\max}\underset{\varepsilon\rightarrow 0}{\lim\sup}\ \underset{\begin{subarray}{c}\mathbf{x},\mathbf{y}\in \mathbb{R}^{d}\\ \mathbf{x}\neq\mathbf{y}\end{subarray}}{\sup}\ \frac{\Big|\frac{\partial h_{\varepsilon}}{\partial x_{m}}(\mathbf{x})-\frac{\partial h_{\varepsilon}}{\partial x_{m}}(\mathbf{y})\Big|}{\|\mathbf{x}-\mathbf{y}\|}
\leq
\underset{m\in\{1,2,\ldots,d\}}{\max}\|\frac{\partial h}{\partial x_{m}}\|_{\textup{Lip}}.
\end{align}
By applying (\ref{Wasserstein_smooth_approx_LHS}) and (\ref{Wasserstein_smooth_approx_RHS})
to (\ref{Wasserstein_smooth_approx}),
we conclude that  (\ref{multidim_stein}) holds for any $h\in \mathcal{H}_{2}$.

\section{Proof of Proposition \ref{prop:multidim_stein_RHS}}\label{sec:proof:prop:multidim_stein_RHS}
The proof of  (\ref{Prop:statement:stein:nu})  consists of three parts.
First, we rewrite $$\mathbb{E}\left[\left(\mathbf{C}_{J,K}(m,n)-\Big\langle DF_{n,\leq K},-DL^{-1}F_{m,\leq K}\Big\rangle\right)^{2}\right]$$ by making use of the Wiener chaos decomposition of $\mathbf{F}_{\leq K}$ in (\ref{thm:S1_chaosexpansion_part1}) as follows.
By the orthogonal property of Wiener-It$\hat{\textup{o}}$ integrals, the covariance $\mathbf C_{J,K}(m,n)$ can be expressed as
\begin{align}\notag
\mathbb{E}\left[F_{m,\leq K}F_{n,\leq K}\right]
=&2^{J}\underset{\ell\in \{2,4,\ldots,K\}}{\sum}\mathbb{E}\left[I_{\ell}\left(s^{(\ell)}_{2^{J}t_{m},j_{m}}\right)
I_{\ell}\left(s^{(\ell)}_{2^{J}t_{n},j_{n}}\right)
\right]
\\\label{def:T:proofuse}=&2^{J}\underset{\ell\in \{2,4,\ldots,K\}}{\sum}\ell! \
s^{(\ell)}_{2^{J}t_{m},j_{m}}\otimes_{\ell}
s^{(\ell)}_{2^{J}t_{n},j_{n}},
\end{align}
where $s^{(\ell)}$ is defined in (\ref{integrand_s}),
and the contraction operator $\otimes_{\ell}$ is defined in (\ref{def:contraction}).
By Lemma \ref{lemma:deltaD} and (\ref{thm:S1_chaosexpansion_part1}),
\begin{align}\label{check:condition2}
DF_{n,\leq K}
=2^{\frac{J}{2}}\underset{\ell\in\{2,4,\ldots,K\}}{\sum}\ \ell I_{\ell-1}\left(s^{(\ell)}_{2^{J}t_{n},j_{n}}\right).
\end{align}
By the property of $L^{-1}$ mentioned in (\ref{def:inverseL}),
\begin{align}\notag
-DL^{-1}F_{m,\leq K}
=&2^{\frac{J}{2}}D\left[
\underset{\ell\in\{2,4,\ldots,K\}}{\sum}\ \frac{1}{\ell} I_{\ell}\left(s^{(\ell)}_{2^{J}t_{m},j_{m}}\right)\right]
\\\label{DLinverseF}=&
2^{\frac{J}{2}}
\underset{\ell\in\{2,4,\ldots,K\}}{\sum}\  I_{\ell-1}\left(s^{(\ell)}_{2^{J}t_{m},j_{m}}\right).
\end{align}
By (\ref{check:condition2}), (\ref{DLinverseF}), and the product formula in Lemma \ref{lemma:itoformula},
\begin{align}\notag
&\langle DF_{n,\leq K},-DL^{-1}F_{m,\leq K}\rangle
=2^{J}
\Big\langle
\underset{\ell\in\{2,4,\ldots,K\}}{\sum}\ \ell I_{\ell-1}\left(s^{(\ell)}_{2^{J}t_{n},j_{n}}\right),\underset{\ell\in\{2,4,\ldots,K\}}{\sum}\  I_{\ell-1}\left(s^{(\ell)}_{2^{J}t_{m},j_{m}}\right)
\Big\rangle
\\\notag=\ &
2^{J}\underset{\ell\in \{2,4,\ldots,K\}}{\sum}\ell! \
s^{(\ell)}_{2^{J}t_{m},j_{m}}\otimes_{\ell}
s^{(\ell)}_{2^{J}t_{n},j_{n}}
\\\notag&+
2^{J}\underset{\ell\in\{2,4,\ldots,K\}}{\sum}\ \ell \overset{\ell-1}{\underset{r=1}{\sum}}(r-1)!\binom{\ell-1}{r-1}^{2}
I_{2\ell-2r}
\left(s^{(\ell)}_{2^{J}t_{m},j_{m}}\otimes_{r}
s^{(\ell)}_{2^{J}t_{n},j_{n}}\right)
\\\label{def:T:cross}
&+
2^{J}\underset{\begin{subarray}{c}
\ell,\ell'\in\{2,4,\ldots,K\} \\
\ell\neq \ell'
\end{subarray}}{\sum}\ell
\overset{\ell\wedge \ell'}{\underset{r=1}{\sum}}(r-1)!\binom{\ell-1}{r-1}\binom{\ell^{'}-1}{r-1}
I_{\ell+\ell^{'}-2r}\left(s^{(\ell)}_{2^{J}t_{m},j_{m}}\otimes_{r}
s^{(\ell)}_{2^{J}t_{n},j_{n}}\right).
\end{align}
By (\ref{def:T:proofuse}),   (\ref{def:T:cross}), and the Minkowski inequality,
\begin{align}\label{estimate_C-DFDLF}
&\sqrt{\mathbb{E}\left[\left(\mathbf{C}_{J,K}(m,n)-\langle DF_{n,\leq K},-DL^{-1}F_{m,\leq K}\rangle\right)^{2}\right]}
\\\notag\leq&2^{J}\underset{\ell\in\{2,4,\ldots,K\}}{\sum}\ \ell \overset{\ell-1}{\underset{r=1}{\sum}}(r-1)!\binom{\ell-1}{r-1}^{2}
\sqrt{\mathbb{E}\left[\left|I_{2\ell-2r}
\left(s^{(\ell)}_{2^{J}t_{m},j_{m}}\otimes_{r}
s^{(\ell)}_{2^{J}t_{n},j_{n}}\right)\right|^{2}\right]}
\\\notag+&
2^{J}\underset{\begin{subarray}{c}
\ell,\ell'\in\{2,4,\ldots,K\} \\
\ell\neq \ell'
\end{subarray}}{\sum}\ell
\overset{\ell\wedge \ell'}{\underset{r=1}{\sum}}(r-1)!\binom{\ell-1}{r-1}\binom{\ell^{'}-1}{r-1}
\sqrt{\mathbb{E}\left[\left|I_{\ell+\ell^{'}-2r}\left(s^{(\ell)}_{2^{J}t_{m},j_{m}}\otimes_{r}
s^{(\ell^{'})}_{2^{J}t_{n},j_{n}}\right)\right|^{2}\right]}.
\end{align}
According to the definition of the contraction operator $\otimes_{r}$ in (\ref{def:contraction}),
$$f:=s^{(\ell)}_{2^{J}t_{m},j_{m}}\otimes_{r}
s^{(\ell^{'})}_{2^{J}t_{n},j_{n}}
$$
is a function defined on $\mathbb{R}^{\ell+\ell^{'}-2r}$.
Denote
\begin{align*}%\label{def:Qpsisym}
\widetilde{f}(\lambda_{1},\ldots,\lambda_{\ell+\ell^{'}-2r})
=\frac{1}{(\ell+\ell^{'}-2r)!} \underset{p\in P[\ell+\ell^{'}-2r]}{\sum}f(\lambda_{p(1)},\lambda_{p(2)},\ldots,\lambda_{p(\ell+\ell^{'}-2r)}),
\end{align*}
where $P[\ell+\ell^{'}-2r]$ represents the set of permutations of $\{1,2,\ldots,\ell+\ell^{'}-2r\}$.
By the fact that
$$I_{\ell+\ell^{'}-2r}(f)=I_{\ell+\ell^{'}-2r}(\widetilde{f})$$
and the isometry property of multiple Wiener integrals \cite[(4.4)-(4.6)]{major1981lecture},
we have
\begin{align}\notag
\mathbb{E}\left[\left|I_{\ell+\ell^{'}-2r}(f)\right|^{2}\right]
=& \mathbb{E}\left[\left|I_{\ell+\ell^{'}-2r}(\widetilde{f})\right|^{2}\right]
\\\label{isometry_ineq}=&(\ell+\ell^{'}-2r)! \|\widetilde{f}\|_{L^{2}}^{2}\leq(\ell+\ell^{'}-2r)! \|f\|_{L^{2}}^{2}.
\end{align}
By applying (\ref{isometry_ineq}) to the expectations in (\ref{estimate_C-DFDLF}), we obtain
\begin{align*}
\mathbb{E}\left[\left|I_{\ell+\ell^{'}-2r}\left(s^{(\ell)}_{2^{J}t_{m},j_{m}}\otimes_{r}
s^{(\ell^{'})}_{2^{J}t_{n},j_{n}}\right)\right|^{2}\right]
\leq
(\ell+\ell^{'}-2r)!\ \left\|s^{(\ell)}_{2^{J}t_{m},j_{m}}\otimes_{r}
s^{(\ell^{'})}_{2^{J}t_{n},j_{n}}\right\|^{2}_{2}
\end{align*}
for all  $\ell,\ell^{'}\in\{2,4,\ldots,K\}$ and $r\in\Lambda_{\ell,\ell^{'}}$, where $\Lambda_{\ell,\ell}=\{1,2,\ldots,\ell-1\}$ and
$\Lambda_{\ell,\ell^{'}}=\{1,2,\ldots,\ell\wedge \ell^{'}\}$ for the case $\ell\neq \ell^{'}$.
Hence, (\ref{estimate_C-DFDLF}) implies that
\begin{align}\notag
&\sqrt{\mathbb{E}\left[\left(\mathbf{C}_{J,K}(m,n)-\Big\langle DF_{n,\leq K},-DL^{-1}F_{m,\leq K}\Big\rangle\right)^{2}\right]}
\\\label{proof:expansion_part1}\leq&2^{J}\underset{\ell\in\{2,4,\ldots,K\}}{\sum}\ \ell \overset{\ell-1}{\underset{r=1}{\sum}}(r-1)!\binom{\ell-1}{r-1}^{2}
\sqrt{(2\ell-2r)!}\ \left\|s^{(\ell)}_{2^{J}t_{m},j_{m}}\otimes_{r}
s^{(\ell)}_{2^{J}t_{n},j_{n}}\right\|_{2}
\\\notag+&
2^{J}\underset{\begin{subarray}{c}
\ell,\ell'\in\{2,4,\ldots,K\} \\
\ell\neq \ell'
\end{subarray}}{\sum}\ell
\overset{\ell\wedge \ell'}{\underset{r=1}{\sum}}(r-1)!\binom{\ell-1}{r-1}\binom{\ell^{'}-1}{r-1}
\sqrt{(\ell+\ell^{'}-2r)!}\ \left\|s^{(\ell)}_{2^{J}t_{m},j_{m}}\otimes_{r}
s^{(\ell^{'})}_{2^{J}t_{n},j_{n}}\right\|_{2}.
\end{align}

Second, we calculate an upper bound
for $\left\|s^{(\ell)}_{2^{J}t_{m},j_{m}}\otimes_{r}
s^{(\ell^{'})}_{2^{J}t_{n},j_{n}}\right\|_{2}$,
where $s^{(\ell)}_{t,j}$ is defined in (\ref{integrand_s}).
To simplify the tedious notation, for any integers $p_1<p_2$, we denote $(\lambda_{p_1},\,\lambda_{p_1+1},\ldots,\lambda_{p_2})$ by $\lambda_{p_1:p_2}$ and $\lambda_{p_1}+\lambda_{p_1+1}+\ldots+\lambda_{p_2}$ by $\lambda^+_{p_1:p_2}$.
From (\ref{integrand_s}),
\begin{align*}
\left|s^{(\ell)}_{2^{J}t_{m},j_{m}}(\lambda_{1:\ell})\right|
\leq&\sigma_{j_{m}}^{\nu-\ell}|c_{\ell}|
\left[\overset{\ell}{\underset{k=1}{\prod}}|\widehat{\psi_{R}}(2^{j_{m}}\lambda_{k})|\sqrt{f_{X}(\lambda_{k})}\right]|\widehat{\phi_{J}}(\lambda_{1:\ell}^{+})|.
\end{align*}
For any $\ell,\ell^{'}\in 2\mathbb{N}$, $r\in \Lambda_{\ell,\ell^{'}}$ and $t_{m},t_{n}\in \mathbb{R}$,
\begin{align}\notag
&\left|s^{(\ell)}_{2^{J}t_{m},j_{m}}\otimes_{r}s^{(\ell^{'})}_{2^{J}t_{n},j_{n}}\left(\lambda_{1:\ell+\ell^{'}-2r}\right)\right|
\\\notag\leq&
\sigma_{j_{m}}^{\nu-\ell} \sigma_{j_{n}}^{\nu-\ell^{'}} |c_{\ell}c_{\ell^{'}}|
\left[\overset{\ell-r}{\underset{k=1}{\prod}}
|\widehat{\psi_{R}}(2^{j_{m}}\lambda_{k})|\sqrt{f_{X}(\lambda_{k})}\right]
\left[\overset{\ell+\ell^{'}-2r}{\underset{k=\ell-r+1}{\prod}}
|\widehat{\psi_{R}}(2^{j_{n}}\lambda_{k})|\sqrt{f_{X}(\lambda_{k})}\right]
\\\notag\times&\int_{\mathbb{R}^{r}}\hspace{-0.15cm}\left[\overset{r}{\underset{k=1}{\prod}}|\widehat{\psi_{R}}(2^{j_{m}}\tau_{k})
\widehat{\psi_{R}}(2^{j_{n}}\tau_{k})|
f_{X}(\tau_{k})\right]
\hspace{-0.15cm}|\widehat{\phi_{J}}(\lambda_{1:\ell-r}^{+}+\tau_{1:r}^{+})\widehat{\phi_{J}}(\lambda_{\ell-r+1:\ell+\ell^{'}-2r}^{+}-\tau_{1:r}^{+})|
d\tau_{1}\cdots d\tau_{r}.
\end{align}
Hence,
\begin{align*}\notag
&\|s^{(\ell)}_{2^{J}t_{m},j_{m}}\otimes_{r}s^{(\ell^{'})}_{2^{J}t_{n},j_{n}}\|_{2}^{2}
\\\notag\leq&
\sigma_{j_{m}}^{2\nu-2\ell} \sigma_{j_{n}}^{2\nu-2\ell^{'}}c_{\ell}^{2}c_{\ell^{'}}^{2}
\int_{\mathbb{R}^{\ell+\ell^{'}}}
\left[\overset{\ell-r}{\underset{k=1}{\prod}}
|\widehat{\psi_{R}}(2^{j_{m}}\lambda_{k})|^{2}f_{X}(\lambda_{k})\right]
\left[\overset{\ell+\ell^{'}-2r}{\underset{k=\ell-r+1}{\prod}}
|\widehat{\psi_{R}}(2^{j_{n}}\lambda_{k})|^{2}f_{X}(\lambda_{k})\right]
\\\notag&\times
\left[\overset{r}{\underset{k=1}{\prod}}|\widehat{\psi_{R}}(2^{j_{m}}\tau_{k})||\widehat{\psi_{R}}(2^{j_{n}}\tau_{k})|
f_{X}(\tau_{k})\right]
\left[\overset{r}{\underset{k=1}{\prod}}|\widehat{\psi_{R}}(2^{j_{m}}\eta_{k})||\widehat{\psi_{R}}(2^{j_{n}}\eta_{k})|
f_{X}(\eta_{k})\right]
\\\notag&\times|\widehat{\phi_{J}}(\lambda_{1:\ell-r}^{+}+\tau_{1:r}^{+})\widehat{\phi_{J}}(\lambda_{\ell-r+1:\ell+\ell^{'}-2r}^{+}-\tau_{1:r}^{+})|
\\&\times|\widehat{\phi_{J}}(\lambda_{1:\ell-r}^{+}+\eta_{1:r}^{+})
\widehat{\phi_{J}}(\lambda_{\ell-r+1:\ell+\ell^{'}-2r}^{+}-\eta_{1:r}^{+})|\ d\tau_{1}\cdots d\tau_{r}\ d\eta_{1}\cdots d\eta_{r}\
d\lambda_{1}\cdots d\lambda_{\ell+\ell^{'}-2r}.
\end{align*}
By considering the change of variables
\begin{align*}%\label{change_variable}
\left\{\begin{array}{l}
u_{k}=\tau_{k},\ k=1,\ldots,r-1\ \ \textup{if}\ r\geq 2,
\\
v_{k}=\eta_{k},\ k=1,\ldots,r-1\ \ \textup{if}\ r\geq 2,
\\
w_{k} = \lambda_{k},\ k =1,\ldots,\ell+\ell^{'}-2r-1,
\\
x=2^{J}\left(\lambda^{+}_{1:\ell-r}+\tau^{+}_{1:r}\right),
\\
y=2^{J}\left(\lambda^{+}_{\ell-r+1:\ell+\ell^{'}-2r}-\tau^{+}_{1:r}\right),
\\
z=2^{J}\left(\lambda^{+}_{1:\ell-r}+\eta^{+}_{1:r}\right),
\end{array}
\right.
\end{align*}
and noting that $\widehat{\phi_{J}}(\cdot) = \widehat{\phi}(2^{J}\cdot)$,
\begin{align}\label{change_variable_integral}
&\|s^{(\ell)}_{2^{J}t_{m},j_{m}}\otimes_{r}s^{(\ell^{'})}_{2^{J}t_{n},j_{n}}\|_{2}^{2}
\\\notag\leq&
2^{-3J}\sigma_{j_{m}}^{2\nu-2\ell} \sigma_{j_{n}}^{2\nu-2\ell^{'}}c_{\ell}^{2}c_{\ell^{'}}^{2}
\int_{\mathbb{R}^{\ell+\ell^{'}}}
\left[\overset{\ell-r}{\underset{k=1}{\prod}}
|\widehat{\psi_{R}}(2^{j_{m}}w_{k})|^{2}f_{X}w_{k})\right]
\left[\overset{\ell+\ell^{'}-2r-1}{\underset{k=\ell-r+1}{\prod}}
|\widehat{\psi_{R}}(2^{j_{n}}w_{k})|^{2}f_{X}(w_{k})\right]
\\\notag&\times
\left[\overset{r-1}{\underset{k=1}{\prod}}|\widehat{\psi_{R}}(2^{j_{m}}u_{k})||\widehat{\psi_{R}}(2^{j_{n}}u_{k})|
f_{X}(u_{k})\right]
\left[\overset{r-1}{\underset{k=1}{\prod}}|\widehat{\psi_{R}}(2^{j_{m}}v_{k})||\widehat{\psi_{R}}(2^{j_{n}}v_{k})|
f_{X}(v_{k})\right]
\\\notag&\times
|\widehat{\psi_{R}}(2^{j_{m}}u^{*})||\widehat{\psi_{R}}(2^{j_{n}}u^{*})|f_{X}(u^{*})
|\widehat{\psi_{R}}(2^{j_{m}}v^{*})||\widehat{\psi_{R}}(2^{j_{n}}v^{*})|f_{X}(v^{*})
|\widehat{\psi_{R}}(2^{j_{n}}w^{*})|^{2}f_{X}(w^{*})
\\\notag&\times|\widehat{\phi}(x)\widehat{\phi_{J}}(y)\widehat{\phi}(z)
\widehat{\phi}(x+y-z)|\ \ du_{1}\cdots du_{r-1}dv_{1}\cdots d_{r-1}dw_{1}\cdots dw_{\ell+\ell^{'}-2r-1}
dx dy dz,
\end{align}
where
$u^{*} = 2^{-J}x-u_{1}-\cdots-u_{r-1}-w_{1}-\cdots-w_{\ell-r}$,
$v^{*} = 2^{-J}z-v_{1}-\cdots-v_{r-1}-w_{1}-\cdots-w_{\ell-r}$, and
$w^{*} = 2^{-J}x+2^{-J}y-w_{1}-\cdots-w_{\ell+\ell^{'}-2r-1}$.
From (\ref{change_variable_integral}), we obtain
\begin{align}\label{proof:sotimess_end}
\|s^{(\ell)}_{2^{J}t_{m},j_{m}}\otimes_{r}s^{(\ell^{'})}_{2^{J}t_{n},j_{n}}\|_{2}^{2}
\leq
C^{2}2^{-3J}c_{\ell}^{2}c_{\ell^{'}}^{2},
\end{align}
%for any $\ell,\ell^{'}\in 2\mathbb{N}$, $j_{1},\ldots,j_{d}\in \mathbb{Z}$, and $r\in \Lambda_{\ell,\ell^{'}}$,
where
\begin{align}\label{proof:part(b):C}
C^{2} =
\|\widehat{\phi}\|_{\infty}
\|\widehat{\phi}\|_{1}^{3}
\underset{1\leq m,n\leq d}{\max}\{\sigma_{j_{m}}^{2\nu-2} \sigma_{j_{n}}^{2\nu-4}
M^{3},\  \sigma_{j_{m}}^{2\nu-2} \sigma_{j_{n}}^{2\nu-2}
M^{2}\}
\end{align}
and
\begin{align*}%\label{fstarp_supnorm}
M=\underset{1\leq m,n\leq d}{\max}\|\widehat{\psi_{R}}(2^{j_{m}}\cdot)\widehat{\psi_{R}}(2^{j_{n}}\cdot)f_{X}(\cdot)\|_{\infty}<\infty
\end{align*}
under the assumption $2\alpha+\beta\geq1$.

\begin{figure}
  \centering
  \includegraphics[scale=0.5]{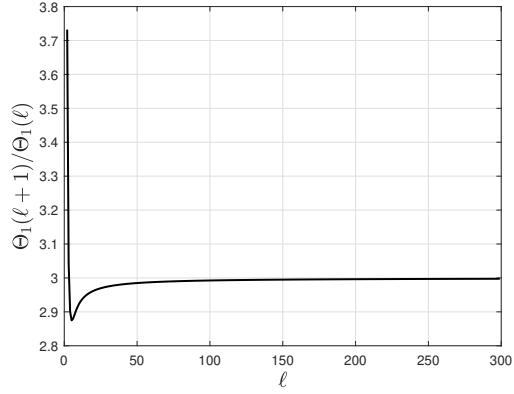}
  \caption{Behavior of the ratio $\Theta_{1}(\ell+1)/\Theta_{1}(\ell)$}
  \label{fig:ratio_theta}
\end{figure}

%{\bf(b)}
Third, we rewrite the summation (\ref{proof:expansion_part1}) in a more compact form by using the inequality (\ref{proof:sotimess_end}).
The summation (\ref{proof:expansion_part1}) can be bounded as follows
%\subsection{Proof of Theorem \ref{corollary:convergence_rate}}
\begin{align}\label{initial_P1+P2}
\sqrt{\mathbb{E}\left[\left(\mathbf{C}_{J,K}(m,n)-\Big\langle DF_{n,\leq K},-DL^{-1}F_{m,\leq K}\Big\rangle\right)^{2}\right]}
\leq C\ 2^{-\frac{1}{2}J}\  \left[P_{1}(K)
+
P_{2}(K)\right],
\end{align}
where
\begin{align}\label{def:P1}
P_{1}(K)=\underset{\ell\in\{2,4,\ldots,K\}}{\sum}\ \ell \overset{\ell-1}{\underset{r=1}{\sum}}(r-1)!\binom{\ell-1}{r-1}^{2}
\sqrt{(2\ell-2r)!}\ c_{\ell}^{2}
\end{align}
and
\begin{align}\label{def:P2}
P_{2}(K)=\underset{\begin{subarray}{c}
\ell,\ell'\in\{2,4,\ldots,K\} \\
\ell\neq \ell'
\end{subarray}}{\sum}\ell
\overset{\ell\wedge \ell'}{\underset{r=1}{\sum}}(r-1)!\binom{\ell-1}{r-1}\binom{\ell^{'}-1}{r-1}
\sqrt{(\ell+\ell^{'}-2r)!}\ |c_{\ell}c_{\ell^{'}}|.
\end{align}
%We {\color{red}rewrite the summations} (\ref{def:P1}) and (\ref{def:P2}) in a more compact form.
Denote
\begin{align*}%\label{def:Theta}
\Theta_{1}(\ell) =\frac{1}{(\ell-1)!} \overset{\ell-1}{\underset{r=1}{\sum}}(r-1)!\binom{\ell-1}{r-1}^{2}\sqrt{(2\ell-2r)!}
\end{align*}
for $\ell\in \mathbb{N}$ with $\Theta_{1}(1) = 0$, which leads to
\begin{align*}
P_{1}(K)=\underset{\ell\in\{2,4,\ldots,K\}}{\sum}\ c_{\ell}^{2}\ \ell!\ \Theta_1(\ell) \,.
\end{align*}
Observe that for every $\ell\in \mathbb{N}$,
\begin{align}\label{bound_theta_v1}
\Theta_{1}(\ell) = \overset{\ell-1}{\underset{r=1}{\sum}}\binom{\ell-1}{r-1}\sqrt{\binom{2\ell-2r}{\ell-r}}
=-1+ \overset{\ell-1}{\underset{k=0}{\sum}}\binom{\ell-1}{k}\sqrt{\binom{2k}{k}}.
\end{align}
Because
$\binom{2k}{k}\leq 4^k$
for all $k\in \mathbb{N}\cup \{0\}$,
(\ref{bound_theta_v1}) implies that
\begin{align}\label{bound_theta_v2}
\Theta_{1}(\ell)
\leq -1+\overset{\ell-1}{\underset{k=0}{\sum}}\binom{\ell-1}{k}2^{k} = -1+3^{\ell-1},\ \ell\in \mathbb{N}.
\end{align}
See Figure \ref{fig:ratio_theta} for the behavior of the ratio
$\Theta_{1}(\ell+1)/\Theta_{1}(\ell)$, which tends to 3 when $\ell\rightarrow\infty$.
By substituting the estimate (\ref{bound_theta_v2}) into (\ref{def:P1}),
\begin{align}\label{proof:P1_final}
P_{1}(K)
\leq& \underset{\ell\in\{2,4,\ldots,K\}}{\sum}c_{\ell}^{2}\ \ell!
\ 3^{\ell-1}.
\end{align}
On the other hand,
\begin{align}\label{boundS2_proof_start}
P_{2}(K) =& \underset{\begin{subarray}{c}
\ell,\ell^{'}\in\{2,4,\ldots,K\}\\
\ell\neq \ell^{'}
\end{subarray}}{\sum}|c_{\ell}c_{\ell^{'}}|\sqrt{\ell!}\sqrt{\ell^{'}!}
\Theta_{2}(\ell,\ell^{'}),
\end{align}
where
\begin{align}\notag
\Theta_{2}(\ell,\ell^{'})=&
\ell\overset{\ell\wedge \ell^{'}}{\underset{r=1}{\sum}}
\frac{(\ell-1)!/\sqrt{\ell!}}{\sqrt{(\ell-r)!}\sqrt{(r-1)!}}\frac{(\ell^{'}-1)!/\sqrt{\ell^{'}!}}{\sqrt{(\ell^{'}-r)!}\sqrt{(r-1)!}}
\frac{\sqrt{(\ell+\ell^{'}-2r)!}}{\sqrt{(\ell-r)!}\sqrt{(\ell^{'}-r)!}}
\\\label{boundS2}=&
\overset{\ell\wedge \ell^{'}}{\underset{r=1}{\sum}}
\frac{r}{\ell^{'}}\frac{\sqrt{\ell!}}{\sqrt{(\ell-r)!}\sqrt{r!}}\frac{\sqrt{\ell^{'}!}}{\sqrt{(\ell^{'}-r)!}\sqrt{r!}}
\frac{\sqrt{(\ell+\ell^{'}-2r)!}}{\sqrt{(\ell-r)!}\sqrt{(\ell^{'}-r)!}}.
\end{align}
Because $r/\ell^{'}\leq1$ and
$
(\ell+\ell^{'}-2r)!\leq \sqrt{(2\ell-2r)!}  \sqrt{(2\ell^{'}-2r)!}
$
for all $\ell, \ell^{'}\in \mathbb{N}$ and $r\in\{1,2,\ldots,\ell\wedge \ell^{'}\}$,
(\ref{boundS2}) implies that
\begin{align}\notag
\Theta_{2}(\ell,\ell^{'}) \leq&
\overset{\ell\wedge \ell^{'}}{\underset{r=1}{\sum}}
\sqrt{\binom{\ell}{\ell-r}}\sqrt{\binom{\ell^{'}}{\ell^{'}-r}}
\left[\binom{2\ell-2r}{\ell-r}\right]^{\frac{1}{4}}
\left[\binom{2\ell^{'}-2r}{\ell^{'}-r}\right]^{\frac{1}{4}}
\\\notag\leq&
\left\{\overset{\ell\wedge \ell^{'}}{\underset{r=1}{\sum}}
\binom{\ell}{\ell-r}\sqrt{\binom{2\ell-2r}{\ell-r}}
\right\}^{\frac{1}{2}}
\left\{\overset{\ell\wedge \ell^{'}}{\underset{r=1}{\sum}}
\binom{\ell^{'}}{\ell^{'}-r}\sqrt{\binom{2\ell^{'}-2r}{\ell^{'}-r}}
\right\}^{\frac{1}{2}}
\\\label{bound:S2K}\leq&
\left\{\overset{\ell\wedge \ell^{'}}{\underset{r=1}{\sum}}
\binom{\ell}{\ell-r}2^{\ell-r}
\right\}^{\frac{1}{2}}
\left\{\overset{\ell\wedge \ell^{'}}{\underset{r=1}{\sum}}
\binom{\ell^{'}}{\ell^{'}-r}2^{\ell^{'}-r}
\right\}^{\frac{1}{2}}
\leq
3^{\frac{\ell}{2}}\ 3^{\frac{\ell^{'}}{2}},
\end{align}
where the second inequality follows from the Cauchy-Schwarz inequality
and the third inequality follows from
the fact $\binom{2\ell-2r}{\ell-r}\leq 4^{\ell-r}$.
By substituting (\ref{bound:S2K}) into (\ref{boundS2_proof_start}),
\begin{align}\label{proof:P2_final}
P_{2}(K) \leq
\underset{\begin{subarray}{c}
\ell,\ell^{'}\in\{2,4,\ldots,K\}\\
\ell\neq \ell^{'}
\end{subarray}}{\sum}|c_{\ell}c_{\ell^{'}}|\sqrt{\ell!}\sqrt{\ell^{'}!}\ 3^{\frac{\ell}{2}}\ 3^{\frac{\ell^{'}}{2}}.
\end{align}
By combining (\ref{proof:P1_final}) and (\ref{proof:P2_final}), we get
\begin{align}\label{mid_P1+P2}
P_{1}(K) + P_{2}(K) \leq \left(\underset{\ell\in\{2,4,\ldots,K\}}{\sum}|c_\ell| \sqrt{\ell!}\ 3^{\frac{\ell}{2}}\right)^{2}.
\end{align}
The inequality  (\ref{Prop:statement:stein:nu}) follows by substituting (\ref{mid_P1+P2}) into (\ref{initial_P1+P2}).

\section{Proof of Theorem \ref{theorem:convergence_rate:rational}}\label{sec:proof:theorem:convergence_rate:rational}
\begin{itemize}

\item For the case $A(r) = r^{\nu}$ with $\nu\in 2\mathbb{N}$,
(\ref{def:can:r_nu}) shows that $c_{A,\frac{\ell}{2}}=0$ if $\ell>\nu$.
It implies that the decomposition (\ref{thm:S1_chaosexpansion}) of
$S^{A}_{J}[j]X$ is only comprised of finite Wiener chaos.
Therefore, the truncation step (\ref{thm:S1_chaosexpansion_part2}) for the decomposition of $\mathbf{F}$
is unnecessary and
Theorem \ref{theorem:convergence_rate:rational}
follows from Proposition \ref{prop:multidim_stein_RHS} with $K=\nu$.

\item For the case $A(r) = r^{\nu}$ with $\nu\in(0,\infty)\setminus 2\mathbb{N}$, by (\ref{def:c_ell}) and (\ref{limit_binom2}),
for any $\varepsilon>0$,
there exists a constant $C_{1}(\nu,\varepsilon)$ such that
\begin{align}\label{proof:c_ell_rnu}
|c_{\ell}| \leq C_{1}(\nu,\varepsilon) 2^{\frac{\ell}{2}}\left(\frac{\ell}{2}!\right)(\ell!)^{-1}
\ell^{-\frac{\nu}{2}-1+\varepsilon}
\end{align}
for all $\ell\in 2\mathbb{N}$.
The inequality (\ref{proof:c_ell_rnu}) also holds for $A(r)=\ln(r)$ with the notation replacement: $\nu=0$ and $\varepsilon=0$.
By Stirling's formula (\ref{Stirling_formula}),
\begin{align*}
2^{\frac{\ell}{2}}\left(\frac{\ell}{2}!\right)(\ell!)^{-\frac{1}{2}}\leq 2 \ell^{\frac{1}{4}}.
\end{align*}
Hence, there exists a constant $C_{2}(\nu,\varepsilon)$ such that
\begin{align*}
\underset{\ell\in\{2,4,\ldots,K\}}{\sum}|c_\ell| \sqrt{\ell!}3^{\frac{\ell}{2}} \leq
C_{2}(\nu,\varepsilon)\overset{K/2}{\underset{\ell=1}{\sum}}\ell^{-\frac{\nu}{2}-\frac{3}{4}+\varepsilon}3^{\ell}.
\end{align*}
For the summation on the right hand side above, for each fixed $\nu\in [0,\infty)\setminus 2\mathbb{N}$, there exists a threshold $K_{\nu}>0$ such that
\begin{align*}
\overset{K/2}{\underset{\ell=1}{\sum}}\ell^{-\frac{\nu}{2}-\frac{3}{4}+\varepsilon}3^{\ell}\leq 3^{1+\frac{K}{2}}\left(K/2\right)^{-\frac{\nu}{2}-\frac{3}{4}+\varepsilon}
\end{align*}
for any even integer $K>K_{\nu}$.
Therefore, for any $\varepsilon>0$, there exists a constant $C_{3}(\nu,\varepsilon)>0$ such that
\begin{align}\label{stein_RHS_special}
\left(\underset{\ell\in\{2,4,\ldots,K\}}{\sum}|c_\ell| \sqrt{\ell!}3^{\frac{\ell}{2}} \right)^{2}\leq
C_{3}(\nu,\varepsilon)3^{K} K^{-\nu-\frac{3}{2}+\varepsilon}
\end{align}
for any even integer $K>K_{\nu}$.
By substituting (\ref{stein_RHS_special}) into (\ref{Prop:statement:stein:nu}),
\begin{align}\label{proof:middlepart_triangle}
|\mathbb{E}[h(\mathbf{F}_{\leq K})]-\mathbb{E}[h(\mathbf{N}_{\mathbf{C}_{J,K}})]|\leq  CC_{3}(\nu,\varepsilon) 2^{-\frac{J}{2}}3^{K} K^{-\nu-\frac{3}{2}+\varepsilon}
\end{align}
for any even integer $K>K_{\nu}$.
Finally, by applying the results of Lemma \ref{lemma:gp_bound} and  (\ref{proof:middlepart_triangle})
to (\ref{triangle_inequality}),
we obtain that
\begin{align*}
|\mathbb{E}[h(\mathbf{F})]-\mathbb{E}[h(\mathbf{N})]|\leq C_{4}(\nu,\varepsilon)\left( K^{-\frac{\nu}{2}-\frac{1}{4}+\varepsilon}+
2^{-\frac{J}{2}}3^{K} K^{-\nu-\frac{3}{2}+\varepsilon}\right)
\end{align*}
for a certain constant $C_{4}(\nu,\varepsilon)$ independent of $J$ and $K$ when $K\geq K_{\nu}$.
The asymptotic behavior (\ref{sWdistance:F_N_O}) follows by choosing $K=2\lfloor\frac{J}{4}\log_{3}2\rfloor$.

\end{itemize}

Finally, because both sides of (\ref{limit_cov}) can be computed by the Wiener chaos expansions (\ref{thm:representationU}) and (\ref{thm:S1_chaosexpansion})
of $U^{A}[j]X$ and $S^{A}_{J}[j]X$,
the verification of (\ref{limit_cov}) is omitted.

\section{Proof of Corollary \ref{corollary:kol}}\label{sec:proof:kol}
The idea of proof is originated from the work \cite[Proposition 2.6]{gaunt2023bounding}, in which more complicated cases were considered.
However, we only need part of it. For the convenience of readers, we will sketch its proof as follows.
\begin{figure}[h]
\centering
\subfigure[][$h$ and $h_{z,\delta}$]
{\includegraphics[scale=0.5]{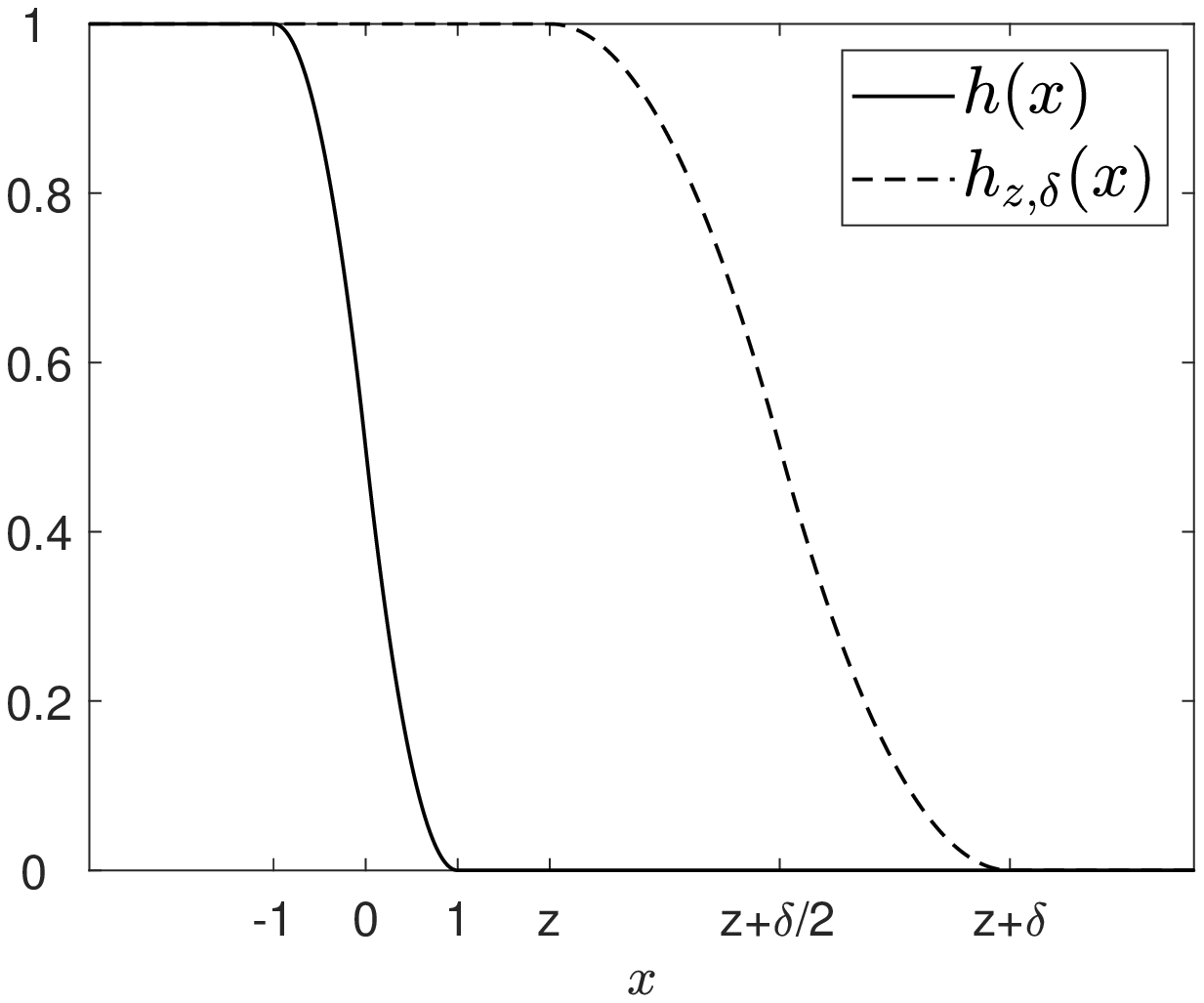}} % assumption5_case1_weak2.eps
\subfigure[][$h_{(-1,-1),2}$]
{\includegraphics[scale=0.5]{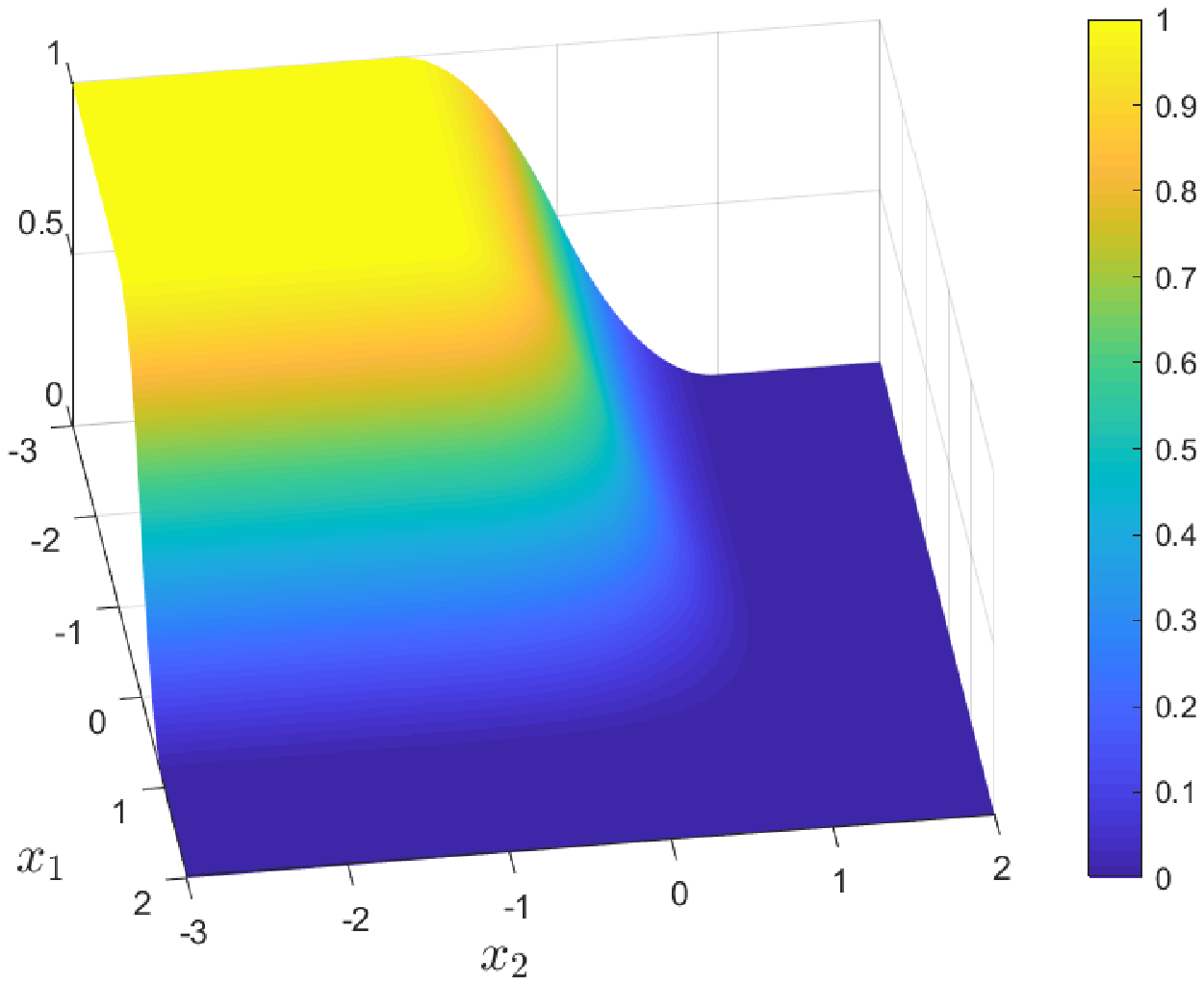}} % assumption5_case2_weak2.eps
\caption{Approximation of indicator functions}
\label{fig:spline}
\end{figure}
First of all, we define
\begin{align*}
h(x) = \left\{
\begin{array}{ll}
1 & \textup{for}\ x\leq -1,\\
1-\frac{1}{2}(1+x)^{2} & \textup{for}\ x\in(-1,0],\\
\frac{1}{2}(1-x)^{2} & \textup{for}\ x\in(0,1],\\
0 & \textup{for}\ x>1,\end{array}\right.
\end{align*}
and $h_{z,\delta}(x) = h\left(\frac{2}{\delta}\left(x-\left(z+\frac{\delta}{2}\right)\right)\right)$
for any $z\in \mathbb{R}$ and $\delta>0$.
See Figure \ref{fig:spline}(a) for the graphs of $h$ and $h_{z,\delta}$.
For any $d\in \mathbb{N}$ and $\mathbf{z}=(z_{1},\ldots,z_{d})\in \mathbb{R}^{d}$,
let
\begin{align*}
h_{\mathbf{z},\delta}(\mathbf{x}) = \overset{d}{\underset{k=1}{\prod}}h_{z_{k},\delta}(x_{k}),\ \mathbf{x}=(x_{1},\ldots,x_{d})\in \mathbb{R}^{d}.
\end{align*}
See Figure \ref{fig:spline}(b) for the graph of  $h_{\mathbf{z},\delta}$ with $\mathbf{z}=(-1,-1)$ and $\delta=2$.
For any fixed $\mathbf{z}\in \mathbb{R}^{d}$,  the indicator function
$1_{\mathbf{x}\leq \mathbf{z}}$ is bounded above by $h_{\mathbf{z},\delta}(\mathbf{x})$,
where  $\mathbf{x}\leq \mathbf{z}$ means that $x_{k}\leq z_{k}$ for all $k=1,\ldots,d$.
Hence,
\begin{align}\notag
\mathbb{P}\left(\mathbf{F}\leq \mathbf{z}\right)-\mathbb{P}\left(\mathbf{N}\leq \mathbf{z}\right)
\leq& \mathbb{E}\left[h_{\mathbf{z},\delta}(\mathbf{F})\right]-\mathbb{E}\left[h_{\mathbf{z},\delta}(\mathbf{N})\right]
+\mathbb{E}\left[h_{\mathbf{z},\delta}(\mathbf{N})\right]-\mathbb{P}\left(\mathbf{N}\leq \mathbf{z}\right)
\\\notag\leq&
\mathbb{E}\left[h_{\mathbf{z},\delta}(\mathbf{F})\right]-\mathbb{E}\left[h_{\mathbf{z},\delta}(\mathbf{N})\right]
+\mathbb{P}(\mathbf{N}\leq \mathbf{z}+\delta )-\mathbb{P}\left(\mathbf{N}\leq \mathbf{z}\right),
\end{align}
where $\mathbf{z}+\delta$ means that $\mathbf{z}+\delta \mathbf{1}_{1\times d}$.
Because $\|h_{z,\delta}\|_{\infty}=1$, $\|h_{z,\delta}^{'}\|_{\infty}=\frac{2}{\delta}$ and
$\|h_{z,\delta}^{'}\|_{\textup{Lip}}=(\frac{2}{\delta})^{2}$, we have
$\left[1+\left(\frac{2}{\delta}\right)^{2}\right]^{-1}h_{\mathbf{z},\delta}
\in \mathcal{H}_{2}$ for any $\mathbf{z}\in \mathbb{R}^{d}$ and $\delta>0$, which implies that
\begin{align}\label{Nazarov_inequalityv1}
\left|\mathbb{E}\left[h_{\mathbf{z},\delta}(\mathbf{F})\right]-\mathbb{E}\left[h_{\mathbf{z},\delta}(\mathbf{N})\right]\right|
\leq \left[1+\left(\frac{2}{\delta}\right)^{2}\right]d_{\mathcal{H}_{2}}\left(\mathbf{F},\mathbf{N}\right).
\end{align}
On the other hand, from Nazarov's anti-concentration inequality \cite[Lemma A.1]{chernozhukov2017central}, we have
\begin{align}\label{Nazarov_inequalityv2}
\mathbb{P}(\mathbf{N}\leq \mathbf{z}+\delta)-\mathbb{P}\left(\mathbf{N}\leq \mathbf{z}\right)
\leq \frac{\delta}{\sqrt{\underset{1\leq m\leq d}{\min}\mathbb{E}[F_{m}^{2}]}}\left(\sqrt{2\ln d}+2\right)
\end{align}
for any $\mathbf{z}\in \mathbb{R}^{d}$ and $\delta>0$.
%where $\sigma = \underset{1\leq k\leq d}{\min}\sqrt{\mathbb{E}[N_{k}^{2}]}=\sqrt{\mathbb{E}[F_{1}^{2}]}$.
Hence,
\begin{align}\label{kol_right_v1}
\mathbb{P}\left(\mathbf{F}\leq \mathbf{z}\right)-\mathbb{P}\left(\mathbf{N}\leq \mathbf{z}\right)
\leq
\left[1+\left(\frac{2}{\delta}\right)^{2}\right]d_{\mathcal{H}_{2}}\left(\mathbf{F},\mathbf{N}\right)+\frac{\delta}{\sqrt{\underset{1\leq m\leq d}{\min}\mathbb{E}[F_{m}^{2}]}}\left(\sqrt{2\ln d}+2\right).
\end{align}
By taking $\delta=2\left[\sqrt{\underset{1\leq m\leq d}{\min}\mathbb{E}[F_{m}^{2}]}(\sqrt{2\ln d}+2)^{-1}d_{\mathcal{H}_{2}}\left(\mathbf{F},\mathbf{N}\right)\right]^{\frac{1}{3}}$, (\ref{kol_right_v1}) leads to the upper bound
\begin{align}\label{upperbound_kol}
\mathbb{P}\left(\mathbf{F}\leq \mathbf{z}\right)-\mathbb{P}\left(\mathbf{N}\leq \mathbf{z}\right)
\leq
3\left(\frac{\sqrt{2\ln d}+2}{\sqrt{\underset{1\leq m\leq d}{\min}\mathbb{E}[F_{m}^{2}]}}\right)^{\frac{2}{3}}\left(d_{\mathcal{H}_{2}}\left(\mathbf{F},\mathbf{N}\right)\right)^{\frac{1}{3}}
+
d_{\mathcal{H}_{2}}\left(\mathbf{F},\mathbf{N}\right).
\end{align}
For finding a lower bound for $\mathbb{P}\left(\mathbf{F}\leq \mathbf{z}\right)-\mathbb{P}\left(\mathbf{N}\leq \mathbf{z}\right)$, we first observe that
\begin{align}\notag
\mathbb{P}\left(\mathbf{N}\leq \mathbf{z}\right)-\mathbb{P}\left(\mathbf{F}\leq \mathbf{z}\right)
\leq& \mathbb{P}\left(\mathbf{N}\leq \mathbf{z}\right)
-
\mathbb{E}\left[h_{\mathbf{z}-\delta,\delta}(\mathbf{N})\right]
+\mathbb{E}\left[h_{\mathbf{z}-\delta,\delta}(\mathbf{N})\right]
-\mathbb{E}\left[h_{\mathbf{z}-\delta,\delta}(\mathbf{F})\right]
\\\label{lower_bound_kol}\leq&
\mathbb{P}\left(\mathbf{N}\leq \mathbf{z}\right)-\mathbb{P}\left(\mathbf{N}\leq \mathbf{z}-\delta\right)
+\mathbb{E}\left[h_{\mathbf{z}-\delta,\delta}(\mathbf{N})\right]
-\mathbb{E}\left[h_{\mathbf{z}-\delta,\delta}(\mathbf{F})\right]
\end{align}
for any $\mathbf{z}\in \mathbb{R}^{d}$ and $\delta>0$.
By applying (\ref{Nazarov_inequalityv1}) and (\ref{Nazarov_inequalityv2})
to (\ref{lower_bound_kol}), we obtain a lower bound for $\mathbb{P}\left(\mathbf{F}\leq \mathbf{z}\right)-\mathbb{P}\left(\mathbf{N}\leq \mathbf{z}\right)$,
which is the negative of the right hand side of (\ref{upperbound_kol}).
Therefore,
\begin{align*}%\label{twoside_bound_kol}
d_{\textup{Kol}}\left(\mathbf{F},\mathbf{N}\right)
\leq
3\left(\frac{\sqrt{2\ln d}+2}{\sqrt{\underset{1\leq m\leq d}{\min}\mathbb{E}[F_{m}^{2}]}}\right)^{\frac{2}{3}}\left(d_{\mathcal{H}_{2}}\left(\mathbf{F},\mathbf{N}\right)\right)^{\frac{1}{3}}
+
d_{\mathcal{H}_{2}}\left(\mathbf{F},\mathbf{N}\right).
\end{align*}
Because Theorem \ref{theorem:convergence_rate:rational} shows that
$\underset{J\rightarrow\infty}{\lim}d_{\mathcal{H}_{2}}\left(\mathbf{F},\mathbf{N}\right)=0$ and
$\underset{J\rightarrow\infty}{\lim}\mathbb{E}[F_{m}^{2}]=\kappa_{m,m}\|\widehat{\phi}\|_{2}^{2}$,
we obtain
\begin{align}\label{summary_kol}
d_{\textup{Kol}}\left(\mathbf{F},\mathbf{N}\right)\leq \mathcal{O}\left(\left(d_{\mathcal{H}_{2}}\left(\mathbf{F},\mathbf{N}\right)\right)^{\frac{1}{3}}\right)
\end{align}
when $J\rightarrow\infty$. Corollary \ref{corollary:kol} follows by combining (\ref{sWdistance:F_N_O}) and (\ref{summary_kol}).

\end{appendix}

\end{document}